\documentclass[12pt, a4paper]{article}
\usepackage{graphicx, theorem, color} \definecolor{webblue}{rgb}{0,0,0.6}
\usepackage[colorlinks=true, linkcolor=webblue, citecolor=webblue,
 urlcolor=webblue, pdfstartview=FitBH,
 pdftitle={The Thurston Algorithm for quadratic matings},
 pdfauthor={Wolf Jung},
 bookmarksopen=true, bookmarksnumbered=true, pdfpagemode=UseNone]{hyperref}
\makeatletter 
\long\def\@makecaption#1#2{%
  \vskip\abovecaptionskip
  {\small\textbf{#1:} #2}\par
  \vskip\belowcaptionskip}
\makeatother
\theoremstyle{break}
\theorembodyfont{\itshape}  \newtheorem{thm}{Theorem}[section]
\newtheorem{dfn}[thm]{Definition} \newtheorem{prop}[thm]{Proposition}
 
\newtheorem{init}[thm]{Initialization}
\theorembodyfont{\slshape}  \newtheorem{rmk}[thm]{Remark}  
\newtheorem{xmp}[thm]{Example}   
\newcommand{\be}{\begin{equation}} \newcommand{\ee}{\end{equation}}
\newcommand{\ban}{\begin{eqnarray}} \newcommand{\ean}{\end{eqnarray}}
\renewcommand{\hat}{\widehat} \renewcommand{\tilde}{\widetilde}
\newcommand{\mybox}{\hspace*{\fill}\rule{2mm}{2mm}}
\DeclareSymbolFont{AMSb}{U}{msb}{m}{n}
\DeclareSymbolFontAlphabet{\mathbb}{AMSb}
\newcommand{\C}{\mathbb{C}}
\newcommand{\Q}{\mathbb{Q}}  \newcommand{\R}{\mathbb{R}}
\newcommand{\Z}{\mathbb{Z}}
\newcommand{\N}{\mathbb{N}}  \newcommand{\disk}{\mathbb{D}}
\newcommand{\e}[1]{\displaystyle {\rm e}^{\displaystyle #1}}
\let\inodot\i  \renewcommand{\i}{\mathrm{i}}
\newcommand\ddoti{\"\inodot}  
\newcommand{\eps}{\varepsilon}  
\renewcommand{\phi}{\varphi}
\renewcommand{\r}{\mathcal{R}}
  \newcommand{\K}{\mathcal{K}}
\newcommand{\M}{\mathcal{M}}  \newcommand{\sM}{{\scriptscriptstyle M}}
\renewcommand{\O}{\mathcal{O}}  \renewcommand{\S}{\mathcal{S}}
\newcommand{\T}{\mathcal{T}}  \newcommand{\W}{\mathcal{W}}
\newcommand{\sC}{{\scriptscriptstyle C}}
\newcommand{\sT}{{\scriptscriptstyle T}}
\newcommand{\sWP}{{\scriptscriptstyle WP}}
\newcommand{\eqth}{\sim} 
\newcommand{\eqr}{\sim} 
\newcommand{\eqg}{\cong} 
\newcommand{\fmate}{\sqcup} 
\newcommand{\tmate}{\coprod} 
\newcommand{\2}{$(2,\,2,\,2,\,2)$}

\date{}
\author{Wolf Jung\\
{\small Gesamtschule Brand, 52078 Aachen, Germany,}\\
{\small and Jacobs University, 28759 Bremen, Germany.}\\
{\small E-mail: \href{mailto:jung@mndynamics.com}{jung@mndynamics.com}}\\[7mm]
\textit{Dedicated to the memory of Tan Lei}}
\title{The Thurston Algorithm for quadratic matings}
\begin{document}
\maketitle
\begin{abstract}\noindent
Mating is an operation to construct a rational map $f$ from two polynomials,
which are not in conjugate limbs of the Mandelbrot set. When the Thurston
Algorithm for the unmodified formal mating is iterated in the case of
postcritical identifications, it will diverge to the boundary of Teichm\"uller
space, because marked points collide. Here it is shown that the colliding
points converge to postcritical points of $f$, and the associated sequence of
rational maps converges to $f$ as well, unless the orbifold of $f$ is of type
\2. So to compute $f$, it is not necessary to encode the topology of
postcritical ray-equivalence classes for the modified mating, but it is
enough to implement the pullback map for the formal mating. The proof combines
the Selinger extension to augmented Teichm\"uller space with local estimates.
\par\noindent
Moreover, the Thurston Algorithm is implemented by pulling back a path in
moduli space. This approach is due to Bartholdi--Nekrashevych in relation
to one-dimensional moduli space maps, and to Buff--Ch\'eritat for slow mating.
Here it is shown that slow mating is equivalent to the Thurston Algorithm for
the formal mating. An initialization of the path is obtained for
repelling-preperiodic captures as well, which provide an alternative
construction of matings. 
\end{abstract}

\section{Introduction} \label{1}
\begin{figure}[h!t!b!]
\unitlength 0.15mm
\begin{picture}(999, 280)
\put(10, 0){\includegraphics[width=0.28\textwidth]{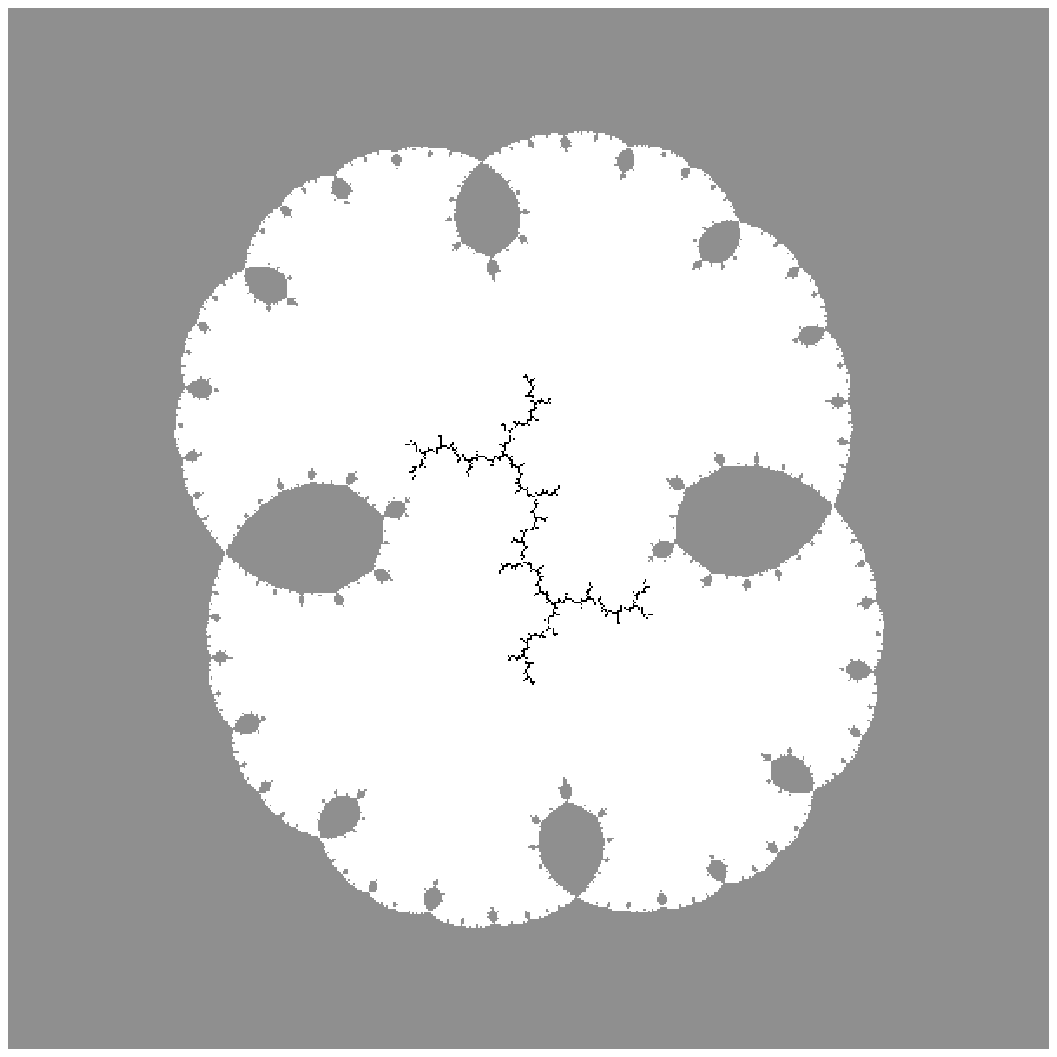}}
\put(360, 0){\includegraphics[width=0.28\textwidth]{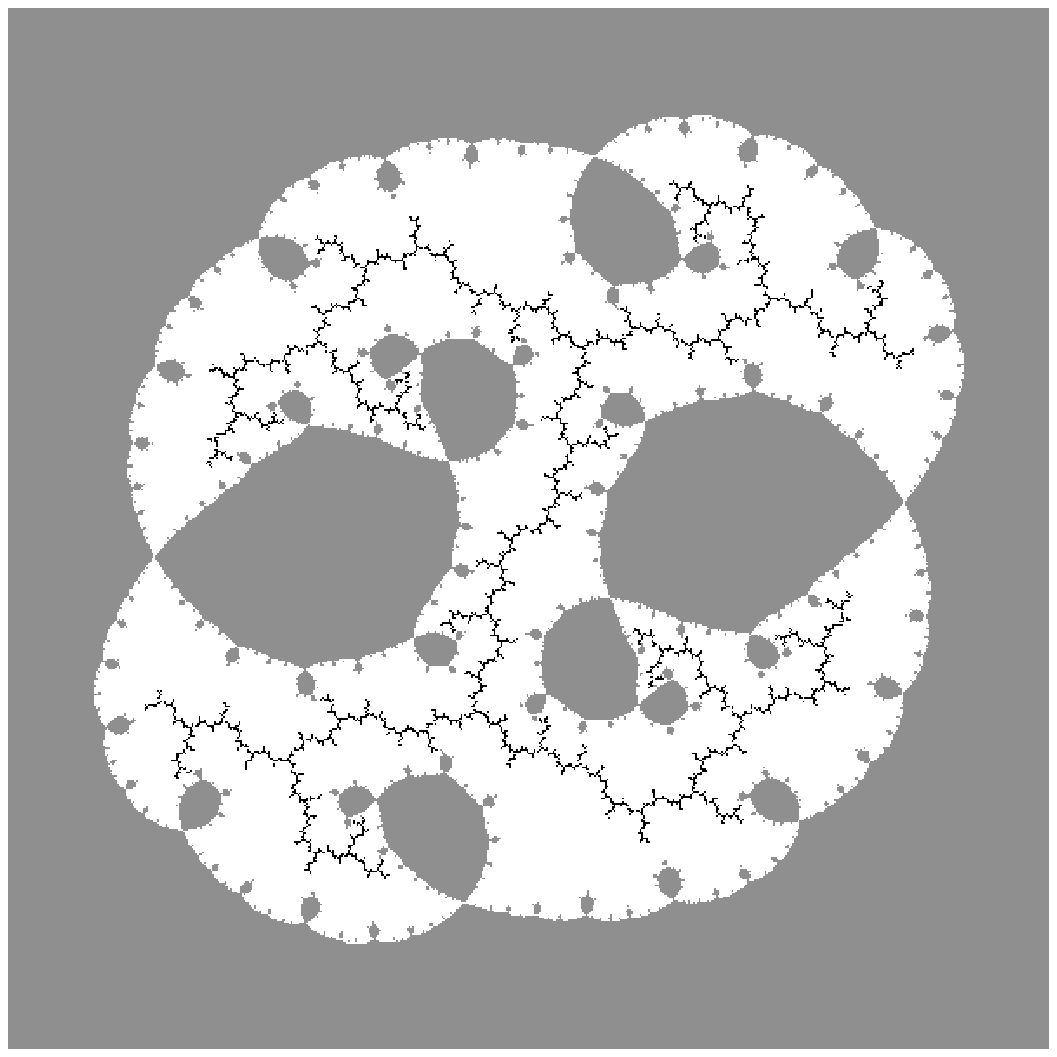}}
\put(710, 0){\includegraphics[width=0.28\textwidth]{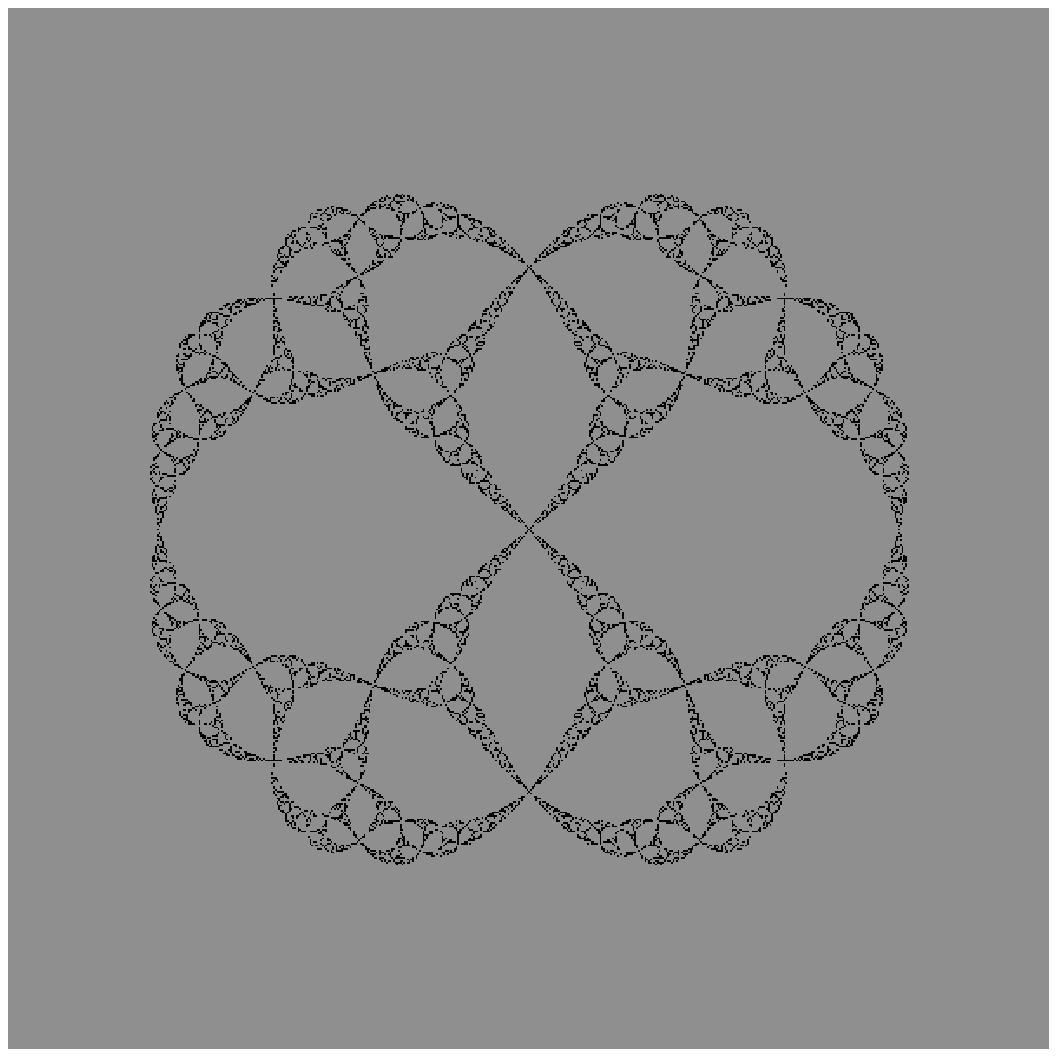}}
\thinlines
\end{picture} \caption[]{\label{Fib}
Various stages of slow mating, illustrated by moving images
$\psi_t(\phi_0(\K_p))$ and $\psi_t(\phi_\infty(\K_q))$ of Julia sets.
Here $P(z)=z^2+\i$ is a Misiurewicz polynomial and $Q(z)=z^2-1$ is the
Basilica polynomial, which has an attracting $2$-cycle. The formulas for
pulling back marked points and rational maps are discussed in
(\ref{eqmibnop}) and in Example~\ref{Xsmib} as well.}
\end{figure}

A postcritically finite quadratic polynomial $f_c(z)=z^2+c$ may be periodic of
satellite type, periodic of primitive type, or critically  preperiodic
(Misiurewicz type). Examples are given by the Basilica $Q(z)=z^2-1$, the
Kokopelli, and by $P(z)=z^2+\i$ in Figures~\ref{Fib} and~\ref{Fkk}. Quadratic
rational maps have two critical orbits and form a two-parameter family. The
dynamics and topology of certain rational maps is understood in terms of one or
two polynomials \cite{rees1}. The operation of mating was introduced by
Douady--Hubbard \cite{sdh}: a rational map $f$ may be described by gluing the
Julia sets of $P$ and $Q$, such that points with conjugate external angles are
identified. According to Rees--Shishikura--Tan \cite{rst, rs}, this
construction works when $P$ and $Q$ are not in conjugate limbs of the
Mandelbrot set: first define the formal mating, where the two Julia sets are in
separate half-spheres. The Thurston Theorem \cite{dh, book2h} shows that there
is an equivalent rational map $f$. Then the topological mating is given by
collapsing all ray-equivalence classes of the formal mating, and it is
conjugate to the geometric mating $f$. Actually, an intermediate step is
required when postcritical points are identified in the mating: then the
formal mating will be obstructed, and an unobstructed essential mating is
constructed by collapsing a finite number of ray-equivalence classes.

The Thurston Algorithm is based on an iteration in Teichm\"uller space, which
consists of isotopy classes of homeomorphisms, These may be represented by
spiders, medusas, or triangulations. Bartholdi--Nekrashevych \cite{bn}
and Buff-Ch\'eritat \cite{cts3} employ a path in moduli space instead.
Using this ``slow'' approach, the algorithm shall be faster,  easier to
implement, and more stable. Explicit initializations are discussed in
Sections~\ref{2p}, \ref{5}, and~\ref{6c}. The slow mating algorithm is
related to equipotential gluing in \cite{emate}. Figure~\ref{Fib} shows
a few snapshots of this process.

In this example, the Thurston Algorithm does not work directly, because the
postcritical $2$-cycle of $P$ needs to be identified with a fixed point of $Q$:
these are connected by external rays, and the ray-equivalence class is
surrounded by a removable Thurston obstruction. The classical approach is to
construct an essential mating, where certain ray-equivalence classes are
collapsed by definition, and to employ the Thurston Algorithm for the modified
map. An alternative approach is suggested here: the divergence of the Thurston
Algorithm has been described by Nikita Selinger \cite{ext, char}
in terms of the augmented Teichm\"uller space. Applying his characterization to
the Thurston Algorithm of the unmodified formal mating, it is shown that marked
points come together automatically in the expected way, and the rational maps
converge to the geometric mating, at least if the orbifold is not of type \2.
The same argument gives convergence of slow mating and
equipotential gluing as well, where no modification is appropriate.
Thus it is possible to obtain matings numerically
without encoding the topology of ray-equivalence classes. In a few more
applications, additional obstructions are created and used to prove convergence
properties \cite{pmate, amate}. Here obstructions do not appear as a
potential problem, but they are turned into an ally: a powerful tool to show
convergence.

The classical Thurston Theorems are discussed in Section~\ref{2}, together
with the implementation by a path in moduli space. See Section~\ref{3} for
examples of canonical obstructions and stabilization of noded Riemann surfaces,
including a relation between core entropy and matings of conjugate polynomials.
A general convergence result is obtained for bicritical maps in a suitable
normalization. The various concepts of mating are developed systematically in
Section~\ref{4}. Convergence of mating and the implementation of slow mating
and of captures is discussed in Sections~\ref{4c}, \ref{5} and~\ref{6c},
respectively. Various algorithms are compared briefly in Section~\ref{7}.

\textbf{Acknowledgment:}
Several colleagues have contributed to this work by inspiring discussions
and helpful suggestions. I wish to thank in particular
Laurent Bartholdi, Walter Bergweiler, Xavier Buff, Arnaud Ch\'eritat,
Dzmitry Dudko, Adam Epstein, Sarah Koch, Michael Mertens,
Carsten Petersen. Mary Rees, Pascale Roesch, Dierk Schleicher,
Nikita Selinger, Tan Lei, Dylan P.~Thurston, and Vladlen Timorin.
And I am grateful to the mathematics department of
Toulouse University Paul Sabatier for their hospitality.

\textsl{When I started to learn about complex dynamics some twenty years ago,
the most impressive phenomenon was the self-similarity of the Mandelbrot set at
Misiurewicz points, and Tan Lei's name is firmly attached to this. She has
worked on parabolic maps, quasi-conformal deformations, and on vector fields
as well. A recurring theme in her work is the Thurston characterization of
rational maps, in particular its application to matings. I remember joyful
discussions in Holb{\ae}k 2007, and in recent years we have had a few
conversations about topological entropy. Tan Lei passed away in April 2016.
This paper is dedicated to her memory.}

\section{The Thurston Algorithm} \label{2}
The Thurston Theorem~\ref{TT1} gives a combinatorial
characterization of branched covers equivalent to rational maps, which is used
to describe and to define rational maps, and to construct them numerically:
the related Thurston Algorithm provides a convergent sequence of rational maps.
An underlying iteration in Teichm\"uller space $\T$ is needed both to define a
unique pullback, and to have analytic tools providing global convergence to
a unique fixed point in $\T$. This fails in the presence of Thurston
obstructions: then certain annuli get big, curves get short, marked points
collide. This process is understood by extending the pullback map to
augmented Teichm\"uller space $\hat\T$; see Section~\ref{3}. The present
section provides an introduction to the classical theory by
William P.~Thurston, and Section~\ref{2p} discusses the implementation of the
Thurston Algorithm in terms of a path in moduli space $\M$.

\subsection{Hyperbolic geometry and Teichm\"uller spaces} \label{2g}
A hyperbolic Riemann surface of finite type and genus $0$ is isomorphic to
the Riemann sphere $\hat\C=\C\cup\{\infty\}$ with $n\ge3$ punctures. Although
the manifold extends analytically to a puncture or \textbf{marked point}, the
hyperbolic metric is infinite there.
We shall deal with homotopy classes of simple closed curves and the
hyperbolic length of geodesics:
\begin{itemize}
\item A simple closed curve in the complement of the marked points is
\textbf{essential}, if each disk in the complement of the curve contains at
least two marked points;
\item \textbf{peripheral}, if one component contains only one marked point; and
\item \textbf{trivial} or null-homotopic, if one component contains no marked
point.
\end{itemize}
Note that some authors say non-peripheral instead of essential or inessential
instead of peripheral. The following properties of hyperbolic geodesics
are fundamental:

\begin{prop}[Hyperbolic geodesics]\label{Pcoll}
Consider the hyperbolic metric on $\hat\C$ with $n\ge3$ punctures:

$1$. For any essential simple closed curve there is a unique geodesic
homotopic to it.

$2$. A simple closed geodesic $\gamma$ has a \textbf{collar} neighborhood, an
embedded annulus of definite width. Collars around disjoint geodesics are
disjoint, and a geodesic crossing the collar of $\gamma$ has an explicit lower
bound on its length, which goes to $\infty$ when $l(\gamma)\to0$. In
particular, all sufficiently short geodesics are disjoint.

$4$. Any annulus around $\gamma$ has modulus bounded above by
$\pi/l(\gamma)$. The collar has modulus bounded below by $\pi/l(\gamma)-1$.

$4$. For a sequence of surfaces in a suitable normalization, two marked points
collide with respect to the spherical metric, if and only if a hyperbolic
geodesic separating them from two other marked points has length going to $0$.
\end{prop}

References for the \textbf{proof:} See \cite{book1h} for item~1 and
\cite{dh, book1h, teich} for items~2 and~3. Item~4 is a standard estimate for
extremal annuli. \mybox

For an implicit $n\ge3$, \textbf{moduli space} $\M$ is the space of Riemann
spheres with marked points up to M\"obius maps or normalization of three
points; in our case of genus $0$, it has an explicit description as a subset of
$\hat\C^{n-3}$ from the positions of marked points. Now
\textbf{Teichm\"uller space} $\T$ is the universal cover of $\M$. It can be
described by isotopy classes of orientation-preserving homeomorphisms
$\psi:\hat\C\to\hat\C$ between spheres with marked points; here the left sphere
is a topological sphere fixed for reference. Although it has no complex
structure, let us write $\hat\C$ instead of $S^2$ nevertheless: this allows to
use explicit coordinates and formulas from $\C$. The \textbf{projection}
$\pi:\T\to\M$ gives the universal cover, and the
\textbf{pure mapping class group} $G$ is the group of deck transformations:
$[h]\in G$ is an isotopy class of homeomorphisms of the topological sphere
fixing the marked points, which acts on $\T$ by
$[h]\cdot[\psi]=[\psi\circ h^{-1}]$\,; $G$ is generated by Dehn twists
\cite{book1h}.

There are various metrics on $\T$, such that $G$ acts by isometries and the
metrics project to $\M$. Actually the definition as Finsler metrics is lifted
from $\M$ to $\T$ locally. The dual tangent space is given by integrable
quadratic differentials, and by using different norms there, the infinitesimal
metrics $\|\mathrm{d}\tau\|_\sT$ and $\|\mathrm{d}\tau\|_\sWP$ are obtained;
integration along shortest curves defines the \textbf{Teichm\"uller metric}
$d_\sT$ and the \textbf{Weil--Petersson metric} $d_\sWP$
\cite{book1h, wwp, mwp}. The Teichm\"uller metric is given equivalently by
$d_\sT([\psi_1],\,[\psi_2])=\inf1/2\log K(\psi)$, where $\psi$ is isotopic to
$\psi_1\circ\psi_2^{-1}$ and $K$ denotes the quasi-conformal dilatation.

\begin{prop}[Basic properties of Teichm\"uller space]\label{Pbptms}
$1$. $\T$ and $\M$ are analytic manifolds and the universal cover
$\pi:\T\to\M$ is analytic. It is a local isometry for both metrics,
$d_\sT$ and $d_\sWP$\,, and $G$ acts by isometries. 

$2$. Both metrics generate the same topology on $\T$. Now $\T$ is complete and
not compact with respect to $d_\sT$\,, incomplete with respect to $d_\sWP$\,.

$3$. For an essential curve $\gamma$ in the topological sphere and $\tau\in\T$,
denote by $l(\gamma,\,\tau)$ the \textbf{length} of the geodesic in the
Riemann surface $\pi(\tau)$, that is homotopic to $\psi(\gamma)$ for
$\psi\in\tau$. This length is continuous on $\T$ with
$|\log l(\gamma,\,\tau')-\log l(\gamma,\,\tau)|\le2d_\sT(\tau',\,\tau)$.

$4$. There are finitely many essential curves $\gamma_i$\,, such that
the collection of length functions $l(\gamma_i\,,\,\tau')$ determines $\tau$
uniquely.

$5$. There is a relative estimate
$\|\mathrm{d}\tau\|_\sWP\le C_*\|\mathrm{d}\tau\|_\sT$
with $C_*=C_*(\T)$.

$6$. For $R>0$ there is $D_*=D_*(R,\,\T)>0$ such that all $\tau$ with shortest
geodesic length $l_*(\tau)\ge R$ satisfy: all $\tau'$ with
$d_\sWP(\tau',\,\tau)\le D_*$ have $d_\sT(\tau',\,\tau)\le1/4$.

$7$. A closed subset of $\M$ is compact, if and only if the length of all
simple closed geodesics is bounded uniformly below.
\end{prop}

References for the \textbf{proof:}
See \cite{book1h, wwp, mwp} for items~1--4 and \cite{mcme} for item~5.

6. According to Lemma~3.22 in \cite{mwp} 
we have the relative estimate
\be\label{eqrelesttwp}
\|\mathrm{d}\tau\|_\sT\le\frac{C}{l_*(\tau)}\|\mathrm{d}\tau\|_\sWP\ ,\ee
where $C$ depends only on the Teichm\"uller space $\T$ and $l_*$ denotes the
length of the shortest hyperbolic geodesic on the Riemann surface $\pi(\tau)$.
So we must choose $D_*$ small enough to have a lower bound for $l_*(\cdot)$ on
the $WP$-geodesic from $\tau'$ to $\tau$\,. Note that it is not sufficient to
have a lower bound at $\tau'$ and $\tau$ only; cf.~Remark~\ref{Raug}.
Suppose $\tau'$ has $l_*(\tau')<R/2$, then for $\eps>0$ there is a subarc
$[\tilde\tau',\,\tilde\tau]$ of the $WP$-geodesic $[\tau',\,\tau]$ and a simple
closed curve $\gamma$, such that $l(\gamma,\,\cdot)$ increases from $R/2+\eps$
to $R$ along this subarc while $l_*(\cdot)\ge R/2$. Now
$|\mathrm{d}l|\le2l\|\mathrm{d}\tau\|_\sT$ and (\ref{eqrelesttwp}) give
\be
d_\sWP(\tau',\,\tau)>d_\sWP(\tilde\tau',\,\tilde\tau)
\ge\frac{1}{2C}\int_{[\tilde\tau',\,\tilde\tau]}
\frac{l_*(\cdot)}{l(\gamma,\,\cdot)}\mathrm{d}l(\gamma,\,\cdot)
\ge\frac{1}{4C}\int_{R/2+\eps}^R\mathrm{d}l \ ,\ee
so if $D_*=\frac{R}{8C}$ and $d_\sWP(\tau',\,\tau)\le D_*$ then
$l_*(\cdot)\ge R/2$ on the $WP$-geodesic $[\tau',\,\tau]$ and
(\ref{eqrelesttwp}) gives $d_\sT(\tau',\,\tau)\le\frac{2C}{R}\cdot D_*=1/4$.

7. This is the Mumford compactness theorem, whose proof is simplified in genus
$0$: the length is bounded below on any bounded ball. Every sequence in $\M$
has a convergent subsequence in $\hat\C^{n-3}$, but if the length of geodesics
may go to $0$, marked points collide according to Proposition~\ref{Pcoll}.4 and
the limit does not belong to $\M$. Note that every compact subset of $\T$ has
length bounded below as well, but the converse is wrong: for general $g\in G$
the sequence $\tau_k=g^k\cdot\tau_0$ does not accumulate in $\T$, but all
$l(\gamma,\,\tau_k)$ are independent of $k$. \mybox

\subsection{Thurston maps and pullback map} \label{2f}
A postcritically finite rational map $f$ is not characterized uniquely by
its ramification portrait. Except for flexible Latt\`es maps, the
additional topological information can be given combinatorially. For
polynomials, Hubbard trees and external angles provide an explicit description.
For rational maps, the combinatorial object is an equivalence class of
Thurston maps:
\begin{itemize}
\item A \textbf{Thurston map} $g:\hat\C\to\hat\C$ is an orientation-preserving
branched cover of degree $d\ge2$ with finite postcritical set $P$ and marked
set $Z\supset P$. Here $P$ contains all forward iterates of critical points
and $Z$ may contain additional critical, preperiodic, and periodic points,
such that $g(Z)\subset Z$.
\item Two Thurston maps $f,\,g$ are Thurston equivalent or
\textbf{combinatorially equivalent}, if there are homeomorphisms
$\psi_0\,,\,\psi_1$ with $\psi_0\circ g=f\circ \psi_1$\,, $\psi_0=\psi_1$ on
$Z_g$\,, $\psi_i(Z_g)=Z_f$\,, and $\psi_1$ is isotopic to $\psi_0$ relative to
$Z_g$\,. So $g$ is isotopic to $\psi_1^{-1}\circ f\circ\psi_1$\,.
\item A pullback map is associated with each Thurston map $g$ as follows: for
any homeomorphism $\psi$ there is a rational map $f$ and another homeomorphism
$\psi'$ with $\psi\circ g=f\circ \psi'$; see \cite{dh, book2h, teich}. The
complex structure defined by $\psi$ is pulled back with $g$ and integrated with
$\psi'$. These functions
are unique up to M\"obius maps, or unique after normalizing three marked points
to $\infty,\,0,\,1$. It turns out that the isotopy class of $\psi$ determines
the isotopy class of $\psi'$, so an analytic \textbf{pullback map}
$\sigma_g:\T\to\T$ is defined by $\sigma_g([\psi])=[\psi']$.
\item The fixed points of $\sigma_g$ in $\T$ correspond to M\"obius
conjugacy classes of rational maps. Under suitable conditions, the pullback map
will be strictly contracting, and the \textbf{Thurston Algorithm} converges in
addition: define $\psi_n\circ g=f_n\circ \psi_{n+1}$ recursively, then
$f_n\to f$ and $[\psi_n]$ converges in $\T$ to the fixed point of $\sigma_g$\,.
\end{itemize}
The ramification portrait of $g$ translates to relations between the images
of marked points, $x_i=\psi(z_i)$ and $x_i'=\psi'(z_i)$, such that
$g(z_i)=z_j$ implies $f(x_i')=x_j$ when $\psi\circ g=f\circ \psi'$. In the
bicritical case with marked critical points, $f$ is determined by $x_i$ and
the $x_i'$ are determined up to the branch of the $d$-th root. For the formal
mating $P\fmate Q$ with $P(z)=z^2+\i$ and $Q(z)=z^2-1$ according to
Figure~\ref{Fib} and Example~\ref{Xsmib}, this reads
\be\label{eqmibnop}
 x_1'=\pm\sqrt{\frac{x_1-x_2}{1-x_2}} \quad
 x_2'=\pm\sqrt{\frac{2x_1}{1+x_1}} \ . \ee
Note that we cannot pull back marked points with these formulas alone, since
the choice of branch is determined by the pullback in Teichm\"uller space. The
following proposition gives contraction properties of $\sigma_g$\,;
see Section~\ref{2t} for the notation \2 and for the relation between
convergence and Thurston obstructions.

\begin{prop}[Thurston--Selinger]\label{Popnorm}
Consider a Thurston map $g$ of degree $d\ge2$ and the pullback map
$\sigma_g$\,.\\[1mm]
$1$. $\sigma_g$ is weakly contracting with respect to the Teichm\"uller
metric.\\[1mm]
$2$. If $g$ has orbifold type not \2, then some iterate of $\sigma_g$ is
strictly contracting. The contraction is uniform on subsets of $\T$, such
that $\pi(\tau)$ varies in a compact subset of $\M$.\\[1mm]
$3$. $\sigma_g$ is Lipschitz continuous on $\T$ with respect to the
Weil-Petersson metric; a factor is given by $\sqrt{d}$.
\end{prop}

References for the \textbf{proof:} 1: Weak contraction follows from
the definition in terms of a minimal dilation, or since the
Teichm\"uller metric is a Kobayashi metric \cite{book1h}.\\
2. The contraction with respect to the infinitesimal Teichm\"uller metric
is obtained from the dual operator, a push-forward of quadratic differentials.
An explicit integral gives strict contraction, unless there is a specific
type \2 of branch portrait. See \cite{dh, book2h}, and \cite{teich} for the
case of additional marked points. Uniform contraction follows from the fact
that the rational maps $f$ depend only on a finite intermediate cover of $\M$;
cf.~Section~\ref{3p}.\\
3. Apply the Cauchy--Schwartz inequality to the same integral representation
of the push-forward operator, to obtain an estimate in the $L^2$-norm
\cite{ext}.
\mybox

\subsection{A path in moduli space} \label{2p}
The pullback of homeomorphisms $\psi_n$ was easy to define, but it is not
computed easily: repeated pullbacks would be defined piecewise, and solving
the Beltrami equation numerically would be impractical as well. The isotopy
classes in Teichm\"uller space are meant to represent only combinatorial
information anyway: we are interested in the pullback of marked points
$x_i(n)\in\pi(\sigma_g^n([\psi_0]))$ and maps $f_n$\,, and the combinatorial
description is needed to make a finite choice between different possible
preimages. This characterization of the topology has been implemented in
terms of spiders \cite{bfh, hss}, medusas \cite{medusa}, and
triangulations \cite{img}. These contain the necessary information from
Teichm\"uller space without being actual homeomorphisms $\psi_n$\,.

Following Bartholdi--Nekrashevych \cite{bn}
and Buff--Ch\'eritat \cite{cts3}, the following alternative method shall be
discussed. It means that Teichm\"uller space is used explicitly only to
check a suitable initialization of a path in moduli space. Afterward the
path is pulled back simply by choosing preimages from continuity. The
application to matings is discussed in Sections~\ref{5} and~\ref{6}. A few
more applications to rational maps are given in \cite{amate}.
The spider algorithm is implemented with a path in \cite{smate} and further
applications to quadratic polynomials are given; twisted polynomials and
Latt\`es maps are discussed in \cite{pmate} as well.

\begin{prop}[Path in moduli space]\label{Ppath}
Suppose $g$ is a Thurston map of degree $d\ge2$, and there is a continuous path
of homeomorphisms $\psi_t:\hat\C\to\hat\C$, $0\le t\le1$, with
$\psi_0\circ g=f_0\circ\psi_1$ for a rational map $f_0$\,.
So $[\psi_1]=\sigma_g([\psi_0])$.

$1$. Using a suitable normalization, there is a unique path of homeomorphisms
$\psi_t$\,, $0\le t<\infty$, with $\psi_t\circ g=f_t\circ\psi_{t+1}$ for
rational maps $f_t$\,, so $[\psi_{t+1}]=\sigma_g([\psi_t])$. It projects to a
continuous path $\pi([\psi_t])$ in moduli space. Note that
$\sigma_g^n([\psi_0])=[\psi_n]$ for $n\in\N$.

$2$. Suppose that $d=2$, or more generally, that $g$ is bicritical. Normalize
the marked points $x_i(t)\in\pi([\psi_t])$ such that $0$ and $\infty$ are
critical and $1$ is postcritical or marked in addition. Then the path $x_i(t)$
in moduli space is computed for $1\le t<\infty$ by pulling back the initial
segment continuously.
\end{prop}

Probably the statement remains true when $g$ is not bicritical, but the
pullback is less explicit, and I am not sure if it is unique in general. Note
that $[\psi_1]=\sigma_g([\psi_0])$ and an initial path $\psi_t$ is projected
to moduli space. If this condition is neglected by choosing an arbitrary path
from $\pi([\psi_0])$ to $\pi([\psi_1])$, the pullback may correspond not to
$g$ but to some twisted version of it. Conditions for convergence of
$\sigma_g^n([\psi_0])$ are discussed in Section~\ref{2t}; in the case of a
non-\2 orbifold, convergence in Teichm\"uller space is equivalent to
convergence in moduli space, and in both spaces, convergence of the sequence
implies convergence of the path as $t\to\infty$. The situation is more involved
for an orbifold of type \2. The implementation in terms of a piecewise linear
path is discussed in \cite{amate, smate}.

\textbf{Proof:} 1. $\sigma_g$ and $\pi$ are continuous. Marked points never
meet under iterated pullback, so $\psi_{t+1}$ is always defined uniquely
up to M\"obius conjugation.

2. In this normalization, we have $f_t(z)=m_t(z^d)$, and the M\"obius
transformation $m_t$ is determined uniquely from the images of $0,\,1,\,\infty$
at time $t$. The path is pulled back uniquely by
$f_t^{-1}(z)=\sqrt[d]{m_t^{-1}(z)}$, since any coordinate is either constant
$0$ or $\infty$, or the argument of the radical is never passing through
$0$ or $\infty$.
\mybox

\begin{xmp}[Misiurewicz polynomial mates Basilica] \label{Xsmib}
The mating of the Misiurewicz polynomial $P(z)=z^2+\i$ and the Basilica
polynomial $Q(z)=z^2-1$ is illustrated in Figure~\ref{Fib}. Consider the
Thurston Algorithm for the formal mating $g$ with a path according to
Initialization~\ref{Ism} and the radius $R_t=\exp(2^{1-t})$. Rescaled to
$f_t(\infty)=1$, the initialization for $0\le t\le1$ reads
\be
 x_1(t)=-\i/R_t^2 \quad
 x_2(t)=\frac{(1-\i)/R_t^2}{1+(1-t){\rm e}^{-4}} \quad
 x_3(t)=\frac{\i/R_t^2}{1+(1-t)2\i{\rm e}^{-4}} \ . \ee
Note that the normalization $x_3(t)=-x_1(t)$ is satisfied for $t\ge1$ only.
For $t\ge0$ we have the following pullback relation, and the formula for
$x_2(t+1)$ simplifies to (\ref{eqmibnop}) when $t\ge1$:
\be
 x_1(t+1)=\pm\sqrt{\frac{x_1(t)-x_2(t)}{1-x_2(t)}} \quad
 x_2(t+1)=\pm\sqrt{\frac{x_1(t)-x_3(t)}{1-x_3(t)}} \quad
 x_3(t+1)=-x_1(t+1) \ , \ee
where the sign is chosen by continuity. According to Theorem~\ref{Tsmrcc},
the rational maps $f_t$ converge to the rescaled geometric mating
$f(z)=(z^2+2)/(z^2-1)$, so $x_1(t)\to-2$, $x_2(t)\to2$, and $x_3(t)\to2$.
Since two postcritical points are identified, the iteration diverges in
moduli space and in Teichm\"uller space.
\end{xmp}

An alternative interpretation of the path reads as follows: by a standard
technique from algebraic topology, the universal cover of moduli space is
constructed as the space of homotopy classes of paths with a fixed starting
point. So that space is isomorphic to Teichm\"uller space. In this sense, the
pullback of the path is a direct implementation of $\sigma_g$\,, and
information on the dynamics of $\sigma_g$ is available from homotopy classes
of paths. See Section~3.3 in \cite{pmate} for an application.

Sarah Koch \cite{endo} gives criteria on $g$ for the existence of a moduli
space map from $\pi(\sigma_g([\psi]))$ to $\pi([\psi])$, which is a critically
finite map in the same dimension as the moduli space.
See also Section~3.2 in \cite{pmate}. Then the path may be chosen within
the Julia set of the moduli space map, which is easily visualized when it is
one-dimensional \cite{bn}. This happens for a NET map, which has four
postcritical points and only simple critical points \cite{net}.
In the quadratic case of NET maps, a moduli space map exists if at least one
critical point is postcritical, and not when $g$ is a Latt\`es map
of type \2.

\begin{xmp}[Obstructed self-mating] \label{Xsmbb}
For the self-mating of the Basilica polynomial $P(z)=Q(z)=z^2-1$, consider
the radius $R_t=\exp(2^{1-t})$ again, and Initialization~\ref{Ism} reads
$x_1(t)=-1/R_t$ for $0\le t\le1$. The normalization is symmetric under
inversion, and the pullback relation $x_1(t+1)=-\sqrt{-x_1(t)}$ has an explicit
solution in this case, which is given by $x_1(t)=-1/R_t$ for $0\le t<\infty$.
So $x_1(t)\to-1$ as $t\to\infty$, and the rational maps
$f_t(z)=(z^2+x_1(t))/(1+x_1(t)z^2)$ degenerate to a constant map. Note that
there is a moduli space map $x_1(t)=-\Big(x_1(t+1)\Big)^2$, and for a different
initialization, the path would be contained in the unit circle.
\end{xmp}

\subsection{Obstructions and the Thurston Theorems} \label{2t}
A \textbf{multicurve} is a nonempty union of pairwise disjoint and
non-homotopic essential curves, or a homotopy class as well.
For a concrete multicurve $\Gamma$, the preimage under a Thurston map $g$ is
$\Gamma'\cup\Gamma''$, where $\Gamma'$ consists of essential curves and the
curves in $\Gamma''$ are peripheral or trivial. The homotopy class of
$\Gamma'$ depends only on the homotopy classes of $\Gamma$ and $g$. Note that
$\Gamma'$ may contain mutually homotopic curves and it may be empty as well.
$\Gamma$ is called invariant, if every curve in $\Gamma'$ is homotopic to a
curve in $\Gamma$; it is completely invariant if the converse holds in
addition.

For $\Gamma=\{\gamma_1\,,\,\dots\,,\,\gamma_n\}$ the Thurston matrix
$M_\Gamma=(m_{ij})$ is defined as $m_{ij}=\sum1/d_{ijk}$\,, where the sum runs
over all preimages of $\gamma_j$ homotopic to $\gamma_i$ and $d_{ijk}$ is the
degree of $g$ on these preimages. Now $\Gamma$ is a Thurston
\textbf{obstruction}, if the leading eigenvalue $\lambda_\Gamma$ of $M_\Gamma$
satisfies $\lambda_\Gamma\ge1$. There are different conventions, whether
invariance is required. Most important are the following kinds of obstructions:
\begin{itemize}
\item An invariant multicurve $\Gamma$ is a \textbf{simple obstruction}, if no
permutation turns $M_\Gamma$ into a lower-triangular block form, such that the
upper left block has leading eigenvalue $<1$. A simple obstruction is always
completely invariant.
\item A multicurve $\Gamma=\{\gamma_1\,,\,\dots\,,\,\gamma_n\}$ is a
\textbf{L\'evy cycle}, if each $\gamma_i$ is homotopic to a preimage of
$\gamma_{i+1\mathop{mod} n}$ and the corresponding degree is one. Then $\Gamma$
need not be invariant, but it can be extended to a simple obstruction. The
converse holds in the quadratic or bicritical case: every simple obstruction
contains a L\'evy cycle. These are classified further in \cite{rst, rees1}.
\end{itemize}
Obstructions are important for the Thurston pullback, because they are related
to the presence of annuli with large modulus and of short geodesics
\cite{dh, book2h, teich}; see also Section~\ref{32} for more explicit
statements.
Note that by adding more additional marked points, $g$ may have more
obstructions. 
The notion of an orbifold is explained in \cite{bookm, book1mcm}. A Thurston
map of orbifold type \2 has four postcritical points, the critical points are
non-degenerate and not postcritical.

\begin{thm}[Thurston--Pilgrim, general case]\label{TT1}
Suppose $g$ is a Thurston map of degree $d\ge2$ with orbifold not of type \2,
possibly with additional marked points. Then:\\
$\bullet$ Either there is no Thurston obstruction, $g$ is combinatorially
equivalent to a rational map $f$, which is unique up to M\"obius conjugation,
and the Thurston pullback $\sigma_g$ converges globally to its unique fixed
point.\\
$\bullet$ Or $g$ is obstructed, it is not equivalent to a rational map, and the
Thurston pullback diverges. There is a unique \textbf{canonical obstruction}
$\Gamma$, such that every $\gamma\in\Gamma$ is pinched, while the length of
every other curve is bounded uniformly from below. 
\end{thm}

Note that the fine print reads as follows: if there is an obstruction,
it need not be pinching, but then it implies the existence of another
obstruction, which is pinching. The pinching obstructions imply that marked
points get identified in the limit, which means divergence to the boundary in
Teichm\"uller space and in moduli space as well. Then the rational maps may
converge to a rational map of degree $d$ or of lower degree.
Note that there is no algorithm to find obstructions in a time bounded
a priori. 

\textbf {Idea of the proof:} According to Proposition~\ref{Popnorm}.2, some
iterate of $\sigma_g$ is strictly contracting. So if $g$ is equivalent to a
rational map $f$, $\sigma_g$ has a unique and globally attracting fixed point.
If $g$ was obstructed, there would be a simple obstruction $\Gamma$. Choosing
$\tau_0$ with suitable annuli of large modulus, these would not shrink under
the pullback, giving a contradiction to the smaller annuli for $f$. Note also
that a rational map is obstructed only if it is of flexible Latt\`es type, or
postcritically infinite with a rotation domain \cite{book1mcm}.

On the other hand, if $\sigma_g$ does not have a fixed point, the pullback
$\tau_n=\sigma_g^n(\tau_0)$ has the property that $\pi(\tau_n)$ leaves every
compact subset of $\M$; there will be short geodesics by Mumford compactness
according to Proposition~\ref{Pbptms}.7. Now a combinatorial analysis gives
a simple obstruction from curves sufficiently short and shorter than other
curves. See \cite{dh, book2h}, and \cite{teich} for the case of additional
marked points. The proof was refined by Pilgrim \cite{kpcan} to show that the
same curves either do get shorter or stay bounded. For an alternative proof by
Selinger \cite{ext}, see Section~\ref{32}. \mybox

Originally this theorem was stated under the assumption of a hyperbolic
orbifold \cite{dh, book2h}. A parabolic orbifold is either of type \2 as
discussed below, or there are only two or three postcritical points. In that
case, the Thurston pullback is undefined or constant, respectively, unless
there are additional marked points. When there are additional marked points,
the pullback map is strictly contracting as in the hyperbolic orbifold case
\cite{teich}.

For quadratic Thurston maps of type \2, some things are the
same, some are different: a fixed point is still unique, when it exists, but it
is not attracting. Every obstruction is pinching, and it excludes a fixed
point, but there are unobstructed maps not equivalent to a rational map as
well. The converse happens when the degree $d\ge4$ is a square: there is a
one-parameter family of flexible Latt\`es maps, which are mutually equivalent
but not M\"obius conjugate. So uniqueness fails, and moreover, there is a
non-pinching obstruction. For Thurston maps of type \2 with
additional marked points, pinching and non-pinching obstructions are
characterized by Selinger--Yampolsky \cite{char, sy}.

\begin{thm}[Thurston, exceptional case]\label{TT2}
$1$. Suppose $g$ is a Thurston map of degree $d\ge2$ with orbifold
of type \2, without additional marked points. Then it is equivalent to a map
covered by a real-affine map on a torus, which is described by an integer
matrix of determinant $d$:\\
$\bullet$ If the eigenvalues are not real, there is an equivalent rational map
 $f$, unique up to M\"obius conjugation. The Thurston pullback map $\sigma_g$
 is represented by a M\"obius transformation of the upper halfplane. It has a
 neutral fixed point, so the Thurston pullback does not converge.\\
$\bullet$ If the matrix is a multiple of the identity, $g$ is equivalent to
 a family of flexible Latt\`es maps.\\
$\bullet$ Otherwise, $g$ is not equivalent to a rational map.

$2$. Suppose $g$ is quadratic of type \2, without additional marked points,
and the integer matrix has trace $t$\,:\\
$\bullet$ If $|t|\le2$, then $g$ is equivalent to a rational map $f$, unique
up to M\"obius conjugation, and $\sigma_g$ has a neutral fixed point.\\
$\bullet$ If $|t|=3$, then $g$ has a pinching obstruction, and the Thurston
pullback diverges to the boundary in Teichm\"uller space and in moduli space
as well.\\
$\bullet$ If $|t|\ge4$, then $g$ is not equivalent to a rational map either,
but there is no Thurston obstruction for $g$. The Thurston pullback diverges
to the boundary in Teichm\"uller space, but it is bounded in moduli space.
\end{thm}

\textbf {Sketch of the proof:} 1. The lift to an affine map is explained in
\cite{dh, book2h, book1mcm, sy, bookbm}, see also
\cite{pmate}. In the flexible Latt\`es case, the choice of lattice
is arbitrary, and in the case of non-real eigenvalues, the lattice can be
chosen such that the real-affine map is holomorphic. The pullback of a constant
Beltrami coefficient is described explicitly by a M\"obius transformation of
the upper halfplane with integer coefficients; a fixed point cannot be
attracting, since otherwise the complex conjugate fixed point would be
repelling.

2. The eigenvalues are computed from $\eta^2-t\eta+2=0$. When $|t|\ge4$,
they are real and not integer, so there is no invariant multicurve according
to \cite{char}. Note also that any possible obstruction would be a
single curve $\gamma$ separating the critical values from the other two marked
points. Its two preimages are mutually homotopic, so we have a pinching
obstruction when they are homotopic to $\gamma$ as well. This cannot happen
when the two critical values are mapped to the same prefixed point, which
is equivalent to $|t|$ being even. See also the example of twisting the
self-mating $f\eqg1/4\tmate1/4$ in Section~3.3 of \cite{pmate}.
\mybox

--- Dylan P.~Thurston \cite{dt1} gives a positive criterion for $g$ to be
equivalent to a rational map, at least if there is a periodic critical point:
$f^{-k}$ is uniformly contracting on a graph, which forms a spine for
$\hat\C\setminus P$. In \cite{kpao}, Kevin Pilgrim gives an algebraic
characterization of obstructions by non-contraction of the virtual endomorphism
of the pure mapping class group. See \cite{bn, bd0} for algebraic descriptions
of Thurston maps in terms of iterated monodromy groups or bisets.

\section{Extension to the augmented Teichm\"uller space} \label{3}
When a Thurston map $g$ is not combinatorially equivalent to a rational map, so
$\sigma_g$ has no fixed point in $\T$, we may understand this by
considering a space larger than $\T$ or a different topology. Except for
orbifold type \2, divergence of the Thurston Algorithm is related to collisions
of marked points and pinching of essential curves. These phenomena are
described by strata of augmented Teichm\"uller space $\hat\T$ and augmented
moduli space $\hat\M$, which parametrize noded Riemann surfaces. Selinger
\cite{ext, char} has extended $\sigma_g$ to $\hat\T$ and obtained a
related characterization of the canonical obstruction $\Gamma$. --- A somewhat
informal introduction to augmented spaces is given in Section~\ref{31}. Results
on extension and characterization are discussed in Sections~\ref{32}--\ref{33},
and applied to obtain a relation between core entropy and matings of conjugate
polynomials. The main result of the present paper is a convergence principle
for obstructed maps, which is proved in Sections~\ref{3p}--\ref{3s}.\\
\textit{Special thanks to Nikita Selinger for
patiently answering my impatient questions.}

\subsection{Augmented Teichm\"uller space and moduli space} \label{31}
Augmented moduli space $\hat\M$ describes noded Riemann surfaces, and augmented
Teichm\"uller space $\hat\T$ has continuous maps from a topological sphere to a
noded sphere; these send certain curves to single points. Boundary strata of
$\hat\T$ are products of lower-dimensional Teichm\"uller spaces. The
notion of noded Riemann surfaces is motivated here by pinching obstructions;
originally they were introduced to compactify $\M$, and to describe
algebraic curves with self-intersections. These constructions are due
to Deligne--Mumford, Bers, Abikoff, and Masur; see the references in
\cite{hkat, wwp, mwp}. The following example shows the degeneration of a
Riemann surface with boundary explicitly: 

\begin{xmp}[Pinching a short geodesic] \label{Xnode}
Consider the Riemann surface
$S_t=\{(x,\,y)\in\C^2\,|\,|x|<1,\,|y|<1,\,x\cdot y=t\}$ for $|t|<1$.
As $t\to0$, $S_t$ becomes the union of two disks intersecting transversely
in the single point $(0,\,0)$. The hyperbolic metric in the annulus
$|t|<|x|<1$ is known explicitly, and seen to converge to the hyperbolic metric
of the punctured disk. When we try to illustrate this process in $\R^2$ or
$\R^3$, either $S_t$ looks disconnected or the limit does not show a
transversal intersection of smooth manifolds.
\end{xmp}

This example shall motivate that we are interested in surfaces consisting of
smooth spheres intersecting transversely. The nodes appear as additional
marked points in the pieces, because the hyperbolic metric is singular there.
An approximate Riemann surface would have long, thin tunnels between thick
components; this is symbolized by connected spheres as in Figure~\ref{Fmab}.

\begin{xmp}[Augmented moduli space] \label{Xdistmod}
1. Consider $\hat\C$ with four marked  points
$x_1=\infty,\,x_2=0,\,x_3=1,\,x_4=a$.
The moduli space is given by $a\in\hat\C\setminus\{\infty,\,0,\,1\}$ and the
augmented moduli space is $\hat\M=\hat\C$: e.g., $a\to0$ corresponds to
pinching a curve separating $x_2=0$ and $x_4=a$ from $x_3=1$ and $x_1=\infty$.
Note that division by $a$ gives a different normalization
$x_1=\infty,\,x_2=0,\,x_3=1/a,\,x_4=1$; now $a\to0$ means $x_3\to x_1$. In
fact this is the same Riemann surface as before, with a node
separating $x_2\,,\,x_4$ from $x_3\,,\,x_1$\,.

2. The case of five marked points is described by
$\M=\{(a,\,b)\in\hat\C^2\,|\,a,\,b\neq\infty,\,0,\,1,\,a\neq b\}$, but now the
topology of $\hat\M$ is more involved than $\hat\C^2$: E.g., one-dimensional
boundary strata are given by $a=0,\,b\neq\infty,\,0,\,1$ or by
$a=b\neq\infty,\,0,\,1$, but we lose information when $a=b=0$: this may be one
of three $0$-dimensional strata, or in a one-dimensional stratum without
information on the relative position of three marked points and a node.
\end{xmp}

So $\hat\T$ contains maps from a topological sphere to a possibly noded
Riemann surface, sending certain curves to nodes, equivalent under an isotopy
in the domain or M\"obius transformations in the range. Boundary strata
$\S_\Gamma\subset\hat\T$ are labeled by homotopy classes of pinched
multicurves. There are only finitely many boundary strata
$\S_{G\cdot\Gamma}\subset\hat\M$; as in the case of $\T$ and $\M$, different
curves in the topological sphere may be mapped to the same short geodesic. A
neighborhood basis for the topology of $\hat\T$ or $\hat\M$ is defined in
terms of maps between noded Riemann surfaces, which map $\eps$-short geodesics
to nodes and which are $(1+\eps)$-quasiconformal outside of the collars. Using
a combination of Fenchel--Nielsen coordinates \cite{book1h} and plumbing
coordinates, which take Example~\ref{Xnode} as a local model, the
infinitely branched covering $\pi:\hat\T\to\hat\M$ is understood locally.
So $\hat\M$ is a compact analytic space \cite{hkat}.
%
%

\begin{prop}[Augmented Teichm\"uller space]\label{Paug}
$\hat\T$ and $\hat\M$ are topological spaces, such that $\pi:\T\to\M$ extends
to a continuous map $\pi:\hat\T\to\hat\M$.

$1$. The Weil--Petersson metric $d_\sWP$ extends to $\hat\T$ and $\hat\M$, such
that $\hat\T$ is the completion of $\T$ and $\hat\M$ is a compactification of
$\M$. On each boundary stratum, the extended $d_\sWP$ is a product of
lower-dimensional Weil--Petersson metrics.

$2$. Each point $\tau\in\hat\T$ is approximated only from finitely many strata.

$3$. Normalizing three marked points, the coordinates of all marked points
extend continuously from $\M$ to $\hat\M$ or from $\T$ to $\hat\T$.

$4$. For every essential simple closed curve $\gamma$, the length
$l(\gamma,\,\tau)\in[0,\,\infty]$ is continuous on $\hat\T$. All length
functions together determine $\tau$ uniquely.
\end{prop}

Explanations, references, and sketch of a \textbf{proof}:
See \cite{hkat, wwp, mwp} for item~1. Note that $\hat\T$ is only a
partial compactification of $\T$: a sequence leaving $\T$ may converge in
$\hat\T$, if closed curves are pinched, but a sequence of the form
$g^k\cdot\tau_0$ will diverge in $\hat\T$ as well.

2. This follows from the definition of neighborhoods given above. 

3. Continuity is obtained from extending the $(1+\eps)$-quasiconformal maps
into approximately round collars, or from a compactness argument and continuity
of length. 
The normalization singles out a sphere, where all marked points have limits,
while marked points in other components converge to nodes of this sphere.
The statement is equivalent to a continuous extension of
cross-ratios; in \cite{FT11} a completion with respect to cross-ratios is
used to construct a space isomorphic to $\hat\T$.

4. See the references above and \cite{ext}. Approaching a lower-dimensional
stratum according to item~2, specific curves have length $\to0$ and
intersecting curves have length $\to\infty$. For all other curves, the
hyperbolic metric converges on each component in a suitable normalization.
Injectivity of lengths follows from Proposition~\ref{Pbptms}.4. \mybox

The pure mapping class group $G$ acts on $\hat\T$ by Weil--Petersson
isometries, but the description of $\hat\M=\hat\T/G$ is more involved:

\begin{rmk}[Action of $G$ and uniqueness of geodesics]\label{Raug}
Near a boundary stratum $\S_\Gamma\subset\hat\T$, distinguish the following
kinds of Dehn twists $g\in G$ about $\gamma$:\\[1mm]
a) If $\gamma$ intersects a curve in $\Gamma$, the action of $g$ would
map to a different stratum.\\[1mm]
b) If $\gamma$ is contained in a component of $\hat\C\setminus\Gamma$,
then $g$ acts from the pure mapping class group of the component space.\\[1mm]
c) If $\gamma\in\Gamma$, then $g$ acts trivially on the stratum but not
trivially in a neighborhood.
For $\tau\in\T$ close to $\S_\Gamma$, $d_\sWP(\tau,\,g^k\cdot\tau)$ is
bounded by the triangle inequality, but $d_\sT(\tau,\,g^k\cdot\tau)\to\infty$,
although all hyperbolic geodesics are bounded below;
cf.~Proposition~\ref{Pbptms}.\\[1mm]
So $\pi$ is infinitely branched and not a
local isometry. According to \cite{wwp}, $\hat\T$ is a unique geodesic
space nevertheless, with open geodesics passing through a unique stratum
of lowest possible dimension.
\end{rmk}

\subsection{Extended pullback map and the canonical obstruction} \label{32}
\begin{center}\textit{Well, the extension to the augmented Teichm\"uller space
exists independently of how or whether we understand it.}
--- Nikita Selinger [private communication] 
\end{center}

To extend the pullback map $\sigma_g$ to augmented Teichm\"uller space,
consider a multicurve $\Gamma$ and the collection $\Gamma'$ of non-homotopic
essential preimages. Then $\sigma_g$ shall map $\S_\Gamma$ to $\S_{\Gamma'}$\,.
The definition is understood by considering the full preimage $g^{-1}(\Gamma)$
first; this defines noded surfaces with possibly non-hyperbolic pieces. So
whenever a disk component contains at most one marked point, it is reduced to
a point, and an annulus between homotopic curves is reduced to a point as well.
This process of stabilization defines a noded Riemann surface with hyperbolic
pieces.

\begin{figure}[h!t!b!]
\unitlength 0.001\textwidth 
\begin{picture}(990, 420)
\put(10, 0){\includegraphics[width=0.42\textwidth]{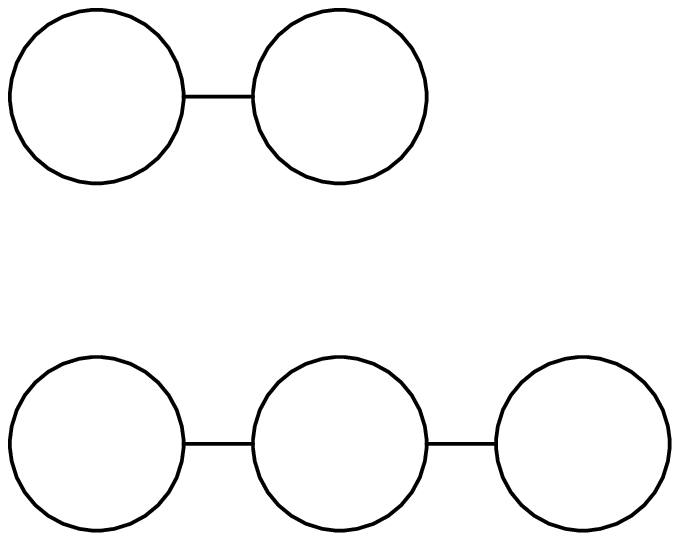}}
\put(570, 0){\includegraphics[width=0.42\textwidth]{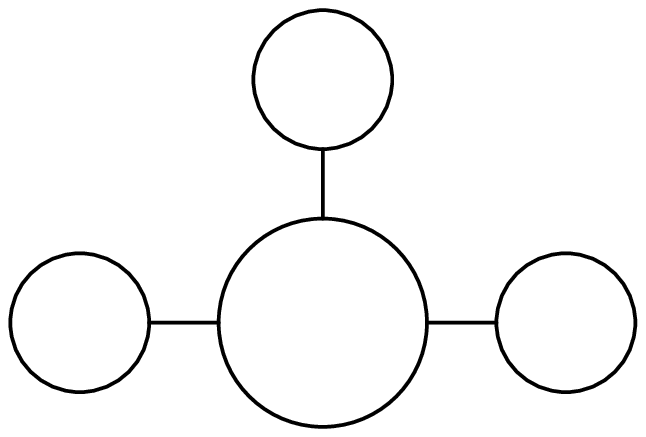}}
\thinlines
\multiput(10, 0)(560, 0){2}{\line(1, 0){420}}
\multiput(10, 420)(560, 0){2}{\line(1, 0){420}}
\multiput(10, 0)(560, 0){2}{\line(0, 1){420}}
\multiput(430, 0)(560, 0){2}{\line(0, 1){420}}
\put(80, 75){\makebox(0, 0)[cb]{$z_1$}}
\put(80, 145){\makebox(0, 0)[ct]{$w_1$}}
\put(220, 75){\makebox(0, 0)[cb]{$z_0$}}
\put(220, 145){\makebox(0, 0)[ct]{$w_0$}}
\put(360, 75){\makebox(0, 0)[cb]{$z_2$}}
\put(150, 125){\makebox(0, 0)[cb]{$\alpha$}}
\put(290, 125){\makebox(0, 0)[cb]{$-\alpha$}}
\put(80, 275){\makebox(0, 0)[cb]{$z_1$}}
\put(80, 345){\makebox(0, 0)[ct]{$w_1$}}
\put(220, 275){\makebox(0, 0)[cb]{$z_0$}}
\put(220, 345){\makebox(0, 0)[ct]{$w_0$}}
\put(255, 310){\makebox(0, 0)[rc]{$z_2$}}
\put(150, 325){\makebox(0, 0)[cb]{$\alpha$}}
\put(780, 95){\makebox(0, 0)[cb]{$z_0$}}
\put(640, 115){\makebox(0, 0)[cb]{$z_2$}}
\put(640, 165){\makebox(0, 0)[ct]{$w_1$}}
\put(780, 255){\makebox(0, 0)[cb]{$w_0$}}
\put(780, 305){\makebox(0, 0)[ct]{$z_1$}}
\put(920, 115){\makebox(0, 0)[cb]{$z_3$}}
\put(920, 165){\makebox(0, 0)[ct]{$w_2$}}
\end{picture} \caption[]{\label{Fmab}
Left: The formal mating $g$ of Airplane and Basilica has a cyclic ray
connection $\Gamma=\{\gamma\}$ between the two $\alpha$-fixed points, which is
the canonical obstruction. The preimage $\Gamma'\cup\Gamma''$ contains a
peripheral curve, so the right sphere in the
lower surface $\hat\C/\Gamma'$ is considered as one point. Then $g_\Gamma$ is
a self-map of the noded surface $\hat\C/\Gamma$ and $z_0$ is no longer
$3$-periodic $z_0\Rightarrow z_1\to z_2\to z_0$\,, but preperiodic
with $z_0\Rightarrow z_1\to z_2\to\alpha\uparrow$.\\
Right: The formal mating $1/4\fmate1/2$ with marked critical points has
three curves in the canonical obstruction, which surround ray connections
with the angles $1/4,\,1/2,\,0$. Due to the identification of points in the
small pieces, the geometric mating $f\eqg1/4\tmate1/2$ satisfies
$f(0)=\infty$.}
\end{figure}

So when $\sigma_g([\psi])=[\psi']$ and $\psi$ maps the curves of $\Gamma$ to
nodes, then $\psi'$ maps curves in $\Gamma'$ to nodes, and certain annuli and
disks to nodes or marked points as well. The pullback map on a product of
lower-dimensional Teichm\"uller spaces is described in terms of homeomorphisms
or Thurston maps $g^C$ on the pieces $C$, which are defined uniquely up to
combinatorial equivalence. Note that for a completely invariant multicurve
$\Gamma$, appropriate identifications must be made to describe a
$\sigma_g$-invariant boundary stratum $\S_\Gamma$\,. The collection $g_\Gamma$
of component maps is defined on the topological model surface $\hat\C/\Gamma$.
The examples in Figure~\ref{Fmab} show that it may be discontinuous, not
surjective, and it may map marked points to nodes.

\begin{thm}[Selinger extension]\label{Text}
For a Thurston map $g$ of degree $d\ge2$, the Thurston pullback $\sigma_g$ has
a unique continuous extension to $\hat\T$. On each boundary stratum, it is
given by a pullback with component maps as described above.
\end{thm}

\textbf{Idea of the proof:} A unique extension is given by completion, using
the uniform Lipschitz estimate from Proposition~\ref{Popnorm}.3. For the
explicit extension above, the length functions of all geodesics are continuous
when a lower-dimensional stratum is approximated from a higher-dimensional
stratum \cite{ext}. By Proposition~\ref{Paug}.4, both extensions agree.
\mybox

The following result is due to Pilgrim \cite{kpcan} in the case of a hyperbolic
orbifold. The proof by Selinger \cite{ext} works for maps of type \2 as well.

\begin{thm}[Canonical obstruction by Pilgrim--Selinger]\label{Tcan}
Suppose $g$ is a Thurston map of degree $d\ge2$, fix $\tau_0\in\T$, and
set $\tau_n=\sigma_g^n(\tau_0)$. There is an $R(\tau_0)>0$ and a
multicurve $\Gamma$, possibly empty, such that:\\
$\bullet$ If $\gamma\in\Gamma$, then $l(\gamma,\,\tau_n)\to0$.\\
$\bullet$ If $\gamma\notin\Gamma$, then $l(\gamma,\,\tau_n)\ge R(\tau_0)$.\\
In augmented moduli space, $\pi(\tau_n)$ accumulates at a compact subset of the
canonical stratum $\S_{G¢\cdot\Gamma}\subset\hat\M$. The
\textbf{canonical obstruction} $\Gamma$ is independent of $\tau_0$\,. If
$\Gamma\neq\emptyset$, it is a simple Thurston obstruction, and the curves of
$\Gamma$ do not intersect another curve from any simple obstruction.
\end{thm}

This implies that every accumulation point of $\tau_n$ belongs to the
canonical stratum $\S_\Gamma\subset\hat\T$, but there need not be
accumulation in $\hat\T$ at all. 
The accumulation sets in $\hat\M$ and $\hat\T$ may depend on
the starting point $\tau_0$\,. While multicurves and obstructions are never
empty, it is customary to say ``$\Gamma$ is empty'' instead of
``there is no $\Gamma$'' here.
Recall the definition of the Thurston matrix $M_\Gamma$ from Section~\ref{2t}.
We have:
\begin{itemize}
\item If $\Gamma$ is a simple obstruction, $M_\Gamma$ has a positive
eigenvector $v$ with eigenvalue $\lambda_\Gamma\ge1$. Suppose that in the
Riemann surface $\pi(\tau)$, there are annuli around the corresponding
geodesics with moduli proportional to $v$, then these moduli will grow at
least by $\lambda_\Gamma$ under the pullback $\pi(\sigma_g(\tau))$. This is
a direct application of the Gr\"otsch inequality \cite{book1h}. By the
collar estimate from Proposition~\ref{Pcoll}.3, a lower bound on the modulus
corresponds to an upper bound on the hyperbolic length of the geodesic.
\item For sufficiently short geodesics, there is a kind of reverse estimate:
when $\Gamma$ is completely invariant but not a simple obstruction, there is a
semi-norm on the vector of inverse lengths, which cannot increase arbitrarily.
Here a preimage annulus is decomposed along parallels to the core curve, and
the new annuli are related to inverse length by the collar theorem again. See
Theorem~7.1 in \cite{dh}, Theorem~10.10.3 in \cite{book2h},
or Lemma~2.6 in \cite{teich}.
\end{itemize}
Sketch of the \textbf{proof} of Theorem~\ref{Tcan}: Let $N\ge0$ be the
maximal number of arbitrarily short geodesics in $\pi(\tau_n)$  as
$n\to\infty$. So there is $R>0$, a subsequence $n_k$\,, multicurves
$\Gamma_k$ with $N$ elements, and $\eps_k\to0$, such that: \\[1mm]
$\bullet$ For each $n$, there are at most $N$ curves $\gamma$ with
$l(\gamma,\,\tau_n)<R$.\\[1mm]
$\bullet$ For all $k$ and all $\gamma\in\Gamma_k$ we have
$l(\gamma,\,\tau_{n_k})<\eps_k$\,.\\[1mm]
Now if $N=0$, the claims are satisfied for $\Gamma=\emptyset$. So assume
$N>0$. For large $k$ we have $\eps_k<<R$ and continuity of $l(\gamma,\cdot)$
together with the reverse inequality above shows that $\Gamma_k$ is completely
invariant. We may assume that all $\Gamma_k$ have the same partition of
marked points $G\cdot\Gamma_k$, the same Thurston matrix $M=M_{\Gamma_k}$\,,
and $\pi(\tau_{n_k})$ has a limit in $\S_{G\cdot\Gamma_k}\subset\hat\M$. If
$\Gamma_k$ was not a simple obstruction, the reverse inequality applied to $M$
would give a lower bound for $\eps_k$\,. This is a contradiction for some
large $k$, and we set $\Gamma=\Gamma_k$\,. Now for all $n\ge n_k$ the
moduli of annuli around $\gamma\in\Gamma$ have non-decreasing lower bounds
and the lengths $l(\gamma,\,\tau_n)$ have non-increasing upper bounds; together
with the assumptions on the subsequence this gives $l(\gamma,\,\tau_n)\to0$.
The lower bound $R$ is satisfied by all other curves, and Mumford
compactness according to Proposition~\ref{Pbptms}.7 applies to all components
of $\S_{G\cdot\Gamma}$. If some $\gamma\in\Gamma$ was intersecting a
$\gamma'$ from another simple obstruction homotopically transversely, there
would be an upper bound for $l(\gamma',\,\tau_n)$ and a lower bound for
$l(\gamma,\,\tau_n)$ by the collar theorem. Finally, for a different
initial $\tau_0$ all length $l(\gamma,\,\tau_n)$ are changed by a factor
bounded above and below, so $\Gamma$ is independent of $\tau_0$\,. \mybox

An alternative \textbf{proof of the Thurston Theorem~\ref{TT1}} based on
Theorem~\ref{Tcan}: Consider $\tau_n=\sigma_g^n(\tau_0)$ for some
$\tau_0\in\T$. If the canonical obstruction is $\Gamma\neq\emptyset$, then
$\pi(\tau_n)$ leaves every compact subset of $\M$, so $\sigma_g$ cannot have
a fixed point in $\T$ and there is no rational map $f$ equivalent to $g$.

Now assume $\Gamma=\emptyset$. Then $\pi(\tau_n)$ stays in a compact subset
of $\M$. If $g$ is of type \2, it may be obstructed or not, equivalent to a
rational map or not. But otherwise some iterate of $\sigma_g$ is uniformly
contracting over the compact set of $\M$ defined by $R(\tau_0)$. So $\tau_n$
converges to a fixed point, which corresponds to a rational map $f$. If $g$
had a non-canonical simple obstruction $\tilde\Gamma$, then $\tau_n$ could
not converge to the fixed point if $\tau_0$ had $l(\gamma,\,\tau_0)$ too
small for $\gamma\in\tilde\Gamma$. \mybox

\begin{figure}[h!t!b!]
\unitlength 0.001\textwidth 
\begin{picture}(990, 420)
\put(10, 0){\includegraphics[width=0.42\textwidth]{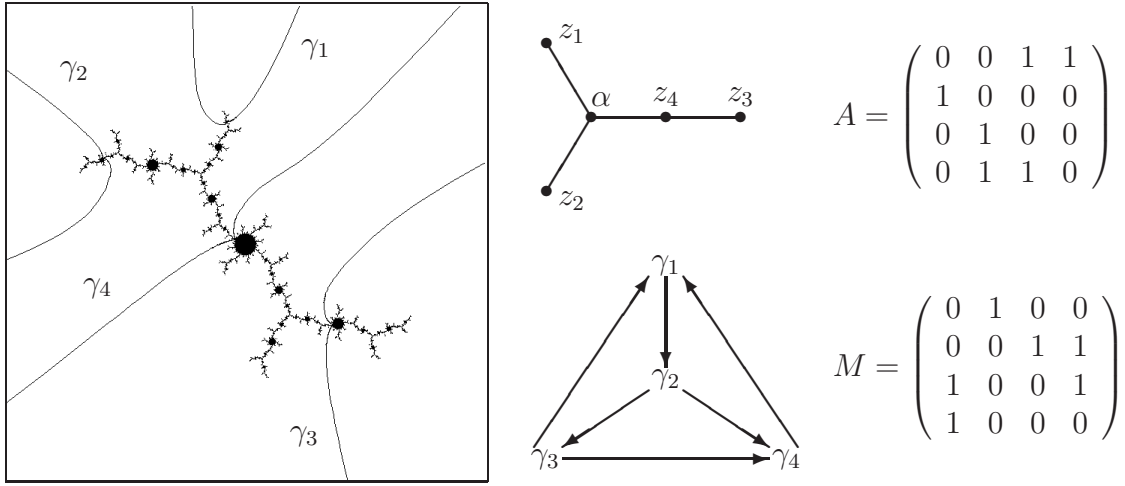}}
\thinlines
\put(10, 0){\line(1, 0){420}}
\put(10, 420){\line(1, 0){420}}
\put(10, 0){\line(0, 1){420}}
\put(430, 0){\line(0, 1){420}}
\put(280, 380){\makebox(0, 0)[cc]{$\gamma_1$}}
\put(70, 360){\makebox(0, 0)[cc]{$\gamma_2$}}
\put(270, 40){\makebox(0, 0)[cc]{$\gamma_3$}}
\put(90, 170){\makebox(0, 0)[cc]{$\gamma_4$}}
\thicklines
\put(470, 230){\begin{picture}(230, 180)
\put(50, 90){\line(-3, 5){39}}
\put(50, 90){\line(-3, -5){39}}
\put(50, 90){\line(1, 0){130}}
\put(11, 25){\circle*{10}} \put(21, 25){\makebox(0, 0)[lt]{$z_2$}}
\put(11, 155){\circle*{10}} \put(21, 155){\makebox(0, 0)[lb]{$z_1$}}
\put(180, 90){\circle*{10}} \put(180, 100){\makebox(0, 0)[cb]{$z_3$}}
\put(115, 90){\circle*{10}} \put(115, 100){\makebox(0, 0)[cb]{$z_4$}}
\put(50, 90){\circle*{10}} \put(50, 100){\makebox(0, 0)[lb]{$\alpha$}}
\end{picture}}
\put(465, 10){\begin{picture}(240, 190)
\put(120, 170){\vector(0, -1){80}}
\put(105, 70){\vector(-3, -2){75}}
\put(135, 70){\vector(3, -2){75}}
\put(5, 20){\vector(2, 3){100}}
\put(235, 20){\vector(-2, 3){100}}
\put(30, 10){\vector(1, 0){180}}
\put(120, 180){\makebox(0, 0)[cc]{$\gamma_1$}}
\put(120, 80){\makebox(0, 0)[cc]{$\gamma_2$}}
\put(15, 10){\makebox(0, 0)[cc]{$\gamma_3$}}
\put(225, 10){\makebox(0, 0)[cc]{$\gamma_4$}}
\end{picture}}
\put(730, 320){\makebox(0, 0)[lc]{$A=\left(\begin{array}{cccc}
0 & 0 & 1 & 1\\
1 & 0 & 0 & 0\\
0 & 1 & 0 & 0\\
0 & 1 & 1 & 0
\end{array}\right)$}}
\put(730, 100){\makebox(0, 0)[lc]{$M=\left(\begin{array}{cccc}
0 & 1 & 0 & 0\\
0 & 0 & 1 & 1\\
1 & 0 & 0 & 1\\
1 & 0 & 0 & 0
\end{array}\right)$}}
\end{picture} \caption[]{\label{Fkk}
The formal mating of the Kokopelli $p=\gamma_\sM(3/15)$ with its conjugate
$q=\overline p$. The $4$-periodic rays define four loops $\gamma_i$ and a noded
surface $\hat\C/\Gamma$ with five pieces. The Thurston matrix $M$ is the
transpose of the Markov matrix $A$, which describes the mapping of edges in the
Hubbard tree $T_p$\,. The leading eigenvalue $\lambda=1.395337$ with
$\lambda^4-2\lambda-1=0$ determines the core entropy
$h(p)=\log\lambda$.\\
Under pullback with $g=P\fmate Q$, the three L\'evy cycles
converge to ray-equivalence classes:
$(\gamma_1\,,\,\gamma_2\,,\,\gamma_3\,,\,\gamma_4)$ gives the original
$4$-periodic rays, $(\gamma_1\,,\,\gamma_2\,,\,\gamma_3)$ converges to a
$3$-cycle of loops with $6$-periodic rays, and
$(\gamma_1\,,\,\gamma_2\,,\,\gamma_4)$ gives the $3$-periodic rays from
$\alpha_p$ to $\alpha_q$\,.}
\end{figure}

\subsection[Augmented Teichm\"uller space and canonical obstructions]%
{Characterization of the canonical obstruction} \label{33}
There is no algorithm to determine obstructions of an arbitrary Thurston
map $g$, but if you have a guess what the canonical obstruction $\Gamma$ might
be, this can be checked with the following criterion. In particular, it shows
that the canonical obstruction of a formal mating from non-conjugate limbs is
given by loops around ray-equivalence classes with at least two postcritical
points; see Section~\ref{4c}. Examples of canonical obstructions are given in
Figures~\ref{Fmab} and~\ref{Fkk} as well.

\begin{thm}[Selinger characterization of the canonical obstruction]
\label{Tchar}
When $g$ is a Thurston map of degree $d\ge2$, consider the family of
multicurves $\tilde\Gamma$, which are simple obstructions or empty, with the
following property: for the map $g_{\tilde\Gamma}$ between components of the
noded surface defined by $\tilde\Gamma$, the first-return map of each periodic
component is\\
$\bullet$ a homeomorphism,\\
$\bullet$ an unobstructed Thurston map, or\\
$\bullet$ a \2-map with a non-pinching obstruction: all curves are essential
with respect to the four postcritical points, and the degree of the map
is a square.\\
Now this family of multicurves $\tilde\Gamma$ has a unique minimal element
with respect to inclusion, which is the canonical obstruction $\Gamma$.
\end{thm}

Idea of the \textbf{proof} from \cite{char}: First, suppose that a first-return
map of $\hat\C/\Gamma$ is obstructed, then $g$ has a non-canonical obstruction
$\Gamma'$. Suitable annuli of maximal modulus have the property that these
moduli are increasing and bounded above. A subsequence of rational maps
converges to a limit map on the component in a suitable normalization. This map
has annuli of invariant maximal modulus, so the subdivision of preimages
happens parallel to the core curves; this fact is used to obtain type \2. Now
obstructions are related to integer eigenvalues of the corresponding matrix
lift, and if the degree was not a square, these eigenvalues would be different
and have a quotient $>1$. But then the Thurston matrix of $\Gamma'$ would have
$\lambda_{\Gamma'}>1$ and $\Gamma'$ would be pinching for $g$ as well.

So $\Gamma$ satisfies the assumptions on $\tilde\Gamma$ in Theorem~\ref{Tchar}.
Conversely, we must see that any simple obstruction $\tilde\Gamma$ with these
properties contains the canonical obstruction $\Gamma$. Since curves of
$\Gamma$ and $\tilde\Gamma$ do not intersect according to Theorem~\ref{Tcan},
it remains to show that no periodic component of $\hat\C\setminus\tilde\Gamma$
contains a curve of $\Gamma$. 
When the first-return map is a Thurston map, this follows from similar arguments
as above. When it is a homeomorphism, there would be a L\'evy cycle
intersecting $\Gamma$ within the component otherwise. \mybox

The assumption that $\tilde\Gamma$ is a simple obstruction is necessary,
because otherwise curves from $\tilde\Gamma$ and $\Gamma$ might intersect.
Consider the example in Figure~8 of \cite{dh}, where the spider map of angle
$5/12$ is mated with its conjugate. Various obstructions are formed by the
curves $\alpha,\,\beta$, which surround $2$-cycles corresponding to fixed
points of $z^2+\gamma_\sM(5/12)$, and by
$\delta_1\,,\,\delta_2\,,\,\delta_3\,,\,\delta_4$\,, which surround conjugate
postcritical points. Denoting the equator by $\eps$,
$\tilde\Gamma=\{\alpha,\,\beta,\,\eps\}$ would be a non-simple obstruction
with unobstructed component maps, and it does not contain the canonical
obstruction $\Gamma=\{\delta_1\,,\,\delta_2\}$.

Eigenvalues of non-negative integer matrices appear in two different
areas of Thurston's work: the combinatorial characterization of rational maps,
and core entropy of quadratic polynomials. $h(c)$ is the topological entropy
of $z^2+c$ on its Hubbard tree $T_c\subset\K_c$\,. It is computed from the
growth factor of iterated preimages, which in the postcritically finite case
is the leading eigenvalue $1\le\lambda\le2$ of the Markov matrix $A$ describing
the mapping of edges under the polynomial. Moreover, it is related to the
Hausdorff dimension of biaccessing angles. See \cite{core} and the references
therein.

\begin{prop}[Mating conjugate polynomials] \label{Pmcent}
Suppose $p\neq0$ is postcritically finite and take the complex conjugate
parameter $q=\overline{p}$. Consider the canonical obstruction $\Gamma$
of the formal mating $g=P\fmate Q$. Then:\\[1mm]
$1$. Each $\gamma\in\Gamma$ passes through a unique edge of the Hubbard
tree $\phi_0(T_p)$ and through the corresponding edge of $\phi_\infty(T_q)$\,;
there is a unique $\gamma\in\Gamma$ for each edge.\\[1mm]
$2$. The Markov matrix $A$ of $T_p$ is the transpose of the Thurston matrix
$M=M_\Gamma$ of $g$, unless different conventions are used in the preperiodic
case, with the critical points marked in the Hubbard tree but not in the formal
mating: then $A$ has an additional eigenvalue of $0$ compared to $M$. The
leading eigenvalue $\lambda$ is equal in any case,
so $h(p)=\log\lambda_\Gamma$\,. 
\end{prop}

\textbf{Proof} of Proposition~\ref{Pmcent}:
The Hubbard tree $T_p$ is a finite tree with expanding
dynamics; critical and postcritical points are marked and there may be
additional branch points \cite{bsa}. Let us assume that the critical
points are marked as well in the formal mating $g=P\fmate Q$. Marking $1$ on
the $0$-ray will not change the canonical obstruction. Now for
each edge of $T_p$ choose an arc in the dynamic plane of $P$ passing
homotopically transversely through this edge; join it with the complex
conjugate arc in the dynamic plane of $Q$ to obtain a simple loop in $\hat\C$
separating marked points of $g$. This works, since the formal mating is
constructed by mapping each dynamic plane to a half-sphere. The curves in
Figure~\ref{Fkk} are dynamic rays in fact; this choice is possible unless $p$
is a direct satellite center, but it is not required here.

1. These loops define a non-empty multicurve $\Gamma$. Under $P$, each edge is
covered by either one or two edges, so $\Gamma$ is invariant: in the former
case, one of the two preimages of the corresponding loop is inessential.
$\Gamma$ is completely invariant in fact, since every edge covers at least
one edge. By construction, the Markov matrix $A$ of $T_p$ and the Thurston
matrix $M_\Gamma$ agree except for transposition. Now $M_\Gamma$ contains at
least one nonzero entry in each row and in each column, so $\Gamma$ is a simple
obstruction.

Looking at noded Riemann surfaces in $\M_{G\cdot\Gamma}$ or at the
noded topological surface defining the pullback map on $\S_\Gamma$\,, the nodes
correspond to edges of $T_p$ and the pieces correspond to the vertices of
$T_p$\,, so they contain a single marked point or branch point of $T_p$ and the
corresponding point of $T_q$\,. In the former case, there are two
marked points of $g$ and at least one node, while in the latter case, there is
no marked point of $g$ and at least three nodes.

In the preperiodic case, all first-return maps are homeomorphisms. In the
periodic case, the first-return maps for marked points of $T_p$ map a sphere
with three or four marked points to itself. In the latter case the two critical
points are fixed and another point is fixed as well, while the fourth one goes
to that fixed point. This map is unobstructed by the core arc argument: an
obstructing curve cannot separate the two critical points, because its preimage
would cover it twice, giving an eigenvalue $1/2<1$. So an obstruction must
surround an arc between the critical points, and its preimages are two curves
separating the critical points from one of the other points each, so they are
not homotopic to the original curve.

By the Characterization Theorem~\ref{Tchar}, $\Gamma$ contains the canonical
obstruction. It remains to show that no proper subset has the same
properties. Assuming that $M_\Gamma$ and thus $A$ is reducible in the sense of
Perron--Frobenius, it has some block structure understood in terms of simple
renormalization according to \cite{core}. For each possible simple
renormalization $p=p'*\hat p$, so $p$ belongs to a small Mandelbrot set
$p'*\M$ \cite{book1mcm}, we have a block-triangular structure of $A$ with two
diagonal blocks, one corresponding to the periodic parameter $p'$, and an
imprimitive one related to the renormalized parameter $\hat p$. Now $A$
maps edges from the big Julia set of $p'$ over edges from the small Julia set
of $\hat p$, and $M_\Gamma$ maps small loops to big loops. So taking only the
small loops would not give an invariant multicurve, and taking only the big
loops would give obstructed component maps: the first-return map
corresponds to a mating of conjugate polynomials again. Conversely, although
not every possible block decomposition is explained in terms of
renormalization, every edge of $T_p$ belongs to a particular level of
renormalization, so removing the loop for this edge amounts to removing a block
from a particular form with two diagonal blocks.

2. By definition, $\lambda_\Gamma$ is the leading eigenvalue of $M_\Gamma$
and $h(p)=\log\lambda$ with the leading eigenvalue $\lambda$ of $A$.
In the preperiodic case, we need not mark the critical points for $g$, and we
might also not mark it in $T_p$\,; then we still have $A$ as the transpose of
$M_\Gamma$\,, since the edge around $0$ is mapped two--to--one to the edge
before $p$ and the corresponding loop has two homotopic preimages. The
matrices will be different using different conventions for $T_p$ and $g$, but
according to \cite{core} the leading eigenvalue is the same.
\mybox

Any multicurve corresponds to a tree in a similar way. A general estimate of
$\lambda_\Gamma$ in terms of entropy is given by Shishikura according to
\cite{st3}. Further applications include the description of
Herman rings \cite{shr}, the construction of maps with Cantor families
of circles in the Julia set \cite{gcc}, and the classification
of rescaling limits \cite{arl}. In recent talks on tropical complex dynamics,
Shishikura has suggested to combine arithmetic surgery on a tree with
quasiconformal surgery on pieces.
%
%

The canonical decomposition of a Thurston map is done in two steps by
Bartholdi--Dudko \cite{bd0}: L\'evy cycles generate a
decomposition  into pieces, such that the first-return maps are homeomorphisms,
expanding, or of type \2. The expanding pieces are decomposed again, such that
the first-return maps are equivalent to rational ones. Here expansion refers to
a suitable metric, which is defined everywhere except at super-attracting
cycles.  See \cite{bookbm, bookhp} for various notions of expansion.

\subsection{The Selinger proof of the Pilgrim Conjecture}\label{3p}
For a Thurston map $g$ with canonical obstruction $\Gamma\neq\emptyset$, the
first-return maps of $g_\Gamma$ are homeomorphisms or Thurston maps. In
\cite{kpcan}, Pilgrim has conjectured that the latter are unobstructed when
the orbifold is hyperbolic. This was proved by Selinger in Theorem~10.3
of \cite{ext} by constructing a subsequence of $(\tau_n)$ with suitable
convergence properties. He obtained the more general Characterization
Theorem~5.6 of \cite{char} using different techniques; see Theorem~\ref{Tchar}
in Section~\ref{33}. We will need convergence properties from the first proof
in Section~\ref{3s}. Maybe the following proof shall be read as a complement to
the original one in \cite{ext}, since it focuses on a few points that took me
some time to understand, in particular the construction of the uniform bound
$D_1$ and dealing with the fact that $\pi_1:\hat\T\to\hat\W$ and
$\pi_2:\hat\W\to\hat\M$ are not covering maps.

The Hurwitz space is defined as $\W=\T/H$, where $H<G$ denotes the subgroup
of liftable homeomorphisms: $h\in H$ if there is an $h'\in G$ with
$g\circ h'=h\circ g$. Then the covering $\pi:\T\to\M$ factorizes as
$\pi=\pi_2\circ\pi_1$ with $\pi_1:\T\to\W$ and $\pi_2:\W\to\M$. Moreover,
$H$ has finite index in $G$ and $\pi_2$ is a finite cover. Since $\W$ can be
represented by triples of rational maps $f_\tau$ and marked points in its
domain and range \cite{dh, book2h, endo, kps}, there is a continuous map
$\tilde\sigma_g:\W\to\M$ with $\tilde\sigma_g\circ\pi_1=\pi\circ\sigma_g$\,.
Now both $G$ and $H$ act by WP-isometries on $\hat\T$, and we have a
completion $\hat\W=\hat\T/H$ of Hurwitz space, which is compact. While
$\pi,\,\pi_1\,\,\pi_2$ are covering maps and local isometries, the extensions
$\pi_1:\hat\T\to\hat\W$ and $\pi_2:\hat\W\to\hat\M$ are weak contractions;
$\pi_1$ is ``infinitely branched'' and $\pi_2$ is finitely branched. The unique
extension $\tilde\sigma_g:\hat\W\to\hat\M$ is Lipschitz continuous with factor
$\sqrt{d}$. The strata of $\hat\W$ are labeled by classes $H\cdot\Gamma$ of
multicurves. We will not need a concrete description of $\hat\W$ and
$\tilde\sigma_g$ by triples.

\begin{prop}[Selinger proof of the Pilgrim Conjecture]\label{Pselprfpc}
Suppose $g$ is a Thurston map of degree $d$ with canonical obstruction
$\Gamma\neq\emptyset$, $C$ is a piece of $\hat\C/\Gamma$ mapped to
itself by $g_\Gamma$\,, and the component map $g^\sC:C\to C$ of degree $\ge2$
does not have orbifold type \2. Fix $\tau_0\in\T$ and set
$\tau_n=\sigma_g^n(\tau_0)$. Then:\\[1mm]
$1$. There is sequence $w_i$ in a compact subset of
$S_{H\cdot\Gamma}\subset\hat\W$ with $\tilde\sigma_g(w_i)=\pi_2(w_{i+1})$ for
$i\ge0$, and a subsequence $\tau_{n_k}$ with
$\pi_1(\sigma_g^i(\tau_{n_k}))\to w_i$ as $k\to\infty$, for all $i\ge0$.\\[1mm]
$2$. Assuming that $\eps(I)$ is sufficiently small and $\eps(I)\searrow0$
sufficiently fast, there are $k(I)$ and
$\hat\tau_I\in\pi_1^{-1}(w_0)\cap S_\Gamma$ with
$d_\sWP(\sigma_g^i(\tau_{n_{k(I)}}),\,\sigma_g^i(\hat\tau_I))<\eps(I)$
for $0\le i\le I$, such that $\pi_1(\sigma_g^i(\hat\tau_I))=w_i$ for
$0\le i\le I$.\\[1mm]
$3$. The component $\sigma_\sC=\sigma_{g^\sC}$ of the extended pullback map has
a unique fixed point $\hat\tau^\sC$ and
$d_\sT^\sC(\sigma_\sC^i(\hat\tau_I^\sC),\,\hat\tau^\sC)\to0$ as $i\to\infty$,
uniformly in $I\ge i$.
\end{prop}

Intuitively, the situation is as follows: imagine a system of coordinates in
a neighborhood of $\S_\Gamma\subset\hat\T$ adapted to a product of three
components. The first coordinate becomes $0$ as $\Gamma$ is pinched. The second
coordinate is related to pieces, where the first-return map is a homeomorphism,
and where we do not have convergence. The third coordinate describes pieces
where the first-return maps converge. Although such a local product
representation of $\hat\T$ is constructed in \cite{mpat}, we do not have any
estimates of $\sigma_g$ in that representation. So the proof will be given
by constructing various subsequences in an interplay between $\T$, $\hat\T$,
and components of $\S_\Gamma$\,, using both the Teichm\"uller metric and the
Weil--Petersson metric at times.
--- Note that we have an accumulation statement instead of convergence for two
reasons: there may be pieces with a least four marked points and nodes permuted
by a homeomorphism, and $\sigma_g$ is not weakly contracting on $\hat\T$.

\textbf{Proof:} 1. By the Canonical Obstruction Theorem~\ref{Tcan},
$\pi(\tau_n)$ accumulates on a compact subset of
$\S_{G\cdot\Gamma}\subset\hat\M$; pick an accumulation point $m_0$\,. Now
$\pi_2:\hat\W\to\hat\M$ is a finite branched cover, so there is a subsequence
$(\tau_{n_k^0})_{k\in\N_0}\subset(\tau_n)_{n\in\N_0}$ and a
$w_0\in\pi_2^{-1}(m_0)$ with $\pi_1(\tau_{n_k^0})\to w_0$ as $k\to\infty$; we
have $w_0\in\S_{H\cdot\Gamma}\subset\hat\W$ since precisely the curves in
$\Gamma$ are pinched as $n_k^0\to\infty$. Then
$\pi(\sigma_g(\tau_{n_k^0}))\to m_1=\tilde\sigma_g(w_0)$ by continuity.
Choose a subsequence $\tau_{n_k^1}$ of $\tau_{n_k^0}$ with
$\pi_1(\sigma_g(\tau_{n_k^1}))\to w_1\in\pi_2^{-1}(m_1)$. Define
$m_i$\,, $w_i$\,, and subsequences $\tau_{n_k^i}$ inductively, then any
subsequence $\tau_{n_k}$ of the diagonal sequence $\tau_{n_k^k}$ satisfies the
claim.

To obtain the bound $D_1$ below, this subsequence is constructed as follows:
For $\tau$ on a smooth curve from $\tau_0$ to $\tau_1$\,, there is a lower
bound $l(\gamma,\,\sigma_g^n(\tau))\ge R$, $\gamma\notin\Gamma$. Take constants
$C_*$ for $\T=\T_g$ according to item~5 of Proposition~\ref{Pbptms} and
$D_*=D_*(R)$ for $\T^\sC=\T_{g^\sC}$ according to item~6. Pick intermediate
points $\tau_0^0=\tau_0\,,\,\tau_0^1\,,\dots,\,\tau_0^U=\tau_1$ on the curve
with $d_\sT(\tau_0^{u-1}\,,\,\tau_0^u)\le D_*/C_*$ for $1\le u\le U$. Now
choose the indices $n_k$ for the subsequence $\tau_{n_k}$ of $\tau_{n_k^k}$
such that there are limits $\pi_1(\sigma_g^{n_k}(\tau_0^u))\to w_0^u$ as
$k\to\infty$.

2. Intuitively, if $\mathcal{N}$ is a small neighborhood of $w_0$ in $\hat\W$,
we have $\pi_1(\tau_{n_k})\in\mathcal{N}$ for sufficiently large $k$, and we
shall find $\hat\tau_0\in\pi_1^{-1}(w_0)$ close to $\tau_{n_k}$ for some $k$.
Then the curves of $\tilde\Gamma\in G\cdot\Gamma$ are short in a neighborhood
of $\hat\tau_0$ and only the curves in $\Gamma$ are short at $\tau_{n_k}$\,, so
$\tilde\Gamma=\Gamma$ and $\hat\tau_0\in\S_\Gamma$\,. A subgroup $G_\Gamma<G$
acts transitively on the fiber $\pi_1^{-1}(w_0)\cap\S_\Gamma$\,; it leaves
$\Gamma$ invariant and it is generated by Dehn twists about curves in pieces of
$\hat\C/\Gamma$ and about curves of $\Gamma$. The latter act trivially on the
stratum but not trivially in a neighborhood, thus $\pi_1$ is not a covering.
See also Remark~\ref{Raug}. In a neighborhood of $\hat\tau_0$\,, $\pi_1$ maps
to the quotient with respect to $H\cap G_\Gamma$\,, $\hat\T$ is a product
of disks and of half-disks plus the center point, and $\pi_1$ is
infinite--to--one locally on the half-disks at those points. So $\pi_1$ is not
a local WP-isometry at $\hat\tau_0$\,: we have
$d_\sWP(\pi_1(\tau'),\,\pi_1(\tau))=\min d_\sWP(h\cdot\tau',\,\tau)\le
 d_\sWP(\tau',\,\tau)$, and arbitrarily close to $\hat\tau_0$ there are
$\tau',\,\tau\in\T$ with
$d_\sWP(\pi_1(\tau'),\,\pi_1(\tau))<d_\sWP(\tau',\,\tau)$.
But $\hat\T$ is a unique geodesic space and the geodesic from $\tau_{n_k}$ to
$\hat\tau_0$ is contained in $\T$ except for the endpoint \cite{wwp}. All
$h\cdot\tau_{n_k}$ close to $\hat\tau_0$ are related by global isometries
fixing $\hat\tau_0$\,, thus
\be \label{eqsamewp}  d_\sWP(\tau_{n_k}\,,\,\hat\tau_0)
=d_\sWP(\pi_1(\tau_{n_k}),\,w_0)=d_\sWP(\pi(\tau_{n_k}),\,m_0) \ . \ee
So we may define $\mathcal{N}$ in terms of a small distance $\eps(0)$ to $w_0$
and have the same radius in components of $\pi_1^{-1}(\mathcal{N})$. Before
stating the actual construction of $\hat\tau_0$\,, note that we want to have
a shadowing property for a finite number $I$ of steps; this is possible
since $\sigma_g$ is Lipschitz continuous, and both $k(I)$ and $\hat\tau_I$ will
depend on $I$.

Assume that $\eps(I)\searrow0$ for $0\le I\to\infty$ and $(\sqrt{d}+1)\eps(I)$
is less than the minimal distance in the fiber 
$\pi_2^{-1}(m_i)\cap\S_{H\cdot\Gamma}$ for $0\le i\le I$. Moreover, the
preimage of an $\eps(I)$-neighborhood of $w_i$ under $\pi_1$ shall have
disjoint components, where
$\pi_1$ is described explicitly as an infinite--to--one map in terms of
$H\cap G_\Gamma$ as explained above. Then for $I\ge0$ there are $k(I)$ and
$\hat\tau_I\in\pi_1^{-1}(w_0)\cap S_\Gamma$ with
$d_\sWP(\sigma_g^i(\tau_{n_{k(I)}}),\,\sigma_g^i(\hat\tau_I))<\eps(I)$
for $0\le i\le I$. A finite induction shows
$\pi_1(\sigma_g^i(\hat\tau_I))=w_i$ for $0\le i\le I$: 
assuming the claim for $i-1$, we have
$\pi(\sigma_g^i(\hat\tau_I))=\tilde\sigma_g(\pi_1(\sigma_g^{i-1}(\hat\tau_I)))
=\tilde\sigma_g(w_{i-1})=m_i$ and 
\ban
& & d_\sWP(w_i\,,\,\pi_1(\sigma_g^i(\hat\tau_I))) \nonumber \\[1mm]
&\le & d_\sWP(w_i\,,\,\pi_1(\sigma_g^i(\tau_{n_k}))) +
d_\sWP(\pi_1(\sigma_g^i(\tau_{n_k}))\,,\,\pi_1(\sigma_g^i(\hat\tau_I)))
\nonumber \\[1mm]
& \le & d_\sWP(m_i\,,\,\pi(\sigma_g^i(\tau_{n_k}))) + \label{eqrevest}
d_\sWP(\sigma_g^i(\tau_{n_k})\,,\,\sigma_g^i(\hat\tau_I)) \\[1mm]
& \le & d_\sWP(\tilde\sigma_g(\pi_1(\sigma_g^{i-1}(\hat\tau_I)))\,,\,
\tilde\sigma_g(\pi_1(\sigma_g^{i-1}(\tau_{n_k})))) +
d_\sWP(\sigma_g^i(\tau_{n_k})\,,\,\sigma_g^i(\hat\tau_I)) \nonumber \\[1mm]
& \le & (\sqrt{d}+1)\,\eps(I) \ . \nonumber
\ean
Concerning the first term in (\ref{eqrevest}), we have discarded $\pi_2^{-1}$,
which is not a weak contraction in general. But in this case
$d_\sWP(w_i\,,\,\pi_1(\sigma_g^i(\tau_{n_k})))
 =d_\sWP(m_i\,,\,\pi(\sigma_g^i(\tau_{n_k})))$
by the same arguments as for (\ref{eqsamewp}). Finally,  by the assumption on
the distance in the fiber $\pi_2^{-1}(m_i)$, the claim is proved for $i$.

3. All maps, sets, and elements related to the stratum $\S_\Gamma\in\hat\T$
have a  product structure describing pieces of the noded Riemann
surfaces; the component for the piece corresponding to $C$ is denoted by a
superscript or subscript $C$. Since the length of hyperbolic geodesics is
continuous on $\hat\T$, we have $l(\gamma,\,\sigma_\sC^i(\hat\tau_I^\sC))\ge R$
for $0\le i\le I$ and all essential simple closed curves $\gamma$ in $C$. The
same lower bound holds for simple closed geodesics in the corresponding piece
of the noded surface defined by $\pi_2(w_0^u)$. Now $\sigma_g$ is weakly
contracting on $\T$ with respect to $d_\sT$\,, so the intermediate points
satisfy
$d_\sT(\sigma_g^{n_k}(\tau_0^{u-1}),\,\sigma_g^{n_k}(\tau_o^u))\le D_*/C_*$ and
$d_\sWP(\sigma_g^{n_k}(\tau_0^{u-1}),\,\sigma_g^{n_k}(\tau_o^u))\le D_*$\,.
For $I\ge I_*$ there are $\hat\tau_I^u\in\pi_1^{-1}(w_0^u)\cap S_\Gamma$ with
$d_\sWP(\sigma_g^{n_{k(I)}}(\tau_0^u),\,\hat\tau_I^u)<\eps(I)\le D_*/2$. Then
$d_\sWP(\hat\tau_I^{u-1},\,\hat\tau_I^u)<2D_*$ and the same estimate holds for
the components related to $C$. Since these components are intermediate points
between $\hat\tau_I^\sC$ and $\sigma_\sC(\hat\tau_I^\sC)$ with geodesic length
$\ge R$, item~6 of Proposition~\ref{Pbptms} gives
\be
d_\sT(\hat\tau_I^\sC,\,\sigma_\sC(\hat\tau_I^\sC))<2U\cdot1/4=U/2\le D_1 \ .\ee
Here $D_1\ge U/2$ is chosen such that the estimate holds not only for
$I\ge I_*$ but for all $I\ge1$. Now suppose $0\le i<I$. The pullback map
$\sigma_\sC$ is weakly contracting, so $\sigma_\sC^i(\hat\tau_I^\sC)$ and
$\sigma_\sC^{i+1}(\hat\tau_I^\sC)$ are connected with an arc of T-length
$\le D_1$\,, on which the length of simple closed geodesics is bounded below
by $R\cdot\e{-D_1}$. This condition defines compact subsets of $\M_\sC$ and
$\W_\sC$ according to the Mumford Proposition~\ref{Pbptms}.7. The
contraction of $\sigma_\sC$ at $\tau^\sC$ depends only on $f_{\tau^\sC}$ or
$\pi_1^\sC(\tau^\sC)$, so some iterate of $\sigma_\sC$ is uniformly
strictly contracting over the compact subset of $\W_\sC$\,, since $g^\sC$ is
not of type \2. For notational convenience, let us assume that the first
iterate suffices. So there is $0<L<1$ with
\be
d_\sT(\sigma_\sC^i(\hat\tau_I^\sC),\,\sigma_\sC^{i+1}(\hat\tau_I^\sC))
\le D_1\cdot L^i \qquad\mbox{and}\qquad
d_\sT(w_i^\sC\,,\,w_{i+1}^\sC) \le D_1\cdot L^i \ee
for $0\le i<I$; the second estimate is true for
all $i\ge0$ since $w_i$ is independent of $I>i$. By completeness of $\W^\sC$
we have limits $w_i^\sC\to\hat w^\sC$ and $m_i^\sC\to\hat m^\sC$ with
$\tilde\sigma_\sC(\hat w^\sC)=\pi_2^\sC(\hat w^\sC)=\hat m^\sC$ and an
estimate
\be \label{eqestfromsum}
d_\sT(w_i^\sC\,,\,\hat w^\sC)\le\frac{D_1}{1-L}\,L^i \ . \ee
Take a lower bound $\hat\eps$ for the T-distance in the fiber
$(\pi_1^\sC)^{-1}(\hat w^\sC)$, and such that an $\hat\eps/3$-neighborhood of
$\hat w^\sC$ has disjoint preimages under the covering $\pi_1^\sC$\,. Choose
$i=i_*$ such that the right hand side of (\ref{eqestfromsum}) is less than
$\hat\eps/3$. Pick $I>i$ and find $\hat\tau^\sC\in(\pi_1^\sC)^{-1}(\hat w^\sC)$
with $d_\sT(\sigma_\sC^i(\hat\tau_I^\sC),\,\hat\tau^\sC)<\eps/3$. Then
\be
d_\sT(\hat\tau^\sC,\,\sigma_\sC(\hat\tau^\sC)) \le
d_\sT(\hat\tau^\sC,\,\sigma_\sC^i(\hat\tau_I^\sC)) +
d_\sT(\sigma_\sC^i(\hat\tau_I^\sC),\,\sigma_\sC^{i+1}(\hat\tau_I^\sC)) +
d_\sT(\sigma_\sC^{i+1}(\hat\tau_I^\sC),\,\sigma_\sC(\hat\tau^\sC)) \ , \ee
which is bounded by $\hat\eps/3+\hat\eps/3+\hat\eps/3$. So $\hat\tau^\sC$ is a
fixed point of $\sigma_\sC$\,. Since the orbifold is not of type \2, $g^\sC$ is
unobstructed, the fixed point is unique, and with $D=D_1/(1-L)$ we have
$d_\sT(\sigma_\sC^i(\hat\tau_I^\sC),\,\hat\tau^\sC)\le D\cdot L^i$ for all
$i_*\le i<I<\infty$. Actually, an estimate of this form is valid as well, when
only some iterate of $\sigma_\sC$ is strictly contracting. \mybox

\subsection{Essential equivalence and convergence properties} \label{3s}
Suppose $g$ is a Thurston map with canonical obstruction $\Gamma\neq\emptyset$.
Considering the Thurston Algorithm $\tau_n=\sigma_g^n(\tau_0)$ with rational
maps $f_n$ sending point configurations at time $n+1$ to those at time $n$, we
know that curves corresponding to $\Gamma$ will be pinched and points will
collide. Normalizing three points to $\infty,\,0,\,1$, a component $C$ of
$\hat\C\setminus\Gamma$ is singled out, and all other components will shrink
to points. We shall see that in a fairly general situation, these points do
not wander, but they have limits as $n\to\infty$. So there is a limiting point
configuration, such that not all points are distinct, and a rational map $f$
with $f_n\to f$. The following Theorem~\ref{Tselconv} has applications to
quadratic matings, anti-matings, precaptures, and spiders. Its proof is based
on the extension of $\sigma_g$ to augmented Teichm\"uller space $\hat\T$; but I
have tried to formulate the assumptions in terms of components of
$\hat\C\setminus\Gamma$ instead of pieces of $\hat\C/\Gamma$, so that it may
be applied in other papers without introducing $\hat\T$.

\begin{dfn}[Essential equivalence]\label{Dselconv}
Consider a bicritical Thurston map $g:S^2\to S^2$ or $g:\hat\C\to\hat\C$ of
degree $d\ge2$, with marked set $Z$ including the critical points. Suppose
there is a completely invariant multicurve $\Gamma\neq\emptyset$ and a
component $C$ of $\hat\C\setminus\Gamma$ such that:\\[1mm]
$\bullet$ All components $\tilde C\neq C$ of $\hat\C\setminus\Gamma$ are disks;
these disks are preperiodic or periodic under $g$ after an isotopic deformation
$\phi$, and the periodic disks are mapped homeomorphically.\\[1mm]
$\bullet$ $C$ is mapped to itself with degree $d$ in the sense explained
below; the essential map $\tilde g$ is equivalent to a rational map $f$.\\[1mm]
Then we shall say that $g$ is \textbf{essentially equivalent} to $f$.
\end{dfn}

Denoting a union of concrete curves by the same symbol $\Gamma$, the
full preimage is $g^{-1}(\Gamma)=\Gamma'\cup\Gamma''$, where $\Gamma'$
consists of essential curves ambient isotopic to $\Gamma$ and the
curves in $\Gamma''$ are inessential. 
Choose a homeomorphism $\phi:\hat\C\to\hat\C$ with $\phi:\Gamma\to\Gamma'$,
isotopic to the identity relative to the marked set $Z$, then
$g\circ\phi:\Gamma\to\Gamma$. We have assumed that $\Gamma$ is not nested:
there is a distinguished component $C$ of $\hat\C\setminus\Gamma$, such that
all of the other components are disks, and mapped to disks by $g\circ\phi$.
Intuitively, the construction of noded topological surfaces is done by
pinching all curves in $\Gamma$, respectively $\Gamma'\cup\Gamma''$, to points
and by chopping off spheres with at most one marked point; suitable
identifications are needed to obtain a self-map $g_\Gamma$\,.

In this situation, the definition of the component map $g^\sC$ according to
\cite{ext} and Section~\ref{32} can be paraphrased
as follows: $C'=(g\circ\phi)^{-1}(C)\subset C$ is a component of
$\hat\C\setminus(\Gamma\cup\phi^{-1}(\Gamma''))$. Considering each boundary
curve as a single point, $g\circ\phi:C'\to C$ defines branched coverings
$S^2\to S^2$ or $\hat\C\to\hat\C$ up to isotopy. More precisely, choose a
continuous $\psi:\hat\C\to\hat\C$, which is sending closed disk components to
points and is injective otherwise. Modify $\psi$ to $\psi'$ in disjoint
neighborhoods of the additional disk components of
$\hat\C\setminus\phi^{-1}(\Gamma'')$, sending these to points as well, without
moving possible single marked points in these disks. Then the essential map
$\tilde g$ with $\psi\circ g\circ\phi=\tilde g\circ\psi'$ is a Thurston
map, whose marked points correspond to marked points of $g$ in $C'$, disk
components of $\hat\C\setminus\Gamma$, or to disks bounded by peripheral curves
in $\phi^{-1}(\Gamma'')$. It is defined uniquely up to combinatorial
equivalence and represents $g^\sC$.

\begin{thm}[Convergence of marked points and rational maps]\label{Tselconv}
Consider a bicritical Thurston map $g$, a multicurve $\Gamma\neq\emptyset$ and
a component $C$ of $\hat\C\setminus\Gamma$, such that $g$ is essentially
equivalent to a rational map $f$ according to Definition~$\ref{Dselconv}$. Use
a normalization of critical points at $0$ and $\infty$, and another marked
point at $1$, which is arbitrary in $C$ or in a disk not containing a critical
point. Normalize $f$ analogously. Then the Thurston Algorithm $\sigma_g$ for
the unmodified map $g$ with any initial $\tau_0\in\T$ satisfies:\\[1mm]
$1$. Precisely the curves of $\Gamma$ are pinched, and
$\tau_n=\sigma_g^n(\tau_0)$ diverges in $\T$.\\[1mm]
$2$. If $f$ is not of type \2, we have $f_n\to f$. The marked points
converge to preperiodic and periodic points of $f$; two points collide if and
only if they belong to the same disk $\tilde C\neq C$ in
$\hat\C\setminus\Gamma$.\\[1mm]
Analogous statements hold for a path $\tau_t$ with
$\tau_{t+1}=\sigma_g(\tau_t)$.
\end{thm}

The notion of essential equivalence is meant to indicate a generalization of
Thurston's combinatorial equivalence, but in this form it is not an equivalence
relation.
This notion shall emphasize that $g$ itself determines the canonical
obstruction $\Gamma$, the essential map $\tilde g$, and the rational map $f$,
and no modification is needed to ensure $f_n\to f$.

A typical example is provided by a formal mating $g=P\fmate Q$ of quadratic
polynomials, having ray connections between postcritical points but no
cyclic ray connections: then $\Gamma$ consists of curves around postcritical
ray-equivalence classes. These are the only obstructions according to
\cite{rst}, so the essential mating $\tilde g$ is unobstructed and
combinatorially
equivalent to a rational map $f$, excluding type \2 for now. Then
Theorem~\ref{Tselconv} gives $f_n\to f$ for the pullback defined by the
unmodified formal mating; see Section~\ref{4c} for details. According to
\cite{pmate}, the convergence statement is wrong in general when $\tilde g$ and
$f$ are of type \2. --- Further applications are given in
\cite{emate, amate, smate}. For simplicity I have considered the
bicritical case only, because $f_n$ is determined explicitly by its critical
values; to obtain a more general result, it should also be checked whether
collisions between critical points are allowed.

\textbf{Proof of Theorem~\ref{Tselconv}.1:} According to the discussion of
Theorem~\ref{Text}, $C$ is considered as a piece of a noded surface and
is related to a component of the invariant stratum
$\S_{G\cdot\Gamma}\subset\hat\M$, and $\tilde g=g^\sC$ is a component map of
$g_\Gamma$\,. All disks correspond to pieces, which are attached directly to
$C$. If one of these pieces contains a critical point, the node will be
critical as well, and the piece is mapped with degree $d$; by assumption it is
preperiodic. Accordingly, a periodic critical point must belong to $C$. Both
critical points are allowed to be in pieces corresponding to disks
$\tilde C\neq C$, but not in
the same piece, because then $\tilde g$ would have degree one.

By assumption, $\Gamma$ consists of one or several L\'evy cycles and all of
their essential preimages, so it is a simple obstruction. The component map
$g^\sC=\tilde g$ is unobstructed, since it is combinatorially equivalent to the
rational map $f$, which is not a flexible Latt\`es map. 
The remaining first-return maps are homeomorphisms, so according to the
Characterization Theorem~\ref{Tchar}, $\Gamma$ contains the canonical
obstruction. Assuming that the canonical obstruction was smaller, we would
enlarge $C$ by omitting one or several L\'evy cycles from $\Gamma$, and $g^\sC$
would be obstructed. So $\Gamma$ is the canonical obstruction of $g$; by
Theorem~\ref{Tcan}, $(\tau_n)$ accumulates at most on $\S_\Gamma\subset\hat\T$
and $\pi(\tau_n)$ accumulates on a compact subset of $\S_{G\cdot\Gamma}$\,.
\mybox

Intuitively, what happens is that point configurations $x_i(n)$ cluster
according to the disks containing $z_i$\,, and the pullback of clusters
determined by $\sigma_g$ should stay close to the pullback of single points
obtained from $\sigma_f$\,. This argument involves interchanging limits, and I
have not been able to prove it with direct estimates. So, convergence will be
proved below in two steps: first use augmented Teichm\"uller
space and the Selinger Proposition~\ref{Pselprfpc} for a global result,
accumulation at the prospective limit $x^\infty$. Then a local result is
applied, attraction according to Proposition~\ref{Pselconv}. The global result
is needed to get close to $x^\infty$ and to distinguish different fixed points
of the extended pullback relation. And the local result is needed because the
global one provides accumulation only, not convergence, as explained in
Section~\ref{3p}.
The pullback of point configurations may be visualized as a kind of movie:
as time $t$ or $n$ flows, the points $x_i(n)$ move within a single copy of
$\hat\C$. To formulate neighborhoods, convergence, or derivatives of point
configurations, it is more convenient to consider
$x=(x_1\,,\,\dots\,,\,x_{|Z|})$ as a single point in $\hat\C^{|Z|}$.

\textbf{Notations and basic properties} of the pullback relation:
The branch portrait of $g$ defines a map $\#$ of indices, such that marked
points of $g$ are mapped as $g(z_i)=z_{i\#}$\,. The Thurston Algorithm provides
sequences of homeomorphisms $\psi_n$ and rational maps $f_n$ with
$f_n\circ\psi_{n+1}=\psi_n\circ g$. The point configuration $x(n)$ at time $n$
has the components $x_i(n)=\psi_n(z_i)$; sometimes these are called marked
points as well. The basic relation is $f_n(x_i(n+1))=x_{i\#}(n)$.

The coordinates depend on the normalization, and we have specific indices
$\alpha',\,\beta',\,\gamma'$ with $x_{\alpha'}(n)=\infty$, $x_{\beta'}(n)=0$,
and $x_{\gamma'}(n)=1$ for all $n$. Assume that $z_{\alpha'}$ and $z_{\beta'}$
are the critical points of $g$. The normalization fixes an embedding
$\pi_3:\M\to\hat\C^{|Z|}$; its range consists of $(x_1\,,\,\dots\,,\,x_{|Z|})$
with pairwise distinct components, and $\infty,\,0,\,1$ at specific
positions, which is open and dense in a subset $\C^{|Z|-3}\subset\hat\C^{|Z|}$.

Denote $\alpha=\alpha'\#,\,\beta=\beta'\#,\,\gamma=\gamma'\#$, then $f_n$ is
determined by $f_n(\infty)=x_\alpha(n)$, $f_n(0)=x_\beta(n)$, and
$f_n(1)=x_\gamma(n)$ as a M\"obius transformation of $z^d$. This gives the
pullback relation in the form of a multi-valued function from $\pi_3(\M)$ into
itself: 
\be\label{eqpbrmv}
f_n^{-1}(z)=\sqrt[{\scriptstyle d}\,]{
\frac{x_\gamma(n)-x_\alpha(n)}{x_\gamma(n)-x_\beta(n)}
\cdot\frac{z-x_\beta(n)}{z-x_\alpha(n)} } \qquad
x_i'=\sqrt[{\scriptstyle d}\,]{
\frac{x_\gamma-x_\alpha}{x_\gamma-x_\beta}
\cdot\frac{x_{i\#}-x_\beta}{x_{i\#}-x_\alpha} }
\ee
The second formula may be considered either as a multi-valued function
$x\mapsto x'$, or as a step of the pullback $x(n)\mapsto x(n+1)$. By the
normalization we have $|Z|-3$ independent variables in the domain and
$|Z|-3$ variables in the range; when some $x_m=\infty$, the corresponding
factors cancel from the radicand. For each value of the index $i$,
the radicand is either constant $\infty,\,0,\,1$ or never $\infty,\,0,\,1$. In
the Thurston pullback, a specific branch of the $d$-th root is chosen by
the isotopy class of $\psi_{n+1}$ or by continuity of the path $x_i(t)$. Note
that unless all marked points $z_i$ of $g$ are periodic, we have pairs of
indices $i\neq k$ with $g(z_i)=g(z_k)$, so $i\#=k\#$. Then $x_i(n+1)$ and
$x_k(n+1)$ are given by different branches of a root with the same radicand,
and the number of variables could be reduced by replacing $x_k$ with
$\zeta\cdot x_i$ for some $\zeta$ with $\zeta^d=1$ determined by $g$; while
this makes sense for a concrete example, it would complicate the
notation for the present discussion.

\textbf{Relating $g$ to $f$:}
To describe the Thurston pullback $\sigma_f$\,, denote the marked points
of $f$ by $\tilde z=(\tilde z_1\,,\,\dots\,,\,\tilde z_{|\tilde Z|})$ and
define indices with $\tilde z_{\tilde\alpha'}=\infty$,
$\tilde z_{\tilde\beta'}=0$, $\tilde z_{\tilde\gamma'}=1$, and
$\tilde z_{\tilde\alpha}=f(\infty)$, $\tilde z_{\tilde\beta}=f(0)$,
$\tilde z_{\tilde\gamma}=f(1)$. Each marked point of $f$ corresponds either to
a unique marked point of $g$ in $C$ or to a unique disk component
$\tilde C\neq C$ of $\hat\C\setminus\Gamma$ with at least two marked points;
this defines a surjection $D:\{1,\,\dots,\,|Z|\}\to\{1,\,\dots,\,|\tilde Z|\}$.
The normalizations are assumed to be compatible, i.e.,
$D(\alpha')=\tilde\alpha'$, $D(\beta')=\tilde\beta'$, and
$D(\gamma')=\tilde\gamma'$.

Denote a specific diagonal of $\hat\C^{|Z|}$ by $\Delta_\Gamma$\,: it contains
all $x=(x_1\,,\,\dots\,,\,x_{|Z|})$ with $x_i=x_k$ if and only if $D(i)=D(k)$,
and $x_{\alpha'}=\infty$, $x_{\beta'}=0$, $x_{\gamma'}=1$. The extension of
$\pi_3$ to $\hat\M$ satisfies $\pi_3(\S_{G\cdot\Gamma})=\Delta_\Gamma$\,; this
is an isomorphism when all disks have at most two marked points, but it is
forgetful otherwise. Note that there is a natural bijection between
$\Delta_\Gamma\subset\hat\C^{|Z|}$ and the configuration space of $\sigma_f$
contained in $\hat\C^{|\tilde Z|}$: the components of $x$ have repetitions,
such that $x$ corresponds to $\tilde x$ with $x_i=\tilde x_{D(i)}$ for all $i$.
In particular, the point configuration $x^\infty$ with repetitions given by
$x^\infty_i=\tilde z_{D(i)}$ is the prospective limit of $x(n)$.

The Thurston pullback $\sigma_f$ defines a multi-valued pullback relation on
its configuration space, which is given by formulas analogous to
(\ref{eqpbrmv}). A specific branch of this pullback relation has a fixed point
at $\tilde x=\tilde z$; it is analytic in a neighborhood of $\tilde z$ in
$\C^{|\tilde Z|-3}\subset\hat\C^{|\tilde Z|}$. A simple but suggestive
observation is the following: under the bijection of $x$ and $\tilde x$ the
conjugate pullback map in a neighborhood of
$x^\infty$ in $\Delta_\Gamma$ is given by choosing a suitable branch in
(\ref{eqpbrmv}). The reason is that when $z_i$ is in the disk corresponding to
$\tilde z_j$\,, then $g(z_i)$ is in the disk corresponding to $f(\tilde z_j)$.
And if $z_i$ and $z_k$ are in the same disk, their images are both in one disk.
So for $x\in\Delta_\Gamma$ the radicands with different indices $i\#$ and $k\#$
agree whenever $D(i)=D(k)$. Note that this does not mean that the local branch
extends to a neighborhood of $x^\infty$ in $\hat\C^{|Z|-3}$: this is the case
precisely when the radicand cannot become $0$ or $\infty$. If, e.g., there is
an index $k\neq\beta'$ with $D(k)=D(\beta')=\tilde\beta'$, so $x_k$ will be
identified with $0$, the radicand may be $0$ within any neighborhood of
$x^\infty$ in $\hat\C^{|Z|-3}$. (On $\Delta_\Gamma$\,, the radicand in
(\ref{eqpbrmv}) is constant $0$ both for $i=k$ and $i=\beta'$,
which is not a problem.)

\begin{prop}[Local attraction]\label{Pselconv}
Consider $g$, $\Gamma$, and $f$ according to Theorem~$\ref{Tselconv}$.
Assume in addition that the marked point normalized to $1$ is chosen more
restrictively: if a marked point is identified with a critical point, then $1$
represents a critical value, whose iterates are not identified with a critical
point. Using the notations introduced above we have:

$1$. If no marked point is identified with a critical point, a branch of the
pullback relation $(\ref{eqpbrmv})$ extends analytically to a neighborhood
of $x^\infty$ in $\C^{|Z|-3}$. The eigenvalues $\lambda$ of the derivative at
the fixed point $x^\infty$ are of the following form: they are eigenvalues of
$D\sigma_f$ at $[\mathrm{id}]$, or $\lambda=0$, or $\lambda^{rp}=\rho^{-r}$
when an $rp$-cycle of $g$ in a $p$-cycle of disks corresponds to a $p$-cycle
of $f$ having multiplier $\rho$.

$2$. If $f$ is not of type \2, there is a neighborhood $\mathcal{N}$ of
$x^\infty$ in $\hat\C^{|Z|-3}$, which is attracting in the following sense:
when $\tau_t$ is a path in $\T$ with $\tau_{t+1}=\sigma_g(\tau_t)$, and
$\pi_3(\pi(\tau_t))\in\mathcal{N}$ for a segment of $t$-length $1$, then the
path in configuration space will stay in $\mathcal{N}\cap\pi_3(\M)$
forever and converge to $x^\infty\in\mathcal{N}\setminus\pi_3(\M)$. 
\end{prop}

The additional assumption restricts the choice of the index
$\gamma'$ with $x_{\gamma'}(n)=\psi_n(z_{\gamma'})=1$, when a critical point
$\omega$ of $g$ is identified with another marked point, i.e., they are in the
same disk component $\tilde C\neq C$ of $\hat\C\setminus\Gamma$. Then $\omega$
is strictly preperiodic and all forward iterates of $\omega$ will belong to
disks as well; preimages of $\omega$ may be marked or not, and in the
former case, may undergo identifications or not. In general we may take the
critical value $g(\omega)$ for $z_{\gamma'}$\,, unless it equals the other
critical point or is identified with it or with a preimage. In that case
we can take the other critical value, which is not identified with a
preimage of $\omega$, since the critical points of $f$ are not periodic.

\textbf{Proof of Proposition~\ref{Pselconv}:} Now after choosing $\gamma'$ and
$\tilde\gamma'=D(\gamma')$, which determines the term of $f$, we shall consider
three cases of increasing complexity:\\[1mm]
Case 1: no marked point is identified with a critical point.\\[1mm]
Case 2: a marked point is identified with a critical point, but no postcritical
point is identified with a critical point.\\[1mm]
Case 3: for some $k\ge1$, $g^k(z_{\alpha'})$ is identified with $z_{\beta'}$\,.
This means that $f^k(\infty)=0$, but $g$ has disjoint critical orbits.\\[1mm]
When $g^k(z_{\beta'})$ is identified with $z_{\alpha'}$ instead, this does not
require separate arguments, since $f$ and $f_n$ are related to case~3 by a
conjugation with $z\mapsto1/z$.

\textbf{Case 1:} Recall that $\tilde z$ denotes the marked points $\tilde z_j$
of $f$ and $x^\infty\in\Delta_\Gamma$ with $x^\infty_i=\tilde z_{D(i)}$ is our
prospective limit of $x(n)$. Since the marked points of $f$ are preimages of
marked points under pullback with suitable branches of $f^{-1}$, for each
$i\neq\alpha',\,\beta'$ there is a unique branch in (\ref{eqpbrmv}), such that
taking components of $x=x^\infty$ for the radicand, gives $x_i'=x_i^\infty$ for
the root. These branches extend analytically to $x$ in a neighborhood of
$x^\infty$ in $\C^{|Z|-3}$, which means $|x_i-x_i^\infty|<\eps$ for
$i=1,\dots,\,|Z|$ with $i\neq\alpha',\,\beta',\,\gamma'$, since the radicands
are varying in small neighborhoods of values distinct from $0$ and $\infty$.

To determine the eigenvalues of the derivative matrix for this branch of the
pullback relation, we shall obtain block matrices by choosing new coordinates
labeled $u=(u_1\,,\,\dots\,,\,u_{|\tilde Z|})$ and
$v=(v_{|\tilde Z|+1}\,,\,\dots\,,\, v_{|Z|})$ as follows:\\[1mm]
$\bullet$
For each $j=1,\dots,\,|\tilde Z|$, choose one index $i$ with $D(i)=j$ and set
$x_i=u_j$\,. If there are $k\neq i$ with $D(k)=j$, set $x_k=u_j+v_m$ for an
unused index $m$. When all variables are defined successively, note that
$x=x^\infty$ corresponds to $u=\tilde z$ and $v=0$.\\[1mm]
$\bullet$
For $j=\tilde\gamma'$, the marked point at $1$ shall be $u_j$\,, i.e.,
$x_{\gamma'}=u_{\tilde\gamma'}$. The corresponding choice for $0$ and $\infty$
is satisfied here anyway, because these are not identified with other marked
points, but it is required explicitly in cases~2 and~3.\\[1mm]
$\bullet$
Renumber $v_{|\tilde Z|+1}\,,\,\dots v_{|Z|}$ such that preperiodic marked
points $u_j+v_m$ appear before periodic ones, higher preperiods before lower
preperiods, and the periodic marked points are grouped according to their
cycles, with a natural order within each cycle. We may renumber the components
of $x$ such that $x_j=u_j$ or $x_i=u_{D(i)}+v_i$\,, respectively.\\[1mm]
Since three components of $u$ are constant, the local branch
$(u,\,v)\mapsto(u',\,v')$ of the pullback relation still has
$(|\tilde Z|-3)+(|Z|-|\tilde Z|)=|Z|-3$ free variables. We shall see that the
derivative at the fixed point $(\tilde z,\,0)$ has the block-triangular form
\be\label{eqmatbl}
D=\left(\begin{array}{cc}A & C\\ 0 & B\end{array}\right)
\qquad\mbox{with}\qquad
B=\left(\begin{array}{cc}R & Q\\ 0 & P\end{array}\right) \ .
\ee
Note that the local manifold with $u\approx\tilde z$ and $v=0$ is
invariant under the pullback relation: $u_j'$ and $u_j'+v_i'$ correspond to
marked points of $g$ in the same disk, so their images belong to one disk as
well. If $v=0$, the radicands determining $u_j'$ and $u_j'+v_i'$ agree, and the
local branches of the roots agree, so $v_i'=0$. This argument shows that the
lower left block of $D$ is $0$, and the upper left block $A$ coincides with
the derivative of the local pullback relation for $\sigma_f$ in configuration
space, which is analytically conjugate to $\sigma_f$ in a neighborhood of its
fixed point in Teichm\"uller space. So $A$ has attracting eigenvalues $\lambda$
unless $f$ is of type \2 --- then one eigenvalue will be neutral.

Before discussing the block structure of $B$, let us consider partial
derivatives of the pullback relation (\ref{eqpbrmv}) more explicitly. A few
components of $(u,\,v)$ determine the bicritical rational map $f_{uv}$ by its
critical values and the image of $1$. For suitable combinations of indices we
have
\be\label{eqfvdv}
f_{uv}(u_j'+v_i')=u_l+v_k \qquad\mbox{and}\qquad
f_{uv}(u_j')=u_l+v_m \ee
in general; either $v_k$ or $v_m$ may be missing. Total differentials at
$(u,\,v)=(\tilde z,\,0)$ read
\be
\dots+f'(\tilde z_j)\cdot(\mathrm{d}u_j'+\mathrm{d}v_i')
=\mathrm{d}u_l+\mathrm{d}v_k \qquad\mbox{and}\qquad
\dots+f'(\tilde z_j)\cdot\mathrm{d}u_j'
=\mathrm{d}u_l+\mathrm{d}v_m \ , \ee
where $\dots$ denotes differentials involving the partial derivatives of
$f_{uv}(z)$ with respect to the parameters $(u,\,v)$. Observing that these
expressions agree for both equations when setting $u=u'=\tilde z$ and $v=v'=0$,
the difference gives
\be\label{eqpdvdv}
f'(\tilde z_j)\cdot\mathrm{d}v_i'=\mathrm{d}v_k-\mathrm{d}v_m \ . \ee
Note again that either $v_k$ or $-v_m$ may be missing. The argument remains
valid, if a component $v_i\,,\,v_k\,,\,v_m$ appears in the parameters of
$f_{uv}$ as well, or if $u_j'$ or $u_l$ is $1$.
Now (\ref{eqpdvdv}) shows again that the lower left block of $D$ is $0$. The
blocks $R$ and $P$ of $B$ in (\ref{eqmatbl}) refer to preperiodic and periodic
marked points, respectively. In (\ref{eqpdvdv}), $v_i'$ has a higher preperiod
than $v_k$ and $v_m$\,, or these refer to a cycle of disks. Thus $R$ is
strictly upper triangular, with $0$ on the diagonal, and its eigenvalues
$\lambda$ vanish.

Finally, consider the blocks of $P$, which are related to periodic cycles of
$g$. Since periodic disks of $\hat\C\setminus\Gamma$ are mapped
homeomorphically by $g\circ\phi$, each disk in a cycle has the same number of
marked points, which are permuted by the first-return map. The obvious examples
are two cycles of periodic points with the same period identified pairwise, or
a single cycle of high period identified such that a cycle of lower period
results for $f$. All possibilities are covered by the following description: a
cycle of $f$ has period $p\ge1$ and there are one or several cycles of $g$ with
periods $rp$ for possibly different values $r\ge1$. Consider two scenarios:

First, $g$ has a $p$-cycle in the disks under consideration, which is labeled
$u_1\,,\,\dots\,,\,u_p$\,. The indices are chosen to illustrate
the order in the $v$-blocks; they do not actually start with $1$. Every
other $rp$-cycle of $g$, $r\ge1$, is described by $v_i$ in the following order:
$u_1+v_1\,\dots\,u_p+v_p\,,\,u_1+v_{p+1}\,\dots\,u_p+v_{2p}\,,\,\dots\,,\,
u_1+v_{(r-1)p+1}\,\dots\,u_p+v_{rp}$\,.
Then $P$ contains an $rp$-block on the diagonal with nonzero entries
$1/f'(\tilde z_j)$ only directly above the diagonal and in the lower left
corner, since $\pm v_m$ is absent from (\ref{eqfvdv})--(\ref{eqpdvdv}). By
rescaling all variables, this block becomes a companion matrix and
the entry in the lower left position is
$(1/f'(\tilde z_1)\cdot\dots\cdot 1/f'(\tilde z_p))^r=\rho^{-r}$. So the
characteristic polynomial is $\lambda^{rp}-\rho^{-r}$. Note that
$|\rho|>1$ because $f$ is postcritically finite and all periodic points are
superattracting or repelling, so $|\lambda|<1$.

Second, all cycles of $g$ within the $p$-cycle of disks have periods $rp$ with
$r>1$, then the $u_j$ must belong to one of these cycles; choose
this cycle to be last among the cycles of $g$ in the current cycle of
components. Starting with index $1$ again for simplicity, it is labeled as
$u_1\,\dots\,u_p\,,\,u_1+v_1\,\dots\,u_p+v_p\,,\,\dots\,,\,
u_1+v_{(r-2)p+1}\,\dots\,u_p+v_{(r-1)p}$\,.
Now $f_{uv}(u_p')=u_1+v_1$ shows that $\pm v_m$ is no longer absent from
(\ref{eqfvdv})--(\ref{eqpdvdv}). We have a block of size $(r-1)p$ on the
diagonal of $P$, with nonzero entries $1/f'(\tilde z_j)$ directly above the
diagonal, and further entries $-1/f'(\tilde z_p)$ in rows
$p,\,2p,\,\dots\,,\,(r-1)p$ of the first column. After rescaling all variables
appropriately, this is the companion matrix of
$\lambda^{(r-1)p}+\rho^{-1}\lambda^{(r-2)p}+\dots+\rho^{-(r-1)}=
\frac{\lambda^{rp}-\rho^{-r}}{\lambda^p-\rho^{-1}}$.
Note that $A$ contains a cyclic $p$-block for $u_1\,,\,\dots\,,\,u_p$ again,
but since $A$ is not block-triangular, it need not have eigenvalues with
$\lambda^p=\rho^{-1}$. If $g$ has further $r'p$-cycles in the same components
of $\hat\C\setminus\Gamma$, these are treated according to the first scenario,
giving $\lambda^{r'p}=\rho^{-r'}$; there will be additional entries above the
diagonal blocks, which do not contribute to the characteristic polynomial of
$P$, but may prevent it from being diagonalizable. --- An alternative approach
to the second scenario would be to modify $g$ isotopically so that it has a
$p$-cycle in the disks, and to mark this cycle in addition.

So if $f$ is not of type \2, all eigenvalues of the extended pullback relation
$x\mapsto x'$ at the fixed point $x^\infty$ are attracting. We shall construct
a norm on $\C^{|Z|-3}$, such that the linearization satisfies
$\|x'-x^\infty\|\le L\|x-x^\infty\|$ for some $L<1$; in a small neighborhood
$\mathcal{N}$ with $\|x-x^\infty\|<\delta$ the pullback map is analytic and
satisfies $\|x'-x^\infty\|\le L'\|x-x^\infty\|$ for some $L<L'<1$. To define
this norm, conjugate the derivative matrix to its Jordan normal form by a
linear change of variables. Rescale components such that the entries $1$ above
the diagonal become $\eps$, and choose $\eps>0$ small such that the new matrix
is contracting with respect to the standard Euclidean norm. The norm
$\|\cdot\|$ in the original coordinates corresponds to this Euclidean norm. 

In the remaining cases~2 and~3, we will not have an analytic branch of
$x\mapsto x'$ in a neighborhood of $x^\infty$, but we shall construct an
attracting neighborhood for the path $x(t)$ nevertheless when $f$ is not of
type \2. In \textbf{case 2}, suppose that $z_m$ is identified with a critical
point $\omega$. If $g^k(z_{\alpha'})=z_{\beta'}$ and $\omega=z_{\alpha'}$, then
$g^k(z_m)$ is identified with $z_{\beta'}$\,, and we shall redefine
$\omega=z_{\beta'}$ and $m$ such that $z_m$ is identified with $\omega$.
Preimages of $z_m$ and $\omega$ may be marked or not, and identified or not.
Define new coordinates $(u,\,v,\,w)$ with $w$ representing all $x_i$\,, such
that $z_i$ is identified with a critical or precritical point, and $u,\,v$
describing the remaining $x_i$ as in case~1. The marked point $x_{\gamma'}$
normalized to $1$ may be precritical but not be identified with a precritical
or critical point. So the rational maps $f_{uv}$ will not depend on $x_m$ and
its preimages; the multi-valued pullback relation
$(u,\,v,\,w)\mapsto(u',\,v',\,w')$ is such that $u'$ and $v'$ do not depend on
$w$. As in case~1, we have a local analytic branch and an attracting
neighborhood $\mathcal{N}_0$ for $(u,\,v)\mapsto(u',\,v')$.

Now the pullback for $x_m$ is asymptotic to
$x_m'\sim\sqrt[{\scriptstyle d}\,]{c\cdot v_j}$ or
$x_m'\sim\sqrt[{\scriptstyle d}\,]{c/v_j}$ when $\omega=z_{\beta'}$ or
$\omega=z_{\alpha'}$\,,
respectively. The branch of the root is defined uniquely along a path, but
there is no analytic branch in a neighborhood of $0$ or $\infty$.
Moreover, this expression does not seem to be attracting, but it is used for
a preperiodic point here. So we only need it to be continuous in the sense that
$v_j\to0$ implies $\sqrt[{\scriptstyle d}\,]{c\cdot v_j}\to0$ or
$\sqrt[{\scriptstyle d}\,]{c/v_j}\to\infty$
for any branch. For $(u,\,v)\in\mathcal{N}_0$\,, $x_m'$ will be in a small
neighborhood of $0$ or $\infty$, and its preimages will be in small
neighborhoods of precritical points of $f$. The product of $\mathcal{N}_0$
with these neighborhoods defines the attracting neighborhood $\mathcal{N}$ for
$(u,\,v,\,w)$. Note that the coordinates $u$ and $v$ converge geometrically
as $\O(L^t)$, and $w$ with $\O(L^{t/d})$, or $\O(L^{t/d^2})$ if
$f^k(\infty)=0$ or $f^k(0)=\infty$. 

In \textbf{case 3}, $f^k(\infty)=0$ for some $k\ge1$, and the postcritical
point $z_m=g^k(z_{\alpha'})$ is identified with $z_{\beta'}$\,. Consequently,
iterates of $z_m$ are identified with corresponding iterates of
$z_{\beta'}$\,. If preimages of $z_{\beta'}$ are identified with preimages of
$z_m$ as well, or if additional non-postcritical marked points are identified
with critical or precritical points, they are labeled $w$ and treated
separately as in case~2 --- these points will be ignored from now on. The
pullback relation is not reducible in case~3: we have
$x_m'\sim\sqrt[{\scriptstyle d}\,]{c\cdot v_j}$ again, and this coordinate
cannot be treated separately, since it is pulled back to the critical value
$x_\alpha$\,. This value appears in the parameters of $f_{uv}$ and influences
the pullback of every point. Postcritical variables $v$ may appear directly in
these parameters, if $k=1$ or $f(0)$ has preperiod $1$. Note that
$x_m'\sim\sqrt[{\scriptstyle d}\,]{c\cdot v_j}$ is the only component of
(\ref{eqpbrmv}) not analytic in a neighborhood of $x=x^\infty$ or
$(u,\,v)=(\tilde z,\,0)$.

Choose the coordinates $(u,\,v)$ such that preperiodic iterates of $z_{\beta'}$
are of type $u$ and preperiodic iterates of $z_m$ are of type $u+v$, including
$x_m=v_m$\,. So $v_m'\sim\sqrt[{\scriptstyle d}\,]{c\cdot v_j}$ and the partial
derivative $\partial v_m'/\partial v_j\to\infty$ as a branch of the root is
continued analytically along the path. In a way, the matrix $D$ in
(\ref{eqmatbl}) has a unique infinite entry, in block $R$ and above the
diagonal. One idea to deal with this is to use the orbifold metric of $f$ for
$x_i\approx x_i^\infty=\tilde z_{D(i)}$ instead of the usual metric on $\C$.
Alternatively, we may lift the path and the pullback relation to new
coordinates $(U,\,V)$ with $U_i=u_i$ for all $i$ and, e.g., $V_m=v_m$ but
$V_j^d=v_j$\,; this gives $V_m'\sim\sqrt[{\scriptstyle d}\,]{c}\cdot V_j$\,,
where the branch of $\sqrt[{\scriptstyle d}\,]{c}$ is determined by the
chosen lift of a concrete path. The $v$-coordinates of iterates must be lifted
as well, but this gives an analytic pullback relation only when $g$ has a
postcritical cycle of the same period $p$ as $f$.
If we are in the second scenario, simply add a $p$-cycle to $g$ within the
cycle of disks; this does not change the pullback of the other points in the
$x$-coordinates, but when the new points are used as $u$-coordinates, this
allows to estimate the $v$-coordinates. --- In the lifted coordinates, we have
attracting eigenvalues and an attracting neighborhood as in case~1.

In any case, the neighborhood $\mathcal{N}$ can be chosen such that its
preimages under different branches of the pullback relation are either
contained in $\mathcal{N}$ or disjoint from it. Since path segments are
appended continuously, the given segment stays in $\mathcal{N}$ forever and
is attracted. \mybox

Note that when $\Gamma$ would be replaced with a non-homotopic multicurve
$\hat\Gamma\in G\cdot\Gamma$, which is grouping marked points in the same way,
the new $\tilde g$ may be obstructed or equivalent to a different $f$. On the
other hand, we have not used information on the homotopy class of $\Gamma$ in
the proof of Proposition~\ref{Pselconv}. Here the key point is the assumption
on the path in item~2: for a different $\hat\Gamma$, the attracting
neighborhood $\mathcal{N}$ would be the same, but there may be no path segment
of length $1$ contained in it. The same remark applies to the usual Thurston
Algorithm without identifications in fact: the pullback relation will have
several attracting fixed points, and one of these is chosen depending on the
isotopy class of $g$ or on the initial path segment.

\textbf{Proof of Theorem~\ref{Tselconv}.2:}
The marked point $z_{\gamma'}$  normalized to $x_{\gamma'}=1$ is chosen such
that together with the critical points at $x_{\beta'}=0$ and
$x_{\alpha'}=\infty$, it singles out the component $C$. The
associated embedding $\pi_3:\M\to\hat\C^{|Z|}$ extends continuously to
$\pi_3:\hat\M\to\hat\C^{|Z|}$ according to Proposition~\ref{Paug}.3; on
$\S_{G\cdot\Gamma}$ it is described as follows: Each $m\in\S_{G\cdot\Gamma}$
defines a noded Riemann surface; the piece corresponding to $C$ is isomorphic
to $\hat\C$. The isomorphism is unique by sending specific marked points and
nodes to $0,\,1,\,\infty$. Now marked points in other pieces are sent to
the same points as the corresponding nodes. So the fixed point $\hat\tau^\sC$
of $\sigma_\sC=\sigma_f$ has the following property: all $\tau\in\S_\Gamma$
with component $\tau^\sC=\hat\tau^\sC$ have $\pi_3(\pi(\tau))=x^\infty$.

Now suppose that the choice of $x_{\gamma'}=1$ was made according to the
restrictions from Proposition~\ref{Pselconv}, and obtain an attracting
neighborhood $\mathcal{N}$. Combine spherical metrics to define a metric $d$
on $\hat\C^{|Z|-3}\subset\hat\C^{|Z|}$ and choose $\delta>0$ such that the 
open ball of radius $2\delta$ around $x^\infty$ is contained in $\mathcal{N}$.
Proposition~\ref{Pselconv} applies with the same notation of $g$,
$\Gamma$, $C$, and $\hat\tau^\sC$. We shall start by constructing a path in
$\T\cup\S_\Gamma$; for suitable $0<i<I<\infty$ it goes from
$\tau_{n_{k(I)}+i}$ to $\sigma_g^i(\hat\tau_I)$, to
$\sigma_g^{i+1}(\hat\tau_I)$, and to $\tau_{n_{k(I)}+i+1}$\,.

There is a T-ball around $\hat\tau^\sC$ in $\T^\sC$, such that all
$\tau\in\S_\Gamma$ with $\tau^\sC$ in this ball satisfy
$d(\pi_3(\pi(\tau)),\,x^\infty)<\delta$. Choose $i$ according to
Proposition~\ref{Pselprfpc}.3 such that $\sigma_\sC^i(\hat\tau_I^\sC)$ and
$\sigma_\sC^{i+1}(\hat\tau_I^\sC)$ belong to this ball for $I>i$. Since
$\hat\M$ is compact, $\pi_3$ is uniformly continuous. Choose $I>i$ such that
$\eps(I)$ is sufficiently small, so $d_\sWP(m',\,m)<\eps(I)$ implies
$d(\pi_3(m'),\,\pi_3(m))<\delta$. Now we have 
$d_\sWP(\tau_{n_{k(I)}+i}\,,\,\sigma_g^i(\hat\tau_I))<\eps(I)$ and
$d_\sWP(\sigma_g^{i+1}(\hat\tau_I),\,\tau_{n_{k(I)}+i+1})<\eps(I)$ according
to Proposition~\ref{Pselprfpc}.2. The first and third segments of our
preliminary path shall be the corresponding WP-geodesics. The middle segment
from $\sigma_g^i(\hat\tau_I)$ to $\sigma_g^{i+1}(\hat\tau_I)$ shall be the
product of T-geodesics in the components of $\S_\Gamma$.

So with $n_*=n_{k(I)}+i$ we have constructed a preliminary path from
$\tau_{n_*}$ to $\tau_{n_*+1}$ in $\T\cup\S_\Gamma$, such that
$d(\pi_3(\pi(\tau)),\,x^\infty)<2\delta$ on this path. Since the ball is open
and the path is compact, we may choose a nearby path from $\tau_{n_*}$
to $\tau_{n_*+1}$ in $\T\cap(\pi_3\circ\pi)^{-1}(\mathcal{N})$.
The pullback of this path interpolates $(\tau_n)_{n\ge n_*}$ and projects to a
path in $\pi_3(\M)$, which stays in $\mathcal{N}$ and converges to $x^\infty$
according to Proposition~\ref{Pselconv}.

In the course of these proofs, several paths were constructed and discarded to
obtain convergence of the sequence $\pi_3(\pi(\tau_n))$. Now suppose
that a path $\tau_t$ is given from the start. Then for $\eps>0$ we want to
find $T\ge0$ with $d(\pi_3(\pi(\tau_t)),\,x^\infty)<\eps$ for $t>T$. This is
done by applying the result for sequences to the pullback of finitely many
intermediate points on the initial segment, which are chosen depending on
$\eps$, such that each smaller segment gives a change $<\eps/2$ in
$\pi_3(\M)$. Note again that the T-distance and thus the WP-distance stays
bounded under the pullback and that $\pi_3$ is uniformly continuous on
$\hat\M$.

Finally, consider the case where the marked point normalized to $1$ does not
satisfy the assumption of Proposition~\ref{Pselconv}. Then we have convergence
in a different normalization, and the two normalizations are related by an
affine rescaling with a convergent factor. \mybox

\begin{rmk}[Rate of convergence]\label{Rselconv}
1. The attracting eigenvalues at $x^\infty$ in configuration space were
related to multipliers of $f$ and to eigenvalues of $D\sigma_f$ in
Proposition~\ref{Pselconv}, where $\sigma_f$ includes both postcritical points
and additional marked points. Note that similar estimates apply to
collisions with critical points, and that additional marked points of $g$
without collisions converge to marked points of $f$ with a rate determined by
multipliers of $f$ as well. 
The orbifold metric \cite{bookm, book1mcm} provides uniform expansion
and uniform estimates for multipliers of $f$, especially $|\rho|\ge k^p$ for
the multiplier of a non-postcritical $p$-cycle. So if $k_f$ is a bound for
the eigenvalues of $D\sigma_f$ without additional marked points, then
$d(\pi_3(\pi(\tau_n)),\,x^\infty)$ asymptotically shrinks exponentially by
$\max(k^{-1},\,k_f)<1$ independently of the number of additional marked points
with or without collisions.

2. This bound on eigenvalues does not directly imply uniform convergence, e.g.,
in the case of a formal mating $g$ with fixed $\psi_0$ but an arbitrary number
of marked points: using a standard distance on $\hat\C^{|Z|-3}$, the norm of
the derivative may be arbitrarily large when eigenvectors are nearly parallel.
Moreover, the number of initial steps to get into a neighborhood of $x^\infty$
may grow with $|Z|$. See \cite{emate} for results on uniform convergence.

3. Under the assumptions of essential equivalence according to
Definition~\ref{Dselconv} and Theorem~\ref{Tselconv}, the leading eigenvalue is
always $\lambda_\Gamma=1$. In a different situation with $\lambda_\Gamma>1$,
collisions shall happen faster than exponentially.
\end{rmk}

\section{Construction and convergence of mating} \label{4}
We shall employ five different notions of mating: the formal mating is
constructed explicitly, modified to an essential mating, and it is
combinatorially equivalent to a rational map, the combinatorial mating. This is
a geometric mating at the same time, since it is conjugate to the topological
mating, which is defined as a quotient of the formal mating or of the
polynomials in turn. While the notion of the geometric mating may be most
natural, the construction best understood starts with the formal and
combinatorial matings in the postcritically finite case.
--- Convergence properties of the formal mating are discussed in
Section~\ref{4c} in a direct application of Theorem~\ref{Tselconv}.

\subsection{Dynamics and combinatorics of quadratic polynomials} \label{41}
The dynamics of a quadratic polynomial $f_c(z)=z^2+c$ is understood as follows:
all $z$ with large modulus are escaping to $\infty$ under the iteration; the
non-escaping points form the filled Julia set $\K_c$\,. By definition, the
parameter $c$ belongs to the Mandelbrot set $\M$, if $\K_c$ is connected, or
equivalently, if the critical orbit does not escape. Then the Boettcher map
$\Phi_c:\hat\C\setminus\K_c\to\hat\C\setminus\overline\disk$ maps dynamic rays
$\r_c(\theta)$ to straight rays with angle $\theta$ \cite{bookm, ser, mer}.

When $\theta$ is periodic or preperiodic under doubling, the landing point
$z=\gamma_c(\theta)\in\partial\K_c$ is periodic or preperiodic under $f_c$ as
well. In the parameter plane, parameters $c=\gamma_\sM(\theta)\in\partial\M$
are defined as landing points of parameter rays with rational angles. If
$\theta$ is periodic, $c$ is the root of a unique hyperbolic component with a
unique center; for that parameter, the critical orbit is periodic. Preperiodic
angles give Misiurewicz parameters, for which the critical value is
preperiodic. Dynamic rays landing together are important for ray connections.
For parameters $c$ in a limb of $\M$, the fixed point $\alpha_c$ has unique
angles and a unique rotation number, while the other fixed point always
satisfies $\beta_c=\gamma_c(0)$.

\subsection{Formal mating and combinatorial mating} \label{4m}
The formal mating $g=P\fmate Q$ of $P(z)=z^2+p$ and $Q(z)=z^2+q$ is a Thurston
map, which is conjugate to $P$ and $Q$ on the lower and upper half-spheres,
respectively. E.g., consider an odd map $\phi_0:\C\to\disk$ with
$\phi_0(r\cdot\e{\i\theta})\to\e{\i\theta}$ as $r\to\infty$ and set
$\phi_\infty(z)=1/\phi_0(z)$; then define $g=\phi_0\circ\ P\circ\phi_0^{-1}\cup
 \phi_\infty\circ Q\circ\phi_\infty^{-1}$\,. A simple explicit choice is given
by $\phi_0(z)=z/\sqrt{|z|^2+1}$, then $g$ will be smooth but not
quasi-regular. 
\textbf{External rays} of $g$ are unions of $\phi_0(\r_p(\theta))$ and
$\phi_\infty(\r_q(-\theta))$ together with a point on the equator;
\textbf{ray-equivalence classes} are maximal connected unions of rays and
landing points in $\phi_0(\partial\K_p)$ and $\phi_\infty(\partial\K_q)$.
Their geometry is described in \cite{rmate}.
According to \cite{rst, rs}, in the postcritically finite case there are:
\begin{itemize}
\item Cyclic ray connections corresponding to non-removable L\'evy cycles,
when the parameters $p$ and $q$ are in conjugate limbs of the Mandelbrot set.
\item Otherwise only trees giving identifications within and between
Julia sets, maybe in several steps. 
\item If postcritical or additional marked points are in the same
ray-equivalence class, these are surrounded by removable L\'evy cyles. Then an
\textbf{essential mating} $\tilde g$ is defined by modifying $g$:
these trees or disks are collapsed to points and the map is modified at
preimages as well, giving an unobstructed Thurston map with a smaller number
of marked points \cite{rst, rs}. 
\end{itemize}

The Thurston algorithm for $g$ gives a sequence of homeomorphisms $\psi_n$
and of rational maps $f_n$ with $\psi_n\circ g=f_n\circ\psi_{n+1}$\,. The
homeomorphisms $\psi_n$ converge up to isotopy, unless $g$ is obstructed
or of type \2. The following result is classical:

\begin{thm}[Combinatorial mating by Rees--Shishikura--Tan]\label{Tfecm}
Suppose the polynomials $P$ and $Q$ are postcritically finite and the
parameters are not in conjugate limbs of the Mandelbrot set. Then the
formal mating $g=P\fmate Q$ does not have a non-removable obstruction, and the
\textbf{combinatorial mating} $f$ is obtained as follows:

$a)$ If the formal mating $g$ does not have a removable obstruction, then $f$
is combinatorially equivalent to $g$, and the Thurston Algorithm for $g$
converges $f_n\to f$.

$b)$ If the formal mating $g$ has a removable obstruction, then $f$
is defined as the rational map equivalent to the essential mating
$\tilde g$. The Thurston Algorithm for $\tilde g$ converges
$f_n\to f$, unless $\tilde g$ has an orbifold of type \2.
\end{thm}

We may speak of ``the'' combinatorial mating, since M\"obius conjugacy
classes are avoided by assuming a normalization:
the critical point $0$ of $f$ corresponds to
$P$, the critical point $\infty$ to $Q$, and $1$ is the fixed point of
argument $0$. Different combinatorial matings might still be conjugate to each
other by marking a different fixed point, or by interchanging $P$ and $Q$.
E.g., the combinatorial matings of $z^2\pm\i$ with the Basilica $z^2-1$ are
distinct, but conjugate by a rotation of the fixed points. This ambiguity is
avoided with the alternative normalization $f(\infty)=1$.

\textbf{Idea of the proof:} In \cite{rst} it is shown that every obstruction
of a quadratic Thurston map contains a removable L\'evy cycle, or there is a
``good'' L\'evy cycle. The curves are homotopic to periodic ray-equivalence
classes. See Figure~\ref{Fkk}. Removable cycles correspond to loops around
simply connected ray-equivalence classes, while cyclic ray connections
indicate the presence of non-removable L\'evy cycles. Then there is a good
L\'evy cycle corresponding to closed ray connections between the two
$\alpha$-fixed points, which exist precisely when the parameters are in
conjugate limbs. Otherwise the essential mating $\tilde g$ can be defined
as a branched covering, which is unobstructed.

When $\tilde g$ has orbifold type not \2, the proof of existence and uniqueness
is completed with the Thurston Theorem~\ref{TT1},
but the case of type \2 is different: here the absence of
obstructions for $\tilde g$ is not sufficient to guarantee that there is an
equivalent rational map $f$. The proof can be given by applying the Shishikura
Algorithm to each essential mating in question, to determine the matrix of the
real-affine lift and to check that the eigenvalues are not real. The cases of
$1/4\fmate1/4$ and $1/6\fmate1/6$ are described in \cite{rst}, five more cases
are discussed in \cite{m14}, and the remaining two cases are settled
in \cite{pmate}. \mybox

For degree $d\ge3$ analogous definitions are used, but there is no
combinatorial characterization of cyclic ray connections in general.
Obstructions need not contain L\'evy cycles, and non-removable obstructions
may exist also when there are no cyclic ray connections \cite{st3}.
Moreover, the combinatorial mating will not be unique in the case of
flexible Latt\`es maps.

The \textbf{topological mating} $P\tmate Q$ is defined by collapsing all
rational and irrational ray-equivalence classes to points. Alternatively, take
the disjoint union of $\K_p$ and $\K_q$ and consider the equivalence relation
generated by $\gamma_p(\theta)\eqr\gamma_q(-\theta)$.
By the Moore Theorem \cite{maten}, we have a Hausdorff space homeomorphic to
the sphere, when the ray-equivalence relation of $P\fmate Q$ is closed and not
separating. 
Then the formal mating $g=P\fmate Q$ descends to a branched covering of the
quotient space, so the topological mating $P\tmate Q$ is a branched covering
of the sphere, which is defined up to conjugation. 
When mating polynomials from conjugate limbs of $\M$, the topological mating
does not exist because the sphere would be pinched. It may happen that there
is not even a Hausdorff space; examples of nested closed ray-equivalence
classes are given in \cite{rmate}.

Now the \textbf{geometric mating} is a rational map $f$  topologically
conjugate to the topological mating, $f\eqg P\tmate Q$.
In the postcritically finite case of mating, the following result from
\cite{rs, cpt} shows that every combinatorial mating from non-conjugate limbs
is a geometric mating in fact. Moreover, the topological mating exists and
there are no cyclic irrational ray connections either. Note that in contrast to
Theorems~\ref{Tfecm} and~\ref{Tsmrcc}, an orbifold of type \2 does not require
special considerations.
%

\begin{thm}[Rees--Shishikura]\label{Trs}
For postcritically finite quadratic polynomials $P$ and $Q$, with $p$ and $q$
not in conjugate limbs of $\M$, consider the formal mating $g=P\fmate Q$
and the essential mating $\tilde g$. According to Theorem~$\ref{Tfecm}$, the
combinatorial mating $f$ is combinatorially equivalent to $\tilde g$. Moreover:

$1$. There is a semi-conjugation $\Psi_\infty$ from $g$ to $f$.

$2$. $\Psi_\infty$ maps ray-equivalence classes to points, and it is
injective otherwise. So the topological mating is defined on a Hausdorff
space homeomorphic to the sphere, and $\Psi_\infty$ descends to a topological
conjugation from the topological mating $P\tmate Q$ to $f$. Now the geometric
mating of $P$ and $Q$ exists and it is given by $f$.

$3$. $\Psi_\infty$ is a uniform limit of homeomorphisms.
\end{thm}

Let us look at details from the proof for later reference in Section~\ref{4c}.
Note that the Thurston Theorem~\ref{TT1} and its application in
Theorem~\ref{Tfecm} used convergence of isotopy classes $[\psi_n]$ for an
arbitrary $\psi_0$\,, but now we need convergence of maps $\Psi_n$ for
a special choice of $\Psi_0$\,.

Idea of the \textbf{proof:} 1. The simplest case concerns preperiodic $P$ and
$Q$ without postcritical identifications. There are isotopic
$\Psi_0\,,\,\Psi_1$ with $\Psi_0\circ g=f\circ\Psi_1$ and a path
$\Psi_t\in[\Psi_0]$ with $\Psi_t\circ g=f\circ\Psi_{t+1}$ for $0\le t<\infty$.
Since $f$ is uniformly expanding with respect to the orbifold metric
\cite{bookm, book1mcm}, the homotopic length of a segment
$\{\Psi_t(z)\,|\,n\le t\le n+1\}$ shrinks exponentially in $n$, uniformly
in $z$. So $\Psi_\infty=\lim\Psi_n$ is continuous, surjective,
and a semi-conjugation.

The second case includes hyperbolic $P$ or $Q$. Then the orbifold metric is
more singular at periodic critical orbits, and exponential shrinking is
uniform away from these only. Although $\Psi_0$ can be chosen as a local
conjugation at these cycles by employing the B\"ottcher conjugation, this does
not guarantee $\Psi_1=\Psi_0$ there. In \cite{rs} the latter property is
obtained by modifying $\Psi_0$ with suitable Dehn twists. In \cite{cpt},
$\Psi_0$ is left unchanged but its pullback is described locally in terms of
Dehn twists.

The third case requires the construction of an essential mating $\tilde g$ by
collapsing a family $Y$ of critical and postcritical ray-equivalence classes.
Then we have $\Psi_t=\tilde\Psi_t\circ\pi_t$\,, where $\tilde\Psi_t$ is a
pseudo-isotopy for the essential mating. For $n\le t<n+1$, $\pi_t$
collapses the ray-equivalence classes in $g^{-n}(Y)$ to points independent of
$t$. So restricted to $n\le t\le n+1$, $\Psi_t$ is a pseudo-isotopy outside of
finitely many ray-trees.

2. The homotopic length of subarcs of rays shrinks exponentially, at least
away from precritical and postcritical classes. 
Moreover, any ray-equivalence class has a finite number of rays, so its
image under $\Psi_t$ shrinks to a point. It is quite involved to show
that distinct classes are mapped to distinct points by $\Psi_\infty$
\cite{rs}.

3. When $\tilde g=g$, $\Psi_\infty$ is the end of the pseudo-isotopy
$\Psi_t$\,. Otherwise both $\tilde\Psi_n$ converges to 
a continuous map and $\Psi_n$ converges to $\Psi_\infty$\,, but
$\Psi_n=\tilde\Psi_n\circ\pi_n$ is not a homeomorphism. Now the projections
$\pi_n$ are approximated by homeomorphisms, observing that ray-equivalence
classes are collapsed successively. \mybox

Now $\gamma_f(\theta)=\Psi\circ\phi_0\circ\gamma_p(\theta)
 =\Psi\circ\phi_\infty\circ\gamma_q(-\theta)$ is a semi-conjugation from the
angle doubling map on $\R/\Z$ to $f$ on its Julia set, which is not injective
but can be approximated by embeddings; it gives a Peano curve when both
$P$ and $Q$ are preperiodic. 
It maps the Brolin measure on $\R/\Z$ to the Lyubich measure on $\hat\C$.
A tiling is obtained from $T=\gamma_f([0,\,1/2])$ as well: then the Julia set
is $T\cup(-T)$ and $T\cap(-T)$ is the image of the spines intersected with the
polynomial Julia sets \cite{m14}. 
In general this gives no finite subdivision rule. Alternative constructions
with a pseudo-equator \cite{meyer} or Hubbard trees \cite{mw1} are possible
in certain cases.

\subsection{Convergence properties of the formal mating} \label{4c}
According to the discussion of Theorem~\ref{Tselconv}, there is no need to
correct a removable obstruction by identifying marked points manually: it will
be removed automatically during the iteration of the unmodified Thurston
Algorithm, in the sense that several marked points have the same limit, at
least in the non-\2 orbifold  case. Then $[\psi_n]$ diverges in
Teichm\"uller space, but the images of marked points and the rational maps
$f_n$ converge. Now the Thurston Algorithm can be implemented for the formal
mating without dealing with the combinatorics and topology of postcritical
ray-equivalence classes; the essential mating is used only as a step in the
proof, but not in the actual pullback. See Section~\ref{5} for a discussion of
slow mating. Actually, the same technique gives identifications for all
repelling periodic and preperiodic points by marking them in addition. As
conjectured in \cite{medusa, cts3}, e.g., the proof is based on the
Selinger extension to augmented Teichm\"uller space in Section~\ref{3}.

\begin{thm}[Convergence of maps \& rational ray-equivalence classes]%
\label{Tsmrcc}
Consider the Thurston
Algorithm $[\psi_n]$ with any initial $\psi_0$ for the formal mating
$g=P\fmate Q$ of postcritically finite quadratic polynomials $P$ and $Q$, with
$p$ and $q$ not in conjugate limbs of the Mandelbrot set. Moreover, assume that
the combinatorial mating $f$ has an orbifold not of type \2.

$1$. If the formal mating $g$ has removable obstructions, it is
essentially equivalent to the combinatorial mating $f$. 
The rational maps $f_n$ from the unmodified Thurston Algorithm converge
to $f$. The images of marked points of $g$ collide according to their
ray-equivalence classes under the iteration, and converge to marked
points of $f$.

$2$. In both cases, when $g$ is combinatorially equivalent or essentially
equivalent to $f$, consider the evolution of any periodic or preperiodic point
$z$, which corresponds to a point in $\partial\K_p$ or $\partial\K_q$\,: then
$x_n=\psi_n(z)$ converges to a periodic or preperiodic point of $f$.
Different points are identified in the limit, if and only if they belong to
the same ray-equivalence class.
\end{thm}

The second item is motivated by the videos of moving Julia sets, which are
computed from the slow mating algorithm and meant to represent equipotential
gluing \cite{emate}; it does not make sense when the formal mating is
considered only up to isotopy with respect to postcritical points. There are
two ways of understanding the statement in the context of the Thurston
Algorithm:
\begin{itemize}
\item The formal mating $g$ is defined such that it is topologically conjugate
to $P$ on the lower hemisphere and to $Q$ on the upper hemisphere. So there
are subsets of the sphere corresponding to the Julia sets $\K_p$ and $\K_q$
and points corresponding to periodic and preperiodic points of $P$ and $Q$.
Pick a homeomorphism $\psi_0$ and consider its lifts with
$f_n\circ\psi_{n+1}=\psi_n\circ g$ in a suitable normalization. When
$\tilde g=g$, the Thurston Theorem~\ref{TT1} shows that the homeomorphisms
$\psi_n$ converge up to homotopy with respect to the postcritical set of $g$,
and $\psi_n(z)$ converges when $z$ is a marked point. So here the latter
statement is extended to other points $z$ for the same sequence $\psi_n$\,, not
for any homotopic sequence.
\item A finite number of these periodic or preperiodic points of $g$ may be
marked in addition, giving a new pullback map on a higher-dimensional
Teichm\"uller space. Then $\psi_n$ may be considered up to homotopy with
respect to the marked set $Z$. If there are no collisions, the Thurston
Theorem~\ref{TT1} gives convergence immediately, but Theorem~\ref{Tselconv}
is needed in general.
\end{itemize}
\textbf{Proof:} Assume a normalization with critical points at $0$ and
$\infty$, and $1=f_n(\infty)$ or the fixed point on the $0$-ray is at $1$.
Actually, in the latter case we may mark $z=1$ and the two $\beta$-fixed points
in addition; so use the former normalization in the proof of item~1 and treat
the second normalization as a special case of item~2. 

1. According to \cite{rst}, all obstructions of $g$ contain removable L\'evy
cycles, which consist of loops around periodic ray-equivalence classes with at
least two marked points. A simple obstruction $\Gamma$ is obtained by adding
all essential preimages, which are loops around preperiodic ray-equivalence
classes containing at least two critical or postcritical points. Define
the essential mating $\tilde g$ by identifying all of these ray-equivalence
classes, or alternatively, all disks bounded by $\gamma\in\Gamma$, to points.
Then modify the map in neighborhoods of preimages containing at most one marked
point as well. This is done without destroying the orbit of a single
marked point within a disk; see Sections~\ref{32} and~\ref{3s}. Note that the
original definition \cite{rst, rs} may involve collapsing a larger number
of ray-equivalence classes with a single critical or postcritical point,
but all possible choices of $\tilde g$ are combinatorially equivalent. Now $g$,
$\Gamma$, and $\tilde g$ satisfy the assumptions of Definition~\ref{Dselconv}:
\begin{itemize}
\item Again by \cite{rst}, $\tilde g$ is unobstructed. Since its orbifold is
not of type \2, there is an equivalent rational map $f$, which defines the
combinatorial mating and the geometric mating in fact.
\item The critical points of $g$ are not identified in $\tilde g$, because then
$\tilde g$ would not be defined properly as a branched covering of degree $2$;
this happens only when there are non-removable obstructions and the parameters
are in conjugate limbs. No critical point is identified with $1$ either, using
a normalization different from $f(\infty)=1$ when $f(\infty)=0$. 
\item Loops bounding small tubular neighborhoods of disjoint simply connected
ray-equivalence classes define disjoint disks, so $\Gamma$ is not nested.
\item When a ray-equivalence class is not mapped homeomorphically, 
it contains a critical point of $P$ or $Q$, so it is preperiodic: periodic
critical points are superattracting and not accessible by external rays.
--- Note that the first-return map of the disk around a periodic
ray-equivalence class gives a homeomorphism of the corresponding piece, which
is always finite-order and not pseudo-Anosov \cite{book2h}, since the
postcritical points are connected by a tree mapped to itself.
\end{itemize}
So the formal mating $g$ is essentially equivalent to $f$, $\Gamma$ is the
canonical obstruction, and Theorem~\ref{Tselconv} gives convergence of
$f_n\to f$, and of colliding postcritical points as well. When $f$ is of type
\2, this statement is wrong in general 
\cite{pmate}.

2. Assume again that a postcritical point is normalized to $1$, which is not in
the same ray equivalence class as $0$ or $\infty$, and the fixed point on the
equator is marked in addition. Its convergence is obtained together with all
marked points, and the normalization can be changed afterward to $1$  by an
affine rescaling with a convergent factor. --- Given a finite number of
periodic or preperiodic points of $P$ and $Q$, add all of their images and all
ray-equivalent points, and consider the corresponding points in $\phi_0(\K_p)$
and $\phi_\infty(\K_q)$ together with the corresponding points on the equator
of $g$. Denote the union of postcritical ray-equivalence classes by $X$ and the
additional classes by $Y$, set $X'=g^{-1}(X)\setminus X$ and
$Y'=g^{-1}(Y)\setminus Y$. Now we have $X\cap Y=\emptyset$ and
$X'\cap Y'=\emptyset$, but we may have $X'\cap Y\neq\emptyset$. 

So there are finitely many disjoint ray-equivalence classes to consider.
Each of these is a tree, since otherwise the topological and geometric
matings would not exist. The essential mating $\tilde g$ shall be defined by
collapsing $X$ and modifying the new map in a small neighborhood of $X'$. The
essential map $\hat g$ for the larger set of marked points is defined by
collapsing $X\cup Y$ and modification in a neighborhood of $X'\cup Y'$. Denote
by $\Gamma$ the union of loops around the ray-equivalence trees in $X\cup Y$;
it is a simple obstruction again.

When a homeomorphism $\psi_0$ is chosen and the Thurston pullback
$f_n\circ\psi_{n+1}=\psi_n\circ g$ is applied, this gives the same
homeomorphisms $\psi_n$ and rational maps $f_n$ as in item~1; these maps do
not depend on the additional marked points, since the three normalized points
are critical or postcritical. So the question is, do the homeomorphisms
converge on the additional marked points, which follows when they converge
in the larger Teichm\"uller space, i.e, up to homotopy with respect to the
larger marked set. To apply Theorem~\ref{Tselconv} we only need to show that
$\hat g$ is unobstructed and equivalent to $f$. Otherwise for the larger set of
marked points, the canonical obstruction of $g$ would contain a loop around
several disks of $\Gamma$ or marked points of $g$.
\begin{itemize}
\item In the case without postcritical identifications in the formal mating, so
$\tilde g=g$, consider $\Psi_n$ from the proof of the Rees--Shishikura
Theorem~\ref{Trs}, which is defined by pulling back a specific homeomorphism
$\Psi_0$\,. Then $\Psi_n\to\Psi_\infty$\,, which is a semi-conjugation mapping
different ray-equivalence classes to different points. So the convergence claim
is true for $\psi_0=\Psi_0$ and $\hat g$ is unobstructed.
\item When $\tilde g\neq g$, we have
$\Psi_n=\tilde\Psi_n\circ\pi_n\to\Psi_\infty$\,, but $\Psi_0$ is not a
homeomorphism. So consider $\tilde\Psi_n$ instead, which are defined by a
pullback with the essential mating $\tilde g$. Since the essential map $\hat g$
is defined by collapsing ray-equivalence classes of $g$ in $X\cup Y$ and
modification around $X'\cup Y'$, an equivalent map can be defined as a
component map from $\tilde g$ as well. The limit of $\tilde\Psi_n$ exists
according to \cite{cpt} and if $\hat g$ was obstructed, then $\Psi_\infty$
would map different ray-equivalence classes to the same point.
\end{itemize}
Now Theorem~\ref{Tselconv} applies and gives convergence for any initial
$\psi_0$\,, in particular for slow mating and for equipotential gluing.
See Remark~\ref{Rselconv} and \cite{emate} for questions of uniform convergence
with respect to an arbitrary number of additional marked points.
\mybox

There are a few related ways to describe, which periodic or preperiodic points
of $g$ converge to which point of $f$\,:
\begin{itemize}
\item According to the Rees--Shishikura Theorem~\ref{Trs}, there is a
semi-conjugation $\Psi=\Psi_\infty$ from $g$ to $f$, which maps each
ray-equivalence class to a unique point.
\item Then $\Psi\circ\phi_0$ and $\Psi\circ\phi_\infty$ are partial
semi-conjugations from $P$ on $\K_p$ or $Q$ on $\K_q$ to restrictions of $f$.
\item $\gamma_f(\theta)=\Psi\circ\phi_0\circ\gamma_p(\theta)
 =\Psi\circ\phi_\infty\circ\gamma_q(-\theta)$ is a semi-conjugation from the
angle doubling map on $\R/\Z$ to $f$ on its Julia set.
\end{itemize}

\section{Implementation of slow mating} \label{5}
Equipotential gluing was defined by Milnor \cite{mquad}. 
For any radius $1<R<\infty$ the equipotential lines of potential $\log R$,
$|\Phi_p(z)|=R$ and $|\Phi_q(z)|=R$, are glued to form a sphere $S_R$ with
conformal images of $\K_p$ and $\K_q$\,. There are associated quadratic
rational maps $f_R: S_{\sqrt R}\to S_R$\,, and it is expected that $f_R$
converges to the conformal mating $f\eqg P\tmate Q$ as $R\to1$.

Slow mating shall denote any implementation of the Thurston Algorithm for the
formal mating $g=P\fmate Q$ of postcritically finite polynomials, which is
based on a pullback of a path in moduli space according to Section~\ref{2p}. A
particular initialization is given below; it is assumed to
approximate equipotential gluing when the initial radius $R_1$ is large, and
the limit $t\to\infty$ corresponds to $R\to1$. Buff and Ch\'eritat have
introduced slow mating to make movies of equipotential gluing, where the images
of $\K_p$ and $\K_q$ are moving and the identification of boundaries is visible
as a process. See the examples on the web pages
\href{http://www.math.univ-toulouse.fr/~cheritat/MatMovies/}%
{www.math.univ-toulouse.fr/$\sim$cheritat/} and on
\href{http://www.mndynamics.com/index2.html}{www.mndynamics.com}~.

Recall that the formal mating $g=P\fmate Q$ is defined as
$g=\phi_0\circ P\circ\phi_0^{-1}$ on the lower half-sphere $|z|<1$, and by
$g=\phi_\infty\circ Q\circ\phi_\infty^{-1}$ on the upper half-sphere $|z|>1$.
Denote the critical orbits by $P:0\Rightarrow p=p_1\to p_2\to\dots$ and
$Q:0\Rightarrow q=q_1\to q_2\to\dots$; when these are finite, the Thurston
Algorithm of the formal mating shall be implemented with a path in moduli space
according to Section~\ref{2p}. The homeomorphisms $\psi_0$ and $\psi_1$ are
chosen such that $x_i(t)=\psi_t(\phi_0(p_i))$ satisfies
$x_i(0)\approx p_i/R_1^2$ and $x_i(1)\approx p_i/R_1$ for some initial radius
$R_1\gg2$; analogously we require that $y_i(t)=\psi_t(\phi_\infty(q_i))$
satisfies $y_i(0)\approx R_1^2/q_i$ and $y_i(1)\approx R_1/q_i$\,. There are
two reasons for this: first, it is considered as an approximation of
equipotential gluing. Second, we will obtain simple explicit formulas, which
allow to check that there is a corresponding path in Teichm\"uller space. Now
$\psi_1$ shall be the pullback of $\psi_0$\,, so
$f_0\circ(\psi_1\circ\phi_0)=(\psi_0\circ\phi_0)\circ P$ and
$f_0\circ(\psi_1\circ\phi_\infty)=(\psi_0\circ\phi_\infty)\circ Q$ for the same
quadratic rational map $f_0$ with $f_0(z)\approx z^2$.

Interpolate the radius as $\log(R_t)=2^{1-t}\log(R_1)$ for $0\le t<\infty$. For
the initial path segment $0\le t\le 1$ we cannot take $x_i(t)=p_i/R_t$ and
$y_i(t)=R_t/q_i$\,, since there would be no quadratic rational map $f_0$
sending $x_{i-1}(1)$ and $y_{i-i}(1)$ to $x_i(0)$ and $y_i(0)$, respectively.
We might focus on $t=0$, define $x_i(0)$ and $y_i(0)$ explicitly and determine
$f_0$\,, $x_i(1)$, and $y_i(1)$ accordingly, but to avoid computing and
estimating square-roots, we shall take the opposite direction: choose the
simplest formula at time $t=1$. Since $z=1$ is assumed to be fixed, consider
\be
 x_i(1)=\frac{p_i}{R_1} \quad,\quad
 y_i(1)=\frac{R_1}{q_i} \quad,\quad\mbox{and}\quad
 f_0(z)=\frac{1+q/R_1^2}{1+p/R_1^2}\cdot
 \frac{z^2+p/R_1^2}{1+q/R_1^2\,z^2} \quad .\ee
From $p_{i-1}^2=p_i-p$ and $q_{i-1}^2=q_i-q$, we obtain at time $t=0$
\ban
 x_i(0)\:=\:f_0(x_{i-1}(1)) &=&
 \frac{1+q/R_1^2}{1+p/R_1^2}\cdot\frac{p_i/R_1^2}{1+q/R_1^4\,(p_i-p)}
 \:\approx\:\frac{p_i}{R_1^2} \quad\quad\mbox{and}\\[1mm]
 y_i(0)\:=\:f_0(y_{i-1}(1)) &=&
 \frac{1+q/R_1^2}{1+p/R_1^2}\cdot
 \frac{R_1^2\Big(1+p/R_1^4\,(q_i-q)\Big)}{q_i}
 \:\approx\:\frac{R_1^2}{q_i}\ .\ean
Note that in the periodic case, this choice of $f_0$ is optimal, because the
periodic critical points stay at $0$ and $\infty$. In the preperiodic case,
there are opposite marked points for $t\ge1$, which are no longer opposite for
$0\le t<1$. This is violating our usual normalization, and in a sense it is
moving the unmarked critical points, but it does not really matter for $t\ge1$.
Now we shall use the following interpolation $x_i(t)$ and $y_i(t)$:

\begin{init}[Slow mating]\label{Ism}
For postcritically finite polynomials $P$ and $Q$ with critical orbits $(p_i)$
and $(q_i)$, the unmodified Thurston Algorithm for the formal mating
$g=P\fmate Q$ can be implemented with a path in moduli space as follows:
Fix $R_1\ge5$ and interpolate the radius as $\log(R_t)=2^{1-t}\log(R_1)$ for
$0\le t\le1$. Set
\ban
 x_i(t) &=&
 \frac{1+(1-t)q/R_1^2}{1+(1-t)p/R_1^2}\cdot
 \frac{p_i/R_t}{1+(1-t)q/R_1^4\,(p_i-p)}
 \:\approx\:\frac{p_i}{R_t} \quad\mbox{and}\\[1mm]
 y_i(t) &=&
 \frac{1+(1-t)q/R_1^2}{1+(1-t)p/R_1^2}\cdot
 \frac{R_t\Big(1+(1-t)p/R_1^4\,(q_i-q)\Big)}{q_i}
 \:\approx\:\frac{R_t}{q_i} \ .\ean
The initial path for $0\le t\le1$ can be pulled back uniquely for $1<t<\infty$,
choosing the sign in \emph{(\ref{eqpblf})} by continuity. For the
homeomorphism $\psi_0$ from the proof, the marked points
$\pi(\sigma_g^n([\psi_0])$ are given by $x_i(n)$ and $y_i(n)$ in the usual
normalization of $0,\,1,\,\infty$.
\end{init}

Actually, when $R_1\ge10^{10}$, this initialization is the same as
$x_i(t)=p_i/R_t$ and $y_i(t)=R_t/q_i$ in $8$-byte precision. ---
The rational maps $f_t$ with $f_t\circ\psi_{t+1}=\psi_t\circ g$ are obtained
from the critical values $x_1(t)$ and $y_1(t)$. Now the pullback is computed
as follows (read $x$ or $y$ for $z$):
\be\label{eqpblf}
 z_i(t+1)=\pm\sqrt{\frac{1-y_1(t)}{1-x_1(t)}\cdot
 \frac{z_{i+1}(t)-x_1(t)}{z_{i+1}(t)-y_1(t)}} \quad\mbox{for}\quad t\ge0\ .\ee
Note that the first factor is dropped here and above for self-matings with
$p=q$. The formulas of initialization are simplified in general, when the
alternative normalization $f_t(\infty)=1$ instead of $f_t(1)=1$ is used.
See the Examples~\ref{Xsmib} and~\ref{Xsmbb}.

\textbf{Proof:} For $0\le t\le1$ we need a path of homeomorphisms $\psi_t$\,.
In a neighborhood of the filled Julia sets $\phi_0(\K_p)$ and
$\phi_\infty(\K_q)$, it is defined in terms of the following M\"obius
transformations, which are chosen such that marked points are mapped according
to the given initialization. For suitable discs within $|z|\le4$ set
\ban
\psi_t(\phi_0(z)) &=&
 \frac{1+(1-t)q/R_1^2}{1+(1-t)p/R_1^2}\cdot\frac{z/R_t}{1+(1-t)q/R_1^4\,(z-p)}
 \:\approx\:\frac{z}{R_t} \quad\mbox{and}\\[1mm]
 \psi_t(\phi_\infty(z)) &=&
 \frac{1+(1-t)q/R_1^2}{1+(1-t)p/R_1^2}\cdot
 \frac{R_t\Big(1+(1-t)p/R_1^4\,(z-q)\Big)}{z}
 \:\approx\:\frac{R_t}{z} \ .\ean
Between these two discs, the maps are interpolated to obtain homeomorphisms of
spheres. We do not want to introduce a twist accidentally, so we shall assume
that $\psi_t$ is close to the identity at the equator, and that it maps
the positive real axis close to itself. All of this makes sense when $R_1\ge5$,
since $|\psi_t(\phi_0(z))|<1$ and $|\psi_t(\phi_\infty(z))|>1$ for
$|z|\le4$, and since the maps are close to $z/R_t$ and $R_t/z$ there.

According to Proposition~\ref{Ppath}, we need to check
$f_0\circ\psi_1=\psi_0\circ g$. Then our initialization is the projection of
a suitable initial path in Teichm\"uller space for $0\le t\le1$, and the
pullback in moduli space is the projection of the desired path in Teichm\"uller
space, which is interpolating $\sigma_g^n([\psi_0])$ for $n\in\N$.

Our choices are compatible such that
$f_0\circ(\psi_1\circ\phi_0)=(\psi_0\circ\phi_0)\circ P$ and
$f_0\circ(\psi_1\circ\phi_\infty)=(\psi_0\circ\phi_\infty)\circ Q$ in
neighborhoods of the filled Julia sets, so $f_0\circ\psi_1=\psi_0\circ g$
there. If $\phi_0$ was extended arbitrarily to the annulus, it might have an
additional twist compared to the correct pullback. By looking at preimages of
the  positive real axis again, we see that this is not the case, so we may
assume that $\psi_0$ and $\psi_t$ are defined such that
$f_0\circ\psi_1=\psi_0\circ g$ globally, and $\sigma_g([\psi_0])=[\psi_1]$.
Note again that $\psi_t$ is not odd for $0\le t<1$, but the pullback will be
odd for $1\le t<\infty$.
\mybox

So in any case of postcritically finite polynomials $P$ and $Q$, there is a
simple numerical method to compute the pullback of a suitable homeomorphism,
or more precisely, the projection $\pi(\sigma_g^n([\psi_0])$ to moduli space.
The discretization of a path segment shall be discussed in \cite{amate, smate}.
The relation of the slow mating algorithm to equipotential gluing and the
visual representation will be explored further in \cite{emate}. Note that
the algorithm itself is quite fast, if you only compute the maps and not the
Julia sets, and the movies slow it down to illustrate the process.

When the formal mating $g$ is unobstructed, then $f_t$ converges to
the geometric mating $f\eqg P\tmate Q$, since the slow mating algorithm is just
an explicit method to compute the pullback. When there are removable
obstructions, the maps converge as well according to Theorem~\ref{Tsmrcc},
at least if the orbifold of $\tilde g$ or $f$ is not of type \2. Marked points
of $g$ are identified automatically in the limit, and there is no need to
encode the topology of postcritical ray-equivalence classes to define the
pullback of the essential mating $\tilde g$ instead.

However, the pullback $\pi(\sigma_g^n([\psi_0]))$ will diverge in general
when $\tilde g$ has an orbifold of type \2.
See Section~5 in \cite{pmate}. Probably the path accumulates on a
subset of the center manifold, which is disjoint from the fixed point, but
there are other possibilities when the path intersects the stable manifold.
--- Even when the formal mating $g$ has non-removable obstructions, the
pullback can be computed from the slow mating algorithm, but it will diverge
in a different way; e.g., $f_t$ may converge to a constant map.

\section{Captures and matings} \label{6}
Captures and precaptures are ways to construct a Thurston map by shifting a
critical value to a preperiodic point; we shall see that precaptures are
related to matings with preperiodic polynomials in fact.

\subsection{Hyperbolic and repelling-preperiodic captures} \label{6c}
These constructions rely on the concept of shifting or pushing a point from
$a$ to $b$ along an arc $C$. This means that a homeomorphism $\phi$ is chosen,
which is the identity outside off a tubular neighborhood of $C$, and such
that $\phi(a)=b$. We may assume that an unspecified point close to $a$ is
mapped to $a$ and $b$ is mapped to an arbitrary point nearby.

\begin{prop}[and definition]\label{Phrc}
Suppose $P$ is a postcritically finite quadratic polynomial and
$z_1\in\K_p$ is preperiodic and not postcritical. Let the new
postcritical set be $P_g=P_P\cup\{P^n(z_1)\,|\,n\ge0\}$. Consider an arc $C$
from $\infty$ to $z_1$ not meeting another point in $P_g$ and choose a
homeomorphism $\phi$ shifting $\infty$ to $z_1$ along $C$, which is the
identity outside off a sufficiently small neighborhood of $C$. Then:\\[1mm]
$\bullet$ $g=\phi\circ P$ is well-defined as a quadratic Thurston map with
postcritical set $P_g$\,. It is a \textbf{capture} if $z_1$ is eventually
attracting and a \textbf{precapture} in the repelling case.\\[1mm]
$\bullet$ The combinatorial equivalence class of $g$ depends only on the
homotopy class of the arc $C$.
\end{prop}

\textbf{Proof:} By construction, $g$ is a postcritically finite branched cover,
when the neighborhood of $C$ does not include any postcritical point except
$z_1$\,. Note that the preimages of $z_1$ under $P$ are mapped to some
arbitrary point by $g$, so if $z_1$ was periodic or postcritical, $g$ would not
be well-defined. Finally, if we have two different homeomorphisms $\phi$ and
$\phi'$ along the same curve or along two homotopic curves, then
$g'=(\phi'\circ\phi^{-1})\circ g$ and the homeomorphism $\phi'\circ\phi^{-1}$
is isotopic to the identity, since the appended path $C'\cdot C^{-1}$ is
contractible relative to $P_g\setminus\{z_1\}$. \mybox

Consider the following applications and possible generalizations:
\begin{itemize}
\item If a capture $g=\phi\circ P$ is combinatorially equivalent to a rational
map $f$, this gives a hyperbolic map of capture type. Let us say that $f$ is a
\textbf{Wittner capture}, if the capture path $C$ is homotopic to a rational
external ray followed by an internal ray of $P$; this construction is due to
Ben Wittner \cite{matew} and Mary Rees \cite{rees1}. Note that Rees denotes
only Wittner captures as captures, while general captures are called maps
of type~III. Maps of this type are never matings, but they may have a
representation as an anti-mating \cite{amate}.
\item Precaptures along external rays are related to matings in the following
Section~\ref{6i}.
\item Precaptures apply not only to polynomials $P$, but to rational maps in
general as long as the other critical orbits are finite. This construction
provides a finite \textbf{regluing} followed by a possible combinatorial
equivalence. In a more general situation, a countable regluing is followed by a
semi-conjugation \cite{trg, mtrg}. 
\item \textbf{Recapture} means that the finite critical value $P(0)$ is
shifted to a preimage of $0$, resulting in a Thurston map equivalent to a
hyperbolic polynomial. Relations to internal addresses and to Dehn twisted maps
are discussed in \cite{smate}.
\end{itemize}

\begin{init}[Captures and precaptures]\label{Ihrc}
Consider a capture or precapture $g=\phi\circ P$ according to
Proposition~$\ref{Phrc}$. Then the Thurston Algorithm is implemented by pulling
back a path in moduli space, which is initialized as follows: normalize $P$
such that the critical points are $0,\,\infty$ and another point in
$P_g\setminus\{z_1\}$ is $1$. For $0\le t\le1$, $x_1(t)$ moves from $\infty$ to
$z_1$ along $C$, while all of the other marked points stay fixed.
\end{init}

Under a non-conjugate-limbs condition, Wittner captures are unobstructed
\cite{rees1} and precaptures along external rays have only obstructions
satisfying the assumptions of Theorem~\ref{Tselconv}; see below. So the
sequence of rational maps converges to a rational map $f$, unless the orbifold
of $f$ is of type \2: then the sequence does not converge in general, but it
might converge for a special choice of $C$.

\textbf{Proof:} Note that when the preperiod of $z_1$ is one, the corresponding
periodic point satisfies $\psi_t(-z_1)=-x_1(t)$ only for $t\ge1$. Since
$\phi^{-1}\circ g=P\circ\mathrm{Id}$ and $P$ is holomorphic, we have
$[\mathrm{Id}]=\sigma_g([\phi^{-1}])$ and we may initialize the Thurston
Algorithm with a path $\psi_t$ from $\psi_0=\phi^{-1}$ to $\psi_1=\mathrm{Id}$.
Now $\phi^{\pm1}$ is the identity outside off a small neighborhood of $C$, so
$\psi_t$ can be chosen such that it moves $x_1(t)=\psi_t(z_1)$ from
$\phi^{-1}(z_1)=\infty$ to $z_1$ along $C$, and leaves the other marked points
untouched. By Proposition~\ref{Ppath} the projection from $\T$ to $\M$ defines
a suitable initialization to compute the Thurston pullback $\pi(\sigma_g^n)$
from an explicit pullback in moduli space. \mybox

To illustrate the process of slow capture or precapture, we may also define a
sequence or path of images $\psi_t(\K_p)$ of the filled Julia set, which is
constant $\K_p$ for $0\le t\le1$. It will show more and more identifications
happening by a piecewise pseudo-isotopy. See also the videos on
\href{http://www.mndynamics.com/index2.html}{www.mndynamics.com}\ . A similar
initialization is used for Dehn twisted maps; see \cite{bn} and the
Examples~3.1 and~3.7 in \cite{pmate}.

\subsection{Precaptures and matings} \label{6i}
The representation of matings by precaptures along external rays is motivated
by remarks in \cite{mtrg, reesc}. In the former paper,
the boundary of a capture component in $V_n$ is described by matings, which are
related to the postcritically finite map of capture type by regluing. This
means that the critical value is shifted from $\infty$ along an external ray
followed by an internal ray, and then moved back along an internal ray. So can
the mating be constructed by shifting the critical value directly from $\infty$
to $z_1=\gamma_p(\theta)$ along the external ray $\r_p(\theta)$\,? This is true
in general when $z_1$ is preperiodic, not only when it is on the boundary of a
hyperbolic component, but we shall not discuss postcritically infinite maps
here. 

\begin{thm}[Matings as precaptures, following Rees]\label{TRprec}
Suppose $P$ is postcritically finite and $\theta$ is preperiodic, such that
$q=\gamma_\sM(-\theta)$ is not in the conjugate limb and
$z_1=\gamma_p(\theta)\in\partial\K_p$ is not postcritical. Then the precapture
$g_\theta=\phi_\theta\circ P$ along $\r_p(\theta)$ is combinatorially
equivalent or essentially equivalent to the geometric mating $f$ defined by
$P\tmate Q$.
\end{thm}

So if $P\tmate Q$ is not of type \2, any implementation of the Thurston
pullback for $g_\theta$ gives a converging sequence of rational maps; e.g.,
Initialization~\ref{Ihrc} applies. The normalization $\beta_p=1$ ensures
$f(1)=1$. Note that the precapture does not work if both $P$ and $Q$ are
hyperbolic; then there is an alternative construction with two paths
\cite{rees1}. 
When only one of the two polynomials is hyperbolic, then either $P\tmate Q$ or
$Q\tmate P$ is a precapture. And when both are critically preperiodic, then
both $P\tmate Q$ and $Q\tmate P$ are precaptures, unless a critical point is
iterated to the other critical point: then $\infty$ shall be iterated to $0$.
--- By choosing precaptures along homotopic external rays, examples of shared
matings are obtained in \cite{rmate}.

Recall the notation $g$ and $\tilde g$ for the formal mating and the essential
mating; we shall see below that there is an essential precapture
$\tilde g_\theta$ as well. Before showing $\tilde g_\theta\eqth\tilde g$ let us
consider a few examples, to see how identifications happen and why they may
happen in different ways for $g$ and $g_\theta$\,:
\begin{itemize}
\item When $g=9/56\fmate1/4$, so $\theta=3/4$, there are no postcritical
identifications: $\tilde g=g$ and $\tilde g_\theta=g_\theta$\,. The precapture
can be constructed from the formal mating by shifting all postcritical points
in $\phi_\infty(\K_q)$ to $\phi_0(\K_p)$ along external rays, so $g_\theta$ and
$g$ are combinatorially equivalent.
\item In reverse order we have $\tilde g=g=1/4\fmate9/56$ again, but
$\tilde g_\theta\neq g_\theta$ for $\theta=47/56$ and $p=\gamma_\sM(1/4)$. Now
$g_\theta(\infty)$ has preperiod and period three, but
$\tilde g_\theta(\infty)$ has period one. The shift $\phi_\theta$ creates a
subset of the lamination with angle $\theta$ in the exterior of $\K_p$\,, so
there is a triangle connecting $3/7$, $5/7$, $6/7$ with a homotopic preimage
under $g_\theta$\,; pinching the surrounding L\'evy-cycle gives
$\tilde g_\theta$\,.
\item The converse happens for $g=1/4\fmate3/14$, so $p=\gamma_\sM(1/4)$ and
$\theta=11/14$. Now both $q=\gamma_\sM(3/14)$ and $g(\infty)$ have preperiod
one and period three, while $\tilde g\neq g$ has period one. But this
identification is immediate in the precapture $g_\theta=\tilde g_\theta$\,,
since $z_1=-\alpha_p$\,.
\item Both phenomena happen at the same time for $g=3/14\fmate3/14$, so
$\theta=11/14$. In $g_\theta$ the $3$-cycle of $P$ is collapsed by a triangle
in the exterior, while the $3$-cycle of $Q$ is identified with $\alpha_p$
immediately. We have $\tilde g_\theta\neq g_\theta\not\eqth g\neq\tilde g$.
\end{itemize}
For longer ray connections, there may be a similar splitting of branch points
and similar immediate identifications, but otherwise the precapture can be
understood in terms of the same ray-equivalence classes as the formal mating:

\textbf{Proof of Theorem~\ref{TRprec}:} Denote by $X$ the union of all
postcritical ray-equivalence classes of the formal mating $g=P\fmate Q$. Define
another Thurston map $g^\theta$ by shifting the critical value $\phi_\infty(q)$
to $\phi_0(z_1)$ along $\r_\theta$\,, without modifying $g$ on $X$. Consider
the extended Hubbard tree $T_p\subset\K_p$\,, which consists of regular arcs
connecting the postcritical points of $g_\theta$\,. Then
$g_\theta:T_p'\to T_p$\,, where $T_p'=T_p$ except for a slight detour at
$P^{-1}(0)$. We may assume that $g^\theta\circ\phi_0=\phi_0\circ g_\theta$ in a
neighborhood of $T_p$\,. So the two maps are combinatorially equivalent, even
if we mark the critical point $\infty$ in addition, since all other marked
points are contained in $T_p$ and $T_p$ is connected.

Now consider a path of Thurston maps $g_t$\,, such that postcritical points of
$P$ stay fixed in $\phi_0(\partial\K_p)$ and all postcritical points of $Q$
move from $\phi_\infty(\partial\K_q)$ to $\phi_0(\partial\K_p)$ along external
rays of $g$. This deformation is a kind of two-sided pseudo-isotopy from $g$
to $g^\theta$, since marked points may collapse in different ways on both ends,
while each component of $X$ is invariant under each $g_t$\,.
By collapsing all components of $X$ to points and modification at preimages,
equivalent quotient maps are obtained for all $g_t$\,, in particular for
$g$ and $g^\theta$, where postcritical points have been identified already
in different ways. So we know that $\tilde g^\theta=\tilde g\eqth f$ and we
may consider $\tilde g^\theta$ as an essential map in the sense of
Definition~\ref{Dselconv}, with $\Gamma$ consisting of loops around those
trees in $X$, which contain at least two postcritical points of $g^\theta$.
So $g^\theta$ is essentially equivalent to $f$, combinatorially equivalent if
$\Gamma=\emptyset$, and the same applies to the original precapture
$g_\theta$\,. \mybox

\section{Conclusion} \label{7}
Suppose $P$ and $Q$ are postcritically finite polynomials, not in conjugate
limbs of the Mandelbrot set. Then the formal mating $g=P\fmate Q$ is
combinatorially equivalent or essentially equivalent to a rational map
$f\eqg P\tmate Q$. Consider the following implementations of the Thurston
Algorithm for the formal mating, which should converge except for numerical
cancellations, unless $f$ is of type \2:
\begin{description}
\item[The medusa algorithm] was developed by Christian Henriksen and others
under the guidance of John Hamal Hubbard \cite{medusa}. Start with a Thurston
map having marked points on two circles at specific angles; it is equivalent to
the formal mating unless there are Misiurewicz parameters of satellite type
--- then the arguments from Section~\ref{6i} apply. A medusa is a graph
connecting the images of marked points, which corresponds to the equator
united with external rays from the equator to these points. Its pullback up to
homotopy with rational maps provides a unique choice of preimages. Since this
is an implementation of the Thurston Algorithm, the marked points and maps
should converge according to Theorem~\ref{Tsmrcc} or to Section~\ref{6i},
unless the orbifold of $f$ has type \2.

However, medusa often seems to be numerically unstable even for simple
examples: it begins to converge but after 50--100 iterations it oscillates
wildly. It is not known whether this is a bug in the implementation, an unlucky
choice of numerical parameters, or an unavoidable feature of this algorithm. I
had expected that instability would be related to long ray connections
converging to periodic points with a multiplier causing spiraling: then the
equator would have to spiral as well and cannot be pruned to a homotopic curve
with few long segments. But this idea was not confirmed by experiments; e.g.,
medusa did converge for $3/7\fmate3/14$ nd $12/31\fmate19/62$, which have
postcritical ray-equivalence classes of length four, and for $31/96\fmate1/3$
and $511/1536\fmate1/3$, which show significant spiraling. On the other hand,
it diverged even in cases without postcritical identifications, e.g.,
$1/14\fmate1/4$ and $19/60\fmate1/3$. Note also that for the matings
$5/28\fmate13/28$ and $7/60\fmate29/60$ of type \2, the Thurston
pullback accumulates on a four-cycle in configuration space according to
\cite{pmate}; medusa shows this behavior initially but oscillates after a few
more iterations.
%
\item[Triangulations of the sphere] are used by Laurent Bartholdi in the
GAP-package IMG \cite{img, gap}.
A Thurston map is represented algebraically as a biset \cite{bn},
and it is easy to combine maps, as in a formal mating, or to apply a Dehn
twist. Then a triangulation is constructed from the biset, which represents
an isotopy class of homeomorphisms. It is pulled back to implement the Thurston
Algorithm, with appropriate refinement and pruning.

When marked points get close to each other, the pullback is interrupted in the
current version and an obstruction is searched instead, based on the assumption
that points will be grouped in the observed way. In the case of formal matings
with removable obstructions, this approach might be modified such that either
the Thurston pullback is restarted with a component map, or such that the
iteration is continued to allow a collapse of marked points according to
Theorem~\ref{Tsmrcc}. In the latter case, the pullback might become unstable,
when a spiraling of marked points requires an excessive refinement of the
triangulation.
\item[Slow mating] according to Section~\ref{5} is much simpler to implement
\cite{mandel}. In the case of postcritical identifications, the path is not
required to follow a spiraling equator. So there is a good chance to converge
with just a small number of segments per marked point, and the algorithm will
still be fast with a large number of segments. However, any discretization of
a continuous path as a polygonal path should check, whether the exact pullback
to piecewise arcs can be replaced homotopically with piecewise line segments
again. This is easy in the case of quadratic polynomials \cite{smate},
but more involved for quadratic rational maps \cite{amate}. Note that there
may be a trade-off as well when using many small segments: homotopy
violations will happen less often, but detecting them will be numerically
less stable.
\item[An initialization by angles] will be more convenient, but slow mating
assumes that the parameters $p$ and $q$ are given as floating-point
approximations. When angles are given instead, either run the spider algorithm
first to determine these parameters, or draw the parameter rays and improve the
endpoints with Newton method. Alternatively, the slow mating algorithm can be
modified such that the marked points are on two circles initially; the pullback
would give the same marked points as medusa, but be more stable.
%
\item[A precapture] can be implemented according to Initialization~\ref{Ihrc}
as well, using an approximation to dynamic rays in this case. 
\end{description}

So this paper suggests to treat removable obstructions by ignoring them, which
is simple and fast: trust the slow mating algorithm to converge nevertheless.
This route is taken naturally by equipotential gluing as well \cite{emate}.
If you want to collapse ray-equivalence classes manually, you can determine
the relevant angles recursively from the conjugate angle algorithm
\cite{rmate}, but the topological part may be harder. When there are only
direct connections between postcritical points of $P$ and $Q$, so a
pseudo-equator exists \cite{meyer}, the modification can be done by taking a
medusa with all points on the equator. Mary Wilkerson \cite{mw1, mw2} has an
alternative implementation in this case: the pullback is controlled by a finite
subdivision rule, which is constructed from Hubbard trees.

\normalsize\newpage\appendix
\section{The spider algorithm} \label{As}
\textsl{In \cite{smate}, the Thurston Algorithm with a path in moduli space is
implemented for quadratic polynomials, including the spider algorithm, twisted
polynomials, precapture and recapture, and slow tuning. This section sketches
the discussion of the spider algorithm, because it is another application of
the convergence Theorem~\ref{Tselconv}; in fact it was the original
motivation for this research.}

For an angle $\theta\in\Q\setminus\Z$, we want to determine the associated
postcritically finite parameter $c$ of a quadratic polynomial $f_c(z)=z^2+c$.
From $\theta$ a Thurston map $g_\theta$ is constructed, and the Thurston
Algorithm shall give $f_c$\,. Denote the iterates of $\theta=\theta_1$ by
$\theta_i=2^{i-1}\theta$, $i\ge1$, and the preperiod and period of $\theta$
is $k$ and $p$. Consider the map $g_\theta=\phi_\theta\circ F$ with $F(z)=z^2$.
Here the homeomorphism $\phi_\theta$ is the identity in most of $\hat\C$; it
shifts the straight ray with angle $\theta_1$ radially out by $1$, and if
$k=0$, it shifts the ray with angle $\theta_p$ in by $1$. So
$\phi_\theta(0)=\e{\i2\pi\theta_1}$, and in the periodic case
$\phi_\theta(e^{\i2\pi\theta_p})=0$. The straight spider is invariant under
$g_\theta$\,.

To apply the Thurston Algorithm, we need to pull back marked points
$x_i(n)$ with quadratic polynomials. The choice of branch for the square roots
is determined by the pullback of an isotopy class of homeomorphisms. The
basic idea of the \textbf{spider algorithm} is: Teichm\"uller space is
represented by spiders, homotopy classes of graphs with legs from $\infty$ to
the marked points, which are pulled back with the polynomials. According to
\cite{hss, bfh, book2h} we may consider these cases:
\begin{description}
\item[Case 1:] The angle $\theta$ is periodic and $c$ is the center
associated to the root $\gamma_\sM(\theta)$. Then $g_\theta$ is
combinatorially equivalent to $f_c$ and unobstructed. Under the equivalence,
the spider legs are homotopic to external rays extended by internal rays, which
will have common points in the satellite case.
\item[Case 2:] The angle $\theta$ is preperiodic and the Misiurewicz point
$c=\gamma_\sM(\theta)$ is an endpoint or of primitive type. Again, $g_\theta$
is unobstructed and equivalent to $f_c$\,. The spider legs correspond to
external rays at the postcritical orbit.
\item[Case 3:] The Misiurewicz point $c=\gamma_\sM(\theta)$ is of satellite
type; the angle $\theta_{k+1}$ has period $p=rq$ and the landing point has
period $q<p$. Now $g_\theta$ has a L\'evy cycle with $q$ curves, each
containing $r$ marked points. By identifying these points manually, or by
extending the spider legs accordingly, a modified Thurston map
$\tilde g_\theta$ is defined; it is unobstructed and combinatorially
equivalent to $f_c$\,.
\end{description}
See \cite{hss} for a convergence proof in the periodic case, which replaces
Teichm\"uller space with a more explicit spider space. The essential spider
map $\tilde g_\theta$ is constructed in \cite{bfh}, and the relation between
obstructions, kneading sequences, and the satellite case is obtained in
\cite{book2h}.
Note that the description above assumes landing properties of
parameter rays according to \cite{ser, mer}, and the spider
algorithm is just a method to compute parameters numerically. Alternatively,
one may discuss the spider map $g_\theta$ directly and conclude the existence
of quadratic polynomials with specific combinatorics. There are several
variants of implementing the spider algorithm:
\begin{itemize}
\item In a pullback step, each leg and endpoint has two preimages under the
quadratic polynomial, or the preimage is the critical point with two legs. To
choose unique preimages, either employ the cyclic order of rays at $\infty$,
which is related to intersection numbers, 
or consider the angles of the legs at $\infty$.
\item Either normalize the position of two finite marked points, or
assume that all polynomials are of the form $z^2+c_n$\,. This increases the
dimension of Teichm\"uller space by one and gives an additional
eigenvalue $\lambda=1/2$.
\item Each leg is encoded as a sequence of points, such that the curve is
homotopic to a polygonal curve with respect to the marked points. Since the
preimages of straight lines are hyperbolas in general, this means that each
hyperbola segment is replaced with a line segment again; we must check
that it is homotopic. When this condition is violated in the current
step for one or more segments, we may either refine the discretization there
(and prune somewhere else), or restart with an overall finer discretization.
\end{itemize}
In the satellite Misiurewicz case~3, Hubbard--Schleicher \cite{hss} observed
that colliding marked points converge to postcritical points of $f_c$ and the
polynomials converge to $f_c$\,. To understand this process in general,
Selinger \cite{ext, char} considered the extension of the Thurston pullback to
augmented Teichm\"uller space and the dynamics on the canonical stratum. 
This phenomenon motivated the research for the convergence
Theorem~\ref{Tselconv} as well. Intuitively, the points must collide
because the unique obstruction is pinched, and since they stay close together
while moving, the pullback of $g_\theta$ shall be similar to the pullback
defined by $\tilde g_\theta$ or $f_c$\,. But this description involves
interchanging limits, so it is not obvious that the marked points get close to
the expected limit and stay there long enough to be attracted.

\begin{thm}[following Hubbard--Schleicher and Selinger]\label{Thss}
For the pullback defined by the unmodified spider map $g_\theta$\,,
the polynomials converge to $f_c$ and the marked points converge to
postcritical points, with suitable collisions in the satellite Misiurewicz
case~$3$.
\end{thm}

\textbf{Proof:} According to the references given above, either
$g_\theta$ or $\tilde g_\theta$ is unobstructed and equivalent to $f_c$\,
In case~3, the Thurston pullback for $g_\theta$ diverges due to the L\'evy
cycle. The essential map $\tilde g_\theta$ is equivalent to $f_c$ and the other
component maps are homeomorphisms. So Theorem~\ref{Tselconv} applies and
gives convergence immediately. 

Recall the following steps of its proof. In the context of
Proposition~\ref{Pselconv} the current situation was called scenario~2: 
the pullback in configuration space extends to a neighborhood of the
prospective limit. The eigenvalues either come from the modified Thurston
pullback, or they are of the form $\lambda^{rq}=\rho^{-r}$,
$\lambda^q\neq\rho^{-1}$, where $\rho$ is the repelling multiplier of the
$q$-cycle of $f_c$\,. The techniques of Selinger show that the points in
configuration space get arbitrarily close to the prospective limit, such that
a segment of an invariant path in Teichm\"uller space projects into an
attracting neighborhood of that configuration. Then it cannot happen that at
some step another branch of the pullback relation becomes active, so the
points do not jump away. \mybox

In contrast to the situation of formal matings, this generalized convergence
property is not crucial from a numerical perspective, since the modification
from $g_\theta$ to $\tilde g_\theta$ is simple and explicit. As a completely
different approach, the parameter $c$ may be obtained by drawing the parameter
ray $\r_\sM(\theta)$ and starting a Newton iteration from the approximate
endpoint.
Now, let us consider an alternative implementation of the spider algorithm,
which pulls back a path in moduli space instead of spiders in Teichm\"uller
space. So the legs are invisible, but the movement of the feet is recorded:

\begin{init}[Spider algorithm with a path]\label{Isp}
Suppose $\theta=\theta_1\in\Q\setminus\Z$ has preperiod $k$ and period $p$.
Define $(x_1(t),\,\dots,\,x_{k+p}(t))$ for $0\le t\le1$ as
\ban
x_1(t) &=& t\cdot\e{\i2\pi\theta_1} \nonumber \\
x_p(t) &=& (1-t)\cdot\e{\i2\pi\theta_p} \; ,\quad \mbox{if}\quad k=0 \\
x_j(t) &=& \e{\i2\pi\theta_j} \quad\quad\quad\quad,
\quad \mbox{otherwise.} \nonumber
\ean
Pull this path back continuously with $x_i(t+1)=\pm\sqrt{x_{i+1}(t)-x_1(t)}$.
Then it converges to the marked points of $f_c$ with appropriate collisions.
\end{init}

\textbf{Proof:} The argument is similar to that given for captures and
precaptures according to Initialization~\ref{Ihrc}, and for twisted maps
according to Examples~3.1 and~3.7 in \cite{pmate}. We may
initialize the Thurston pullback for $g_\theta=\phi_\theta\circ F$ by
$\psi_0=\phi_\theta^{-1}$ and $\psi_1=\mathrm{Id}$. There is an obvious
deformation $\psi_t$ along one or two rays, which projects to the path defined
in moduli space. By Proposition~\ref{Ppath}, this shows that the pullback of
the path agrees with the projection of the pullback in Teichm\"uller space.
Note that for $k=0$, we have $x_p(t)=0$ only for $t\ge1$. Likewise, for $k=1$
the relation $x_{k+p}(t)=-x_k(t)$ is satisfied for $t\ge1$ only. \mybox

This algorithm gives the same marked points as the spider
algorithm with legs, and it converges unless there are floating-point
cancellations or problems with the discretization: again, the path is
represented by a polygonal curve, and there is an explicit check for homotopy
violations by the simultaneous deformation of hyperbola segments to line
segments; if that happens, refine or restart. Since only a path
of length $|n\le t\le n+1|=1$ needs to be stored instead of full legs,  we may
take a large number of line segments easily, but there is a trade-off: there
will be little need for refinement, because small hyperbola segments are close
to small line segments, but there is a loss of precision by subtracting numbers
that are approximately equal.

For exponential functions with preperiodic singular value, spiders and modified
spiders are constructed in [SZ, LSV], and convergence of unobstructed
pullback maps follows from [HSS]. The alternative implementation with
a path in moduli space is straightforward, but a check for homotopy violations
will be harder. Examples show convergence of colliding marked points
analogously to Theorem~\ref{Thss}. While the local analysis at the prospective
limit is the same, the extension to augmented Teichm\"uller space is unknown
and so the global analysis is incomplete.
\\[5mm]\small
[HSS] J.~Hubbard, D.~Schleicher, M.~Shishikura, Exponential Thurston maps
   and limits of quadratic differentials, J.~AMS~\textbf{22}, 77--117 (2009).
\\[2mm]
[LSV] B.~Laubner, D.~Schleicher, V.~Vicol, A Combinatorial Classification
   of Postsingularly Finite Complex Exponential Maps,
   Discrete cont.~dyn.~systems~\textbf{22}, 2008. 
\\[2mm]
[SZ] D.~Schleicher, J.~Zimmer, Periodic points and dynamic rays of exponential
   maps, Ann.~Acad.~Scient.~Fenn.~\textbf{28}, 327--354  (2003). 
\end{document}